\documentclass{amsart}
\usepackage[utf8]{inputenc}
\usepackage[T1]{fontenc}
\usepackage{lmodern}
\usepackage{amsmath}
\usepackage{amssymb}
\usepackage{mathtools}
\usepackage{latexsym}
\usepackage{amsrefs}
\usepackage{nicefrac}
\usepackage{microtype}
\usepackage{color}
\usepackage{tikz-cd}
\usepackage{MnSymbol}
\usepackage{enumitem} 
\setlist[enumerate,1]{label=\textup{(\arabic*)}}

\usepackage[pdftitle={Non-Archimedean analytic cyclic homology},
pdfauthor={Guillermo Cortiñas, Ralf Meyer and Devarshi Mukherjee},
pdfsubject={Mathematics}]{hyperref}

\BibSpec{book}{%
  +{}  {\PrintPrimary}                {transition}
  +{,} { \textit}                     {title}
  +{.} { }                            {part}
  +{:} { \textit}                     {subtitle}
  +{,} { \PrintEdition}               {edition}
  +{}  { \PrintEditorsB}              {editor}
  +{,} { \PrintTranslatorsC}          {translator}
  +{,} { \PrintContributions}         {contribution}
  +{,} { }                            {series}
  +{,} { \voltext}                    {volume}
  +{,} { }                            {publisher}
  +{,} { }                            {organization}
  +{,} { }                            {address}
  +{,} { \PrintDateB}                 {date}
  +{,} { }                            {status}
  +{}  { \parenthesize}               {language}
  +{}  { \PrintTranslation}           {translation}
  +{;} { \PrintReprint}               {reprint}
  +{.} { }                            {note}
  +{.} {}                             {transition}
  +{} { \PrintDOI}                   {doi}
  +{} { available at \url}            {eprint}
  +{}  {\SentenceSpace \PrintReviews} {review}
}

\renewcommand*{\PrintDOI}[1]{\href{http://dx.doi.org/\detokenize{#1}}{doi: \detokenize{#1}}}

\newcommand{\comment}[1]{}  

\theoremstyle{plain}
\newtheorem{theorem}[equation]{Theorem}
\newtheorem{lemma}[equation]{Lemma}
\newtheorem{corollary}[equation]{Corollary}
\newtheorem{proposition}[equation]{Proposition}
\theoremstyle{remark}
\newtheorem{remark}[equation]{Remark}
\theoremstyle{definition}
\newtheorem{definition}[equation]{Definition}
\newtheorem{example}[equation]{Example}
\numberwithin{equation}{section}

\newcommand\C{\mathbb C}
\newcommand\N{\mathbb N}
\newcommand\Q{\mathbb Q}
\newcommand\R{\mathbb R}
\newcommand\Z{\mathbb Z}

\newcommand\Mat{\mathbb M}

\newcommand\reg{\mathrm{reg}}
\newcommand{\congto}{\xrightarrow{\simeq}}

\newcommand{\bdd}{\mathcal B}
\newcommand{\coma}{\widehat}
\newcommand{\comb}[1]{\overbracket[.7pt][1.4pt]{#1}{}}

\newcommand*{\Fil}{\mathcal F}
\newcommand*{\AU}[1]{#1^+}

\newcommand*{\tub}[2]{\mathcal{U}(#1,#2)}
\newcommand*{\tens}{\mathsf{T}}
\newcommand*{\jens}{\mathsf{J}}

\newcommand*{\tans}{\mathcal{T}}
\newcommand*{\jans}{\mathcal{J}}

\newcommand{\updagger}{\textup{\tiny\!\!\dagger}}

\newcommand{\idealin}{\mathrel{\triangleleft}} 

\newcommand{\diff}{\mathrm{d}}
\newcommand{\Diff}{\mathrm{D}}
\newcommand{\Cont}{\mathrm{C}}
\newcommand{\defeq}{\mathrel{:=}} 

\newcommand*{\into}{\rightarrowtail}
\newcommand*{\onto}{\twoheadrightarrow}

\newcommand*{\ling}[1]{#1_\mathrm{lg}}
\newcommand*{\relling}[2]{#1_{\mathrm{lg}(#2)}}
\newcommand*{\torf}[1]{#1_\mathrm{tf}}

\DeclarePairedDelimiter{\abs}{\lvert}{\rvert}
\DeclarePairedDelimiter{\norm}{\lVert}{\rVert}
\DeclarePairedDelimiter{\floor}{\lfloor}{\rfloor}
\DeclarePairedDelimiter{\ceil}{\lceil}{\rceil}
\DeclarePairedDelimiterX{\setgiven}[2]{\{}{\}}{#1\,{:}\,\mathopen{}#2}
\DeclarePairedDelimiterX{\braket}[2]{\langle}{\rangle}{#1\,{,}\,\mathopen{}#2}

\newcommand{\rig}{\mathrm{rig}}
\newcommand{\MW}{\mathrm{MW}}
\newcommand\hot{\mathbin{\comb{\otimes}}}
\newcommand\haotimes{\mathbin{\coma{\otimes}}}

\DeclareMathOperator{\coker}{coker}

\DeclareMathOperator*{\holim}{holim}

\DeclareMathOperator{\HA}{HA}

\newcommand*{\an}{\mathrm{an}}

\newcommand{\spec}{\mathrm{sp}}

\newcommand\ccirc{\circledcirc}
\newcommand{\Fed}{\odot}

\newcommand{\op}{\mathrm{op}}
\newcommand{\ev}{\mathrm{ev}}
\newcommand{\even}{\mathrm{ev}}
\newcommand{\odd}{\mathrm{odd}}
\newcommand{\id}{\mathrm{id}}
\newcommand{\nb}{\nobreakdash}
\newcommand{\dvr}{V}
\newcommand{\dvgen}{\pi}
\newcommand{\dvf}{F}
\newcommand{\resf}{\mathbb F}

\newcommand{\HAC}{\mathbb{HA}}

\DeclareMathOperator{\car}{char}
\DeclareMathOperator{\Hom}{Hom}
\DeclareMathOperator{\HP}{HP}

\begin{document}
\title{Non-Archimedean analytic cyclic homology}

\author{Guillermo Corti\~nas}
\email{gcorti@dm.uba.ar}
\address{Dep.\ Matem\'atica-IMAS, FCEyN-UBA\\
  Ciudad Universitaria Pab 1\\
  1428 Buenos Aires\\
  Argentina}


\author{Ralf Meyer}
\email{rmeyer2@uni-goettingen.de}

\author{Devarshi Mukherjee}
\email{devarshi.mukherjee@mathematik.uni-goettingen.de}

\address{Mathematisches Institut\\
  Georg-August Universit\"at Göttingen\\
  Bunsenstra\ss{}e 3--5\\
  37073 Göttingen\\
  Germany}

\begin{abstract}
  Let~$\dvr$ be a complete discrete valuation ring with fraction
  field~$\dvf$ of characteristic zero and with residue
  field~\(\resf\).  We introduce analytic cyclic homology of
  complete torsion-free bornological algebras over~\(\dvr\).  We
  prove that it is homotopy invariant, stable, invariant under
  certain nilpotent extensions, and satisfies excision.  We use
  these properties to compute it for tensor products with dagger
  completions of Leavitt path algebras.  If~\(R\) is a smooth
  commutative $\dvr$\nb-algebra of relative dimension~$1$, then we
  identify its analytic cyclic homology with Berthelot's rigid
  cohomology of $R\otimes_{\dvr}\resf$.
\end{abstract}

\thanks{The first named author was supported by a Humboldt Research
  Award and by CONICET, and partially supported by grants UBACyT
  256BA, PICT 2017-1395 and MTM2018-096446-B-C21.}

\maketitle

\section{Introduction}

Analytic cyclic homology of complete bornological algebras over $\R$
and~$\C$ was introduced in~\cite{Meyer:Analytic} as a bivariant
generalisation from Banach to bornological algebras of the entire
cyclic cohomology defined in~\cite{Fillmore-Khalkhali:Entire}.  It
was shown to be stable under tensoring with algebras of nuclear
operators and invariant under differentiable homotopies and under
analytically nilpotent extensions and to satisfy excision with
respect to semi-split extensions~\cite{Meyer:HLHA}.

Let~$\dvr$ be a complete discrete valuation ring whose fraction
field~$\dvf$ has characteristic zero.  Let~\(\dvgen\) be a
uniformiser and let $\resf \defeq \dvr/\dvgen\dvr$ be the residue
field.  In this article, we define and study an analytic cyclic
homology theory for complete, torsion-free bornological
\(\dvr\)\nb-algebras (see Section~\ref{sec:preparations} for the
definitions of these terms).  For example, if~\(R\) is a
torsion-free, finitely generated, commutative $\dvr$\nb-algebra,
then its Monsky--Washnitzer dagger completion~$R^\updagger$ introduced
in~\cite{Monsky-Washnitzer:Formal} is such a complete bornological
algebra (see \cites{Cortinas-Cuntz-Meyer-Tamme:Nonarchimedean,
  Meyer-Mukherjee:Bornological_tf}).

We prove that analytic cyclic homology is invariant under dagger
homotopies and under certain nilpotent extensions, that it is
stable, and that it satisfies excision with respect to semi-split
extensions.  We use these properties to compute the analytic cyclic
homology for dagger completed Leavitt and Cohn path algebras of
countable graphs.  For finite graphs, we also compute the analytic
cyclic homology for tensor products with such algebras.  In
particular, it follows that the analytic cyclic homology of the
completed tensor product of~\(R\) with \(\dvr[t,t^{-1}]^\updagger\)
is isomorphic to the direct sum \(\HA_*(R) \oplus \HA_*(R)[1]\),
where~\(\HA_*\) denotes analytic cyclic homology.  This is a variant
of the fundamental theorem in algebraic K-theory.

We also compute $\HA_*(R^\updagger)$ for a smooth, commutative
\(\dvr\)\nb-algebra~$R$ of relative dimension~$1$.  Namely, it is
isomorphic to the de Rham cohomology of~\(R^\updagger\).  If~\(\resf\)
has finite characteristic, then this agrees with Berthelot's rigid
cohomology of $R\otimes\resf$
(see~\cite{Cortinas-Cuntz-Meyer-Tamme:Nonarchimedean}).  Partial
results that we have for smooth, commutative \(\dvr\)\nb-algebras of
higher dimension have not been included because we have not been
able to prove that analytic and periodic cyclic homology coincide in
this generality.

Monsky--Washnitzer cohomology and Berthelot's rigid cohomology are
defined for varieties in finite characteristic by lifting them to
characteristic zero.  In order to define analogous theories for
noncommutative \(\resf\)\nb-algebras, it is natural to replace de
Rham cohomology by cyclic homology.  Indeed,
in~\cite{Cortinas-Cuntz-Meyer-Tamme:Nonarchimedean}, Berthelot's
rigid chomology for commutative \(\resf\)\nb-algebras is linked to
the \emph{periodic} cyclic homology of suitable dagger completed
commutative \(\dvr\)\nb-algebras.  When we allow noncommutative
algebras, however, then the dagger completion process forces us to
replace periodic cyclic homology by the analytic cyclic homology
that is studied here.

In work in progress, we are going to use the theory defined in this
article in order to define an analytic cyclic homology theory for
algebras over the residue field~\(\resf\).  We want to prove
\(\HA_*(A) \cong \HA_*(R^\updagger)\) whenever~\(R\) is a
torsion-free \(\dvr\)\nb-algebra and \(A\cong R/\dvgen R\) is its
reduction to an \(\resf\)\nb-algebra; the crucial point is that this
should not depend on the choice of~\(R\), and this is where we need
analytic instead of periodic cyclic homology.

Several groups of authors have recently been studying cohomology
theories for varieties in finite characteristic with different
approaches.  We mention, in particular, the work of Petrov and
Vologodsky~\cite{Petrov-Vologodsky:Periodic_topological} that uses
topological cyclic homology.

\bigskip

This paper is organised as follows.  Some notational conventions
used throughout the article are reviewed at the end of this
introduction.

In Section~\ref{sec:preparations}, we start by recalling some basic
notions from bornological analysis and from the Cuntz--Quillen
approach to cyclic homology theories.  In particular, we introduce
dagger completions relative to an ideal
(Section~\ref{sec:relative_dagger}), and review the appropriate
notions of extension of bornological modules, noncommutative
differential forms, tensor algebra, and $X$\nb-complex for
bornological algebras.

Section~\ref{sec:definition_theory} introduces the analytic cyclic
pro-complex $\HAC(R)$ of a complete, torsion-free bornological
algebra~$R$.  It is defined as the $X$\nb-complex of the
scalar extension $\tans R\otimes_\dvr \dvf$ of a certain projective
system~$\tans R$ of complete bornological $\dvr$\nb-algebras
functorially associated to~$R$.  Hence, by definition,
$\HAC(R)= \bigl(\HAC(R)_m\bigr)_{m\ge 1}$ is a pro-supercomplex
(that is, a projective system of $\Z/2$-graded chain complexes) of
complete bornological vector spaces over~$\dvf$.
The analytic cyclic homology of~$R$ is
defined as the homology of the homotopy limit of $\HAC(R)$,
\[
\HA_*(R)\defeq H_*(\holim\HAC(R)).
\]

Section~\ref{sec:aqf} is concerned with analytic nilpotence.
Analytically nilpotent pro-algebras and analytically nilpotent
extensions of algebras and pro-algebras are introduced.  A
pro-algebra~$R$ is called analytically quasi-free if every
semi-split, analytically nilpotent extension of~$R$ splits.  In
particular, the analytic tensor pro-algebra~$\tans R$
(see Definition~\ref{def:tans}) is analytically quasi-free
and is part of a semi-split, analytically nilpotent extension
\[
\jans R\into \tans R\onto R.
\]
We define dagger homotopy of (pro-)algebra homomorphisms using the
dagger completion~$\dvr[t]^\updagger$, and we show that any semi-split
analytically nilpotent extension $N\into E\onto R$ with analytically
quasi-free~$E$ is dagger homotopy equivalent to the extension above.
We use this and the invariance of the $X$\nb-complex under dagger
homotopies
to show that~$\HA$ is invariant under dagger homotopies.
This implies that~$\HA$ is invariant under analytically nilpotent
extensions
and that $\HAC(R)$ is homotopy equivalent to $X(R\otimes \dvf)$
if~\(R\) is analytically quasi-free.

Section~\ref{sec:excision} is devoted to the proof of the Excision
Theorem, which says that if
\[
K \overset{i}\into E \overset{p}\onto Q
\]
is a semi-split pro-algebra extension, then there is a natural exact
triangle
\[
\HAC(K)\overset{i}\to \HAC(E)\overset{p}\to
\HAC(Q)\overset{\delta}\to \HAC(K)[-1].
\]
Applying $\holim$ and taking homology, this gives a natural
$6$\nb-term exact sequence
\[
  \begin{tikzcd}
    \HA_0(K) \arrow[r, "i_*"] &
    \HA_0(E) \arrow[r, "p_*"] &
    \HA_0(Q) \arrow[d, "\delta"] \\
    \HA_1(Q) \arrow[u, "\delta"] &
    \HA_1(E) \arrow[l, "p_*"] &
    \HA_1(K).  \arrow[l, "i_*"]
  \end{tikzcd}
\]
The proof of the excision theorem follows the structure of its
archimedean version in \cites{Meyer:HLHA, Meyer:Excision}, and
adapts it to the present case.

The stability of~$\HA$ under matricial embeddings is proved in
Section~\ref{sec:stability}.  Any pair $X,Y$ of torsion-free
bornological \(\dvr\)\nb-modules with a surjective bounded linear
map $\braket{\cdot}{\cdot}\colon Y\otimes X\to \dvr$ gives rise to
an algebra $\mathcal{M}(X,Y)$ with underlying bornological
\(\dvr\)\nb-module $X\otimes Y$.  We show in
Proposition~\ref{prop:stable} that~$\HA$ is invariant under
tensoring with the dagger completion~$\mathcal{M}(X,Y)^\updagger$.  For
example, the algebra of finite matrices~$\Mat_n$ with $n\le \infty$
and the algebra of matrices with entries going to zero at infinity
are of the form $\mathcal{M}(X,Y)^\updagger$ for suitable $X$ and~$Y$.
Thus~$\HA$ is invariant under tensoring with such algebras.

Section~\ref{sec:Leavitt} is concerned with Leavitt path algebras.
For a directed graph~$E$ with finitely many vertices and a complete
bornological algebra~$R$, Theorems \ref{the:uncompleted_Leavitt}
and~\ref{the:completed_Leavitt} compute $\HA(R\hot L(E)^\updagger)$ in
terms of $\HA(R)$ and a matrix~$N_E$ related to the incidence matrix
of~$E$:
\[
  \HAC(R\hot L(E)^\updagger)
  \simeq (\coker(N_E)\oplus\ker(N_E)[1]) \otimes \HAC(R).
\]
For trivial~\(R\), the homotopy equivalence
\(\HAC(L(E)^\updagger) \simeq (\coker(N_E)\oplus\ker(N_E)[1])\) is
shown also for graphs with countably many vertices.  If~\(E\) is the
graph with one vertex and one loop, it follows that~$\HA$ satisfies
a version of Bass' fundamental theorem:
\[
\HAC(R\hot\dvr[t,t^{-1}]^\updagger)\simeq \HAC(R)\oplus\HAC(R)[-1].
\]
We also compute \(\HAC(R\hot C(E)^\updagger)\) for the Cohn path
algebra if~\(E\) has finitely many vertices, and
\(\HAC(C(E)^\updagger)\) if~\(E\) has countably many vertices.

In Section~\ref{sec:filtered_Noether} we show that if~$R$ is smooth
commutative of relative dimension one, then the analytic cyclic
homology of its dagger completion is the same as the rigid
cohomology of its reduction modulo~$\dvgen$ (see
Theorem~\ref{the:non-singular_1-dim}).  That is,
\[
  \HA_n(R^\updagger) \cong H_\rig^n(R/\dvgen R)
\]
for $n=0,1$.  We outline the idea of the proof.
By~\cite{Cortinas-Cuntz-Meyer-Tamme:Nonarchimedean},
\(H_\rig^n(R/\dvgen R)\) is isomorphic to the periodic cyclic
homology of $R^\updagger\otimes \dvf$.  By
Corollary~\ref{cor:analytically_quasi_free_computation}, $\HA$ and
$\HP(\cdot\otimes \dvf)$ agree on analytically quasi-free
bornological \(\dvr\)\nb-algebras.  It is well known that a smooth
algebra~$R$ of relative dimension~$1$ is quasi-free in the sense
that any square-zero extension of is splits or, equivalently, that
the bimodule $\Omega^1(R)$ of noncommutative differential
$1$\nb-forms admits a connection.  We show in
Theorem~\ref{thm:qfdagger} that if~$R$ is a torsion-free, complete
bornological algebra and~$\nabla$ is a connection on~$\Omega^1(R)$
that satisfies an extra condition, then~$R^\updagger$ is analytically
quasi-free.  We prove that a smooth commutative algebra with
the fine bornology admits such a connection (see
Lemma~\ref{filtered_Noetherian_finite_connection}).

\subsubsection*{Acknowledgement.}

A first draft of this manuscript was completed during a 6-month stay
of the first author at the Universit\"at M\"unster, supported by a
Humboldt Research Award.  He thanks his host, Joachim Cuntz and all
the members of the Operator Algebra group, particularly its head,
Wilhelm Winter, for their wonderful hospitality.  He also thanks the
Humboldt Foundation for its support.

\subsection{Some notation}
\label{sec:conventions}

Throughout this article, we shall use the following notation.
Let~\(\N^*\) be the set of nonzero natural numbers.  Let~\(\dvr\) be
a complete discrete valuation ring, \(\dvgen\in \dvr\) a
uniformiser, \(\resf\) the residue field \(\dvr/(\dvgen)\)
of~\(\dvr\), and~\(\dvf\) the fraction field of~\(\dvr\).  While our
definitions work in complete generality, our homotopy invariance,
stability and excision theorems only work if~\(\dvf\) has
characteristic zero.  All tensor products~\(\otimes\) are taken
over~\(\dvr\).  By convention, algebras are allowed to be non-unital
throughout this article.  An \emph{ideal} in a possibly non-unital
\(\dvr\)\nb-algebra means a two-sided ideal that is also a
\(\dvr\)\nb-submodule.

\section{Preparations}
\label{sec:preparations}

\subsection{Bornologies}
\label{sec:bornologies}
\numberwithin{equation}{subsection}

As in~\cite{Cortinas-Cuntz-Meyer-Tamme:Nonarchimedean}, bornological
\(\dvr\)\nb-algebras play a crucial role.  We first recall some
basic terminology about bornologies from
\cites{Cortinas-Cuntz-Meyer-Tamme:Nonarchimedean,
  Meyer-Mukherjee:Bornological_tf}.

\begin{definition}
  A \emph{bornology} on a set~\(S\) is a set~\(\bdd\) of subsets,
  called \emph{bounded subsets}, such that finite unions and subsets
  of bounded subsets are bounded and finite subsets are bounded.  A
  \emph{bornological set} is a set with a bornology.
\end{definition}

\begin{definition}
  A map \(f\colon S_1 \to S_2\) between bornological sets is
  \emph{bounded} if it maps bounded subsets to bounded subsets.  It
  is a \emph{bornological embedding} if it is injective and \(T\subseteq S_1\) is
  bounded if and only if \(f(T)\subseteq S_2\) is bounded.  It is a
  \emph{bornological quotient map} if it is bounded and any bounded
  subset \(T\subseteq S_2\) is the image of a bounded subset
  of~\(S_1\).
\end{definition}

\begin{definition}
  A \emph{bornological \(\dvr\)\nb-module} is a
  \(\dvr\)\nb-module~\(R\) with a bornology such that any bounded
  subset is contained in a bounded \(\dvr\)\nb-submodule or,
  equivalently, the \(\dvr\)\nb-submodule generated by a bounded
  subset is again bounded.  A \emph{bornological
    \(\dvr\)\nb-algebra} is a bornological \(\dvr\)\nb-module~\(R\)
  with a multiplication \(R\times R \to R\) that is bounded in the
  sense that \(S\cdot T\) is bounded if \(S,T \subseteq R\) are
  bounded.
\end{definition}

\begin{definition}
  A bornological \(\dvr\)\nb-module is \emph{complete} if any
  bounded subset is contained in a bounded \(\dvr\)\nb-submodule
  that is \(\dvgen\)\nb-adically complete.  The
  \emph{completion}~\(\comb{M}\) of a bornological
  \(\dvr\)\nb-module~\(M\) is a complete bornological
  \(\dvr\)\nb-module with a bounded map \(M \to \comb{M}\) that is
  universal in the sense that any bounded map from~\(M\) to a
  complete bornological \(\dvr\)\nb-module factors uniquely through
  it (see
  \cite{Cortinas-Cuntz-Meyer-Tamme:Nonarchimedean}*{Definition~2.14}).
\end{definition}

\begin{example}
  \label{exa:fine}
  Let~\(M\) be a \(\dvr\)\nb-module.  The \emph{fine bornology}
  on~\(M\) consists of those subsets of~\(M\) that are contained in
  a finitely generated \(\dvr\)\nb-submodule.  It is the smallest
  \(\dvr\)\nb-module bornology on~\(M\).  It is the only bornology
  on~\(M\) if~\(M\) itself is finitely generated.  We equip the
  fraction field~\(\dvf\) with the fine bornology.  If~\(R\) is a
  \(\dvr\)\nb-algebra, then the fine bornology makes it a
  bornological \(\dvr\)\nb-algebra.  The fine bornology is
  automatically complete.
\end{example}

\begin{definition}
  Let \(M_1\) and~\(M_2\) be bornological \(\dvr\)\nb-modules.  The
  \emph{tensor product bornology} on the \(\dvr\)\nb-module
  \(M_1 \otimes M_2\) consists of all subsets that are contained in
  \(S_1 \otimes S_2\) for bounded bornological
  \(\dvr\)\nb-submodules \(S_j \subseteq M_j\) for \(j=1,2\).  The
  \emph{complete bornological tensor product} \(M_1 \hot M_2\) is
  defined as the bornological completion of \(M_1 \otimes M_2\) with
  the tensor product bornology.
\end{definition}

The universal property of tensor products easily implies the
following:

\begin{proposition}
  Let \(M_1\), \(M_2\) and~\(N\) be bornological
  \(\dvr\)\nb-modules.  Bounded \(\dvr\)\nb-linear maps
  \(M_1 \otimes M_2 \to N\) are in natural bijection with bounded
  \(\dvr\)\nb-bilinear maps \(M_1 \times M_2 \to N\).
\end{proposition}

\begin{corollary}
  Let \(M_1\), \(M_2\) and~\(N\) be complete bornological
  \(\dvr\)\nb-modules.  Bounded \(\dvr\)\nb-linear maps
  \(M_1 \hot M_2 \to N\) are in natural bijection with bounded
  \(\dvr\)\nb-bilinear maps \(M_1 \times M_2 \to N\).
\end{corollary}

\begin{example}
  \label{exa:fine_tensor}
  Continuing Example~\ref{exa:fine}, let \(M_1\) be a
  \(\dvr\)\nb-module with the fine bornology and let~\(M_2\) be a
  complete bornological \(\dvr\)\nb-module.  Then the tensor product
  bornology on \(M_1 \otimes M_2\) is already complete because the
  tensor product of a \(\dvgen\)\nb-adically complete
  \(\dvr\)\nb-module with a finitely generated \(\dvr\)\nb-module is
  complete.
  Thus
  \(M_1 \hot M_2 = M_1 \otimes M_2\) in this case.  This applies, in
  particular, if \(M_1=\dvf\).  If both \(M_1\) and~\(M_2\) carry
  the fine bornology, then the tensor product bornology on
  \(M_1 \hot M_2 = M_1 \otimes M_2\) is the fine bornology as well.
\end{example}

\begin{definition}[\cite{Meyer-Mukherjee:Bornological_tf}*{Definition~4.1}]
  A bornological \(\dvr\)\nb-module~\(M\) is
  \emph{\textup{(}bornologically\textup{)} torsion-free} if
  multiplication by~\(\dvgen\) is a bornological embedding.
\end{definition}

\begin{remark}
  Let~\(M\) be a bornological \(\dvr\)\nb-module.  If
  \(S\subseteq M\), then define
  \[
    \dvgen^{-1} S \defeq \setgiven{x\in M}{\dvgen\cdot x\in S}.
  \]
  This depends on~\(M\) and not just on~\(S\).  By definition, \(M\)
  is torsion-free if and only if multiplication by~\(\dvgen\) is
  injective and~\(\dvgen^{-1} S\) is bounded for all bounded subsets
  \(S\subseteq M\).
\end{remark}

\begin{proposition}[\cite{Meyer-Mukherjee:Bornological_tf}*{Proposition~4.3}]
  \label{pro:torf_embed_to_fraction_field}
  A bornological \(\dvr\)\nb-module~\(M\) is torsion-free if and
  only if the canonical map \(M \to M \otimes \dvf\) is a
  bornological embedding.
\end{proposition}

\begin{example}
  A \(\dvr\)\nb-module~\(M\) with the fine bornology is torsion-free
  if and only if~\(M\) is torsion-free in the usual sense.
\end{example}

\begin{definition}
  Let~\(M\) be any bornological \(\dvr\)\nb-module and define
  \(\torf{M} \subseteq M\otimes \dvf\) as the image of the canonical
  map \(M \to M\otimes \dvf\), equipped with the restriction of the
  bornology of \(M\otimes \dvf\).
\end{definition}

\begin{proposition}[\cite{Meyer-Mukherjee:Bornological_tf}*{Proposition~4.4}]
  The canonical map \(M\to \torf{M}\) is the universal map
  from~\(M\) to a torsion-free bornological \(\dvr\)\nb-module.
\end{proposition}

\begin{definition}
  \label{def:semi-dagger}
  A bornological \(\dvr\)\nb-algebra~\(R\) is \emph{semi-dagger} if
  any bounded subset \(S\subseteq R\) is contained in a bounded
  \(\dvr\)\nb-submodule \(T\subseteq R\) with
  \(\dvgen\cdot T\cdot T \subseteq T\) (see
  \cite{Meyer-Mukherjee:Bornological_tf}*{Proposition~3.4}).
  Let~\(R\) with the bornology~\(\bdd\) be a bornological
  \(\dvr\)\nb-algebra.  There is a smallest semi-dagger bornology
  on~\(R\) that contains~\(\bdd\).  It is denoted~\(\ling{\bdd}\)
  and called the \emph{linear growth bornology} on~\(R\); we
  write~\(\ling{R}\) for~\(R\) with the linear growth bornology (see
  \cite{Meyer-Mukherjee:Bornological_tf}*{Definition~3.5 and
    Lemma~3.6}).
\end{definition}

\begin{definition}
  A \emph{dagger algebra} is a bornological \(\dvr\)\nb-algebra that
  is complete, (bornologically) torsion-free, and semi-dagger.  The
  \emph{dagger completion} of a bornological
  \(\dvr\)\nb-algebra~\(R\) is a dagger algebra~\(R^\updagger\) with
  a bounded \(\dvr\)\nb-algebra homomorphism \(R\to R^\updagger\)
  that is universal in the sense that any bounded homomorphism
  from~\(R\) to a dagger algebra factors uniquely through it.
\end{definition}

\begin{theorem}[\cite{Meyer-Mukherjee:Bornological_tf}*{Theorem~5.3}]
  \label{the:dagger_completion}
  If~\(R\) is already torsion-free, then~\(R^\updagger\) is the
  completion of~\(\ling{R}\).  In general, it is the completion of
  \(\ling{(\torf{R})}\).
\end{theorem}

\begin{example}
  \label{exa:fine_dagger}
  The dagger completion~\(R^\updagger\) of a torsion-free, finitely
  generated, commutative \(\dvr\)\nb-algebra is usually defined as
  the weak completion of~\(R\) by Monsky and
  Washnitzer~\cite{Monsky-Washnitzer:Formal}.  This agrees with our
  definition of~\(R^\updagger\) by
  \cite{Cortinas-Cuntz-Meyer-Tamme:Nonarchimedean}*{Theorem~3.2.1}:
  the dagger completion of the fine bornology on~\(R\) is naturally
  isomorphic to the weak completion of~\(R\), equipped with a
  canonical bornology.
\end{example}

\begin{proposition}[\cite{Cortinas-Cuntz-Meyer-Tamme:Nonarchimedean}*{Proposition~3.1.25}]
  \label{pro:tensor_commutes_dagger}
  Let \(A\) and~\(B\) be torsion-free, complete bornological
  algebras.  Then
  \(\ling{(A\otimes B)} \cong \ling{A} \otimes \ling{B}\) and
  \((A\otimes B)^\updagger \cong A^\updagger \hot B^\updagger\).
\end{proposition}

\begin{corollary}
  \label{cor:tensor_dagger}
  A completed tensor product of two dagger algebras is again a
  dagger algebra.
\end{corollary}

\begin{proof}
  A completed tensor product is complete by definition.  It remains
  semi-dagger by Proposition~\ref{pro:tensor_commutes_dagger}, and
  torsion-free by
  \cite{Meyer-Mukherjee:Bornological_tf}*{Proposition~4.12}.
\end{proof}

\subsection{Relative dagger completions}
\label{sec:relative_dagger}

We shall define analytic cyclic homology for torsion-free,
complete bornological \(\dvr\)\nb-algebras~\(R\) that need not be
dagger algebras.  This uses a variant of the linear growth bornology
relative to an ideal.

Let~\(R\) be a \(\dvr\)\nb-algebra and let \(M\) and~\(N\) be
\(\dvr\)\nb-submodules of~\(R\).  Let \(M N \subseteq R\) be the
\(\dvr\)\nb-submodule generated by all products~\(x y\) with
\(x\in M\) and \(y \in N\).  Let
\begin{equation}
  \label{eq:diamond}
  M^\diamond \defeq \sum_{i=0}^\infty \dvgen^i M^{i+1},\qquad
  M^{(n)} \defeq \sum_{i=1}^n M^i.
\end{equation}
A subset of~\(R\) has linear growth if and only if it is contained
in~\(M^\diamond\) for some bounded \(\dvr\)\nb-submodule~\(M\)
of~\(R\) (with the present definitions, this is
\cite{Meyer-Mukherjee:Bornological_tf}*{Lemma~3.6}).

\begin{lemma}
  \label{lem:diam}
  Let~\(R\) be a \(\dvr\)\nb-algebra and let \(M,N \subseteq R\) be
  \(\dvr\)-submodules.  Then
  \begin{enumerate}
  \item \label{lem:diam_1}%
    \(M^\diamond + N^\diamond \subseteq (M+N)^\diamond\);
  \item \label{lem:diam_2}%
    \(M \cdot N^\diamond \subseteq ((M \cdot N + N)^{(2)})^\diamond\) and
    \(N^\diamond \cdot M \subseteq ((N\cdot M + N)^{(2)})^\diamond\);
  \item \label{lem:diam_3}%
    \(\dvgen\cdot M^\diamond\cdot M^\diamond \subseteq M^\diamond\);
  \item \label{lem:diam_4}%
    \(M^\diamond \cdot N^\diamond \subseteq
    ((M+N)^{(2)})^\diamond\);
  \item \label{lem:diam_5}%
    \((M^\diamond)^\diamond = M^\diamond\).
  \end{enumerate}
\end{lemma}

\begin{proof}
  The definition of~\(M^\diamond\) immediately
  implies~\ref{lem:diam_1}.  The following computation shows the
  first assertion of~\ref{lem:diam_2}:
  \begin{align*}
    M\cdot N^\diamond
    &=\sum_{i\ge 1}\dvgen^{2i-1} M N^{2i} + \sum_{i\ge 0} \dvgen^{2i} M N^{2i+1}
    \\&= \sum_{i\ge 1} \dvgen^{2i-1}(M N)N^{2i-1}
    + \sum_{i\ge 0}\dvgen^{2i}(MN)(N^2)^i
    \\&\subseteq (M N+N)^\diamond+(M N+N^2)^\diamond
      \subseteq ((M N+N)^{(2)})^\diamond.
  \end{align*}
  Similar calculations give the second assertion of~\ref{lem:diam_2}
  and~\ref{lem:diam_4}.  Statement~\ref{lem:diam_3} follows because
  \(\dvgen\cdot \dvgen^i M^{i+1} \cdot \dvgen^j M^{j+1} =
  \dvgen^{i+j+1} M^{i+j+1+1}\) for all \(i,j\in\N\).  Then
  \(\dvgen^i\cdot (M^\diamond)^{i+1} \subseteq M^\diamond\) follows
  by induction on~\(i\).  This implies~\ref{lem:diam_5}.
\end{proof}

\begin{definition}
  \label{def:relling}
  Let~\(R\) be a bornological \(\dvr\)\nb-algebra and
  \(I \idealin R\) an ideal.  Let \(\relling{\bdd}{I}\) be the set
  of all subsets of~\(R\) that are contained in \(M + N^\diamond\)
  for bounded \(\dvr\)\nb-submodules \(M \subseteq R\) and
  \(N \subseteq I\).  This is a bornology on~\(R\), called the
  \emph{linear growth bornology relative to~\(I\)}.  Let
  \(\relling{R}{I}\) be~\(R\) with this bornology.
\end{definition}

\begin{example}
  By definition, \(\relling{\bdd}{0} = \bdd\) and
  \(\relling{\bdd}{R}\) is the usual linear growth bornology
  on~\(R\).  So \(\relling{R}{0} = R\) and
  \(\relling{R}{R} = \ling{R}\).
\end{example}

\begin{lemma}
  \label{lem:relling}
  The bornology \(\relling{\bdd}{I}\) is an algebra bornology, and
  its restriction to~\(I\) is semi-dagger.  Let~\(S\) be a
  bornological \(\dvr\)\nb-algebra.  A homomorphism
  \(f\colon R\to S\) is bounded for the bornology
  \(\relling{\bdd}{I}\) if and only if \(f(N)\) has linear growth
  in~\(S\) for all bounded subsets \(N\subseteq I\) and \(f(M)\) is
  bounded in~\(S\) for all bounded subsets \(M\subseteq R\).
\end{lemma}

\begin{proof}
  Since~\(I\) is an ideal, Lemma~\ref{lem:diam} implies that
  \(\relling{\bdd}{I}\) makes~\(R\) a bornological
  \(\dvr\)\nb-algebra.  And a subset of~\(I\) belongs to
  \(\relling{\bdd}{I}\) if and only if it is contained
  in~\(N^\diamond\) for some bounded \(\dvr\)\nb-submodule
  \(N\subseteq I\).  The restriction of~\(\relling{\bdd}{I}\)
  to~\(I\) is semi-dagger by Lemma~\ref{lem:diam}.  If \(M\)
  and~\(N\) are as in Definition~\ref{def:relling}, then
  \(f(M+N^\diamond) = f(M) + f(N)^\diamond\).  This is bounded
  in~\(S\) if and only if \(f(M)\) is bounded and \(f(N)\) has
  linear growth.
\end{proof}



\begin{lemma}
  \label{lem:rel=abs}
  Let~\(R\) be a bornological algebra and let \(I\) and~\(J\) be
  ideals in \(R\) with \(I \subseteq J\) and
  \(R/I = \relling{(R/I)}{J/I}\).  Then
  \(\relling{R}{J} = \relling{R}{I}\).  In particular, if \(R/I\) is
  semi-dagger, then \(\relling{R}{I} = \ling{R}\).
\end{lemma}

\begin{proof}
  By Lemma~\ref{lem:relling},
  the bornology \(\relling{\bdd}{J}\) on~\(R\) is the smallest one
  that contains the given bornology and makes~\(J\) semi-dagger, and
  similarly for~\(I\).  And the assumption
  \(R/I = \relling{(R/I)}{J/I}\) says that \(J/I\subseteq R/I\) is
  semi-dagger in the quotient bornology on~\(R/I\).  This is the
  same as the quotient bornology induced by~\(\relling{\bdd}{I}\).
  \cite{Meyer-Mukherjee:Bornological_tf}*{Theorem~3.7} says that an
  extension of semi-dagger algebras remains semi-dagger.  This
  theorem applied to the extension \(I \into J \onto J/I\), equipped
  with the restrictions of the bornology \(\relling{\bdd}{I}\) on
  \(I\) and~\(J\) and the resulting quotient bornology on~\(J/I\)
  shows that~\(J\) is semi-dagger also in the bornology
  \(\relling{\bdd}{I}\).  Then
  \(\relling{\bdd}{J}\subseteq \relling{\bdd}{I}\).  And
  \(\relling{\bdd}{I}\subseteq \relling{\bdd}{J}\) is trivial.
\end{proof}

\begin{lemma}
  \label{lem:relgtf}
  Let~\(R\) be a bornological algebra and \(I \idealin R\) an ideal.
  If~\(R\) is torsion-free, then so is \(\relling{R}{I}\).
\end{lemma}

\begin{proof}
  Let \(S \subseteq \dvgen R\) be a bounded subset in
  \(\relling{R}{I}\).  By definition, there are bounded submodules
  \(M \subseteq R\) and \(N\subseteq I\) with
  \(S \subseteq M + N^\diamond\).  And
  \[
    M + N^\diamond
    = M + N + \sum_{i=1}^\infty \dvgen^i N^{i+1}
    = M + N + \dvgen\cdot
    \biggl(\sum_{i=0}^\infty \dvgen^i N^{i+2}\biggr).
  \]
  Since \(\dvgen^i N^{i+2} \subseteq \dvgen^i (N^{(2)})^{i+1}\) for
  all \(i\ge0\), the subset \(\sum_{i=0}^\infty \dvgen^i N^{i+2}\)
  belongs to \(\relling{\bdd}{I}\).  Since \(M + N\) is bounded
  in~\(R\) and~\(R\) is torsion-free,
  \(\dvgen^{-1}\cdot (M+N)\) is bounded.  Then
  \[
    \dvgen^{-1} S
    \subseteq \dvgen^{-1}(M+N^\diamond)
    \subseteq \dvgen^{-1}(M+N) + \sum_{i=0}^\infty \dvgen^i N^{i+2}
    \in \relling{\bdd}{I}.\qedhere
  \]
\end{proof}

\begin{definition}
  Let~\(R\) be a torsion-free bornological algebra and
  \(I \idealin R\) an ideal.  The \emph{dagger completion of~\(R\)
    relative to~\(I\)} is the completion
  \((R,I)^\updagger \defeq \comb{\relling{R}{I}}\).
\end{definition}

We shall never apply (relative) dagger completions when~\(R\) is not
already bornologically torsion-free.  In general, the correct
definition of the relative dagger completion of \((R,I)\) would be
\((\torf{R},\torf{I})^\updagger\), where~\(\torf{I}\) is identified
with its image in~\(\torf{R}\) (compare
Theorem~\ref{the:dagger_completion}).

\begin{proposition}
  \label{pro:relative_dagger}
  Let \(R\) and~\(S\) be torsion-free bornological
  \(\dvr\)\nb-algebras, \(I\subseteq R\) an ideal, and
  \(f\colon R \to S\) a bounded algebra homomorphism.  Assume~\(S\)
  to be complete.  There is a bounded algebra homomorphism
  \((R,I)^\updagger \to S\) extending~\(f\), necessarily unique, if
  and only if \(f(M)\) has linear growth for each bounded
  \(\dvr\)\nb-submodule~\(M\) of~\(I\).
\end{proposition}

\begin{proof}
  Use Lemma~\ref{lem:relling} and the universal property of the
  completion.
\end{proof}

There seems to be no analogue of
Proposition~\ref{pro:tensor_commutes_dagger} for relative dagger
completions.

\subsection{Extensions of bornological modules}

An \emph{extension} of \(\dvr\)\nb-modules is a diagram of
\(\dvr\)\nb-modules
\[
  K \overset{i}\into E \overset{p}\onto Q,
\]
that is algebraically exact and such that~\(i\) is a bornological
embedding and~\(p\) is a bornological quotient map.  Equivalently,
\(i\)~is a kernel of~\(p\) and \(p\)~is a cokernel of~\(i\) in the
additive category of bornological \(\dvr\)\nb-modules.  This
following elementary lemma says that this category is
\emph{quasi-abelian}:

\begin{lemma}
  \label{lem:pull-back_push-out_bornological_extension}
  Let \(K\overset{i}\into E \overset{p}\onto Q\) be an extension of
  bornological \(\dvr\)\nb-modules.  Let \(K \overset{f}\to K'\) and
  \(Q''\overset{g}\to Q\) be bounded \(\dvr\)\nb-module maps.  The
  pushout of \(i,f\) and the pullback of \(p,g\) exist and are part
  of morphisms of extensions
  \[
    \begin{tikzcd}
      K \arrow[r,rightarrowtail, "i"] \arrow[d,"f"] &
      E \arrow[r, twoheadrightarrow, "p"] \arrow[d, "\hat{f}"] &
      Q \arrow[d, equal] \\
      K' \arrow[r,rightarrowtail,"i'"] &
      E' \arrow[r, twoheadrightarrow, "p'"] &
      Q,
    \end{tikzcd} \qquad
    \begin{tikzcd}
      K \arrow[r,rightarrowtail, "i''"] \arrow[d, equal] &
      E'' \arrow[r, twoheadrightarrow, "p''"] \arrow[d, "\hat{g}"] &
      Q'' \arrow[d, "g"] \\
      K \arrow[r,rightarrowtail, "i"] &
      E \arrow[r, twoheadrightarrow, "p"] &
      Q.
    \end{tikzcd}
  \]
  Here
  \[
    E'\defeq \frac{K' \oplus E}{\setgiven{(f(k),-i(k))}{k\in K}}, \qquad
    E'' \defeq \setgiven{(e,q'') \in E \times Q''}{p(e)=g(q'')},
  \]
  equipped with the quotient and the subspace bornology,
  respectively, and \(\hat{f}(e) = [(0,e)]\), \(i'(k') = [(k',0)]\),
  \(p'[(k',e)]=p(e)\), \(\hat{g}(e,q'')=e\), \(p''(e,q'') = q''\),
  and \(i''(k)=(i(k),0)\) for \(e \in E\), \(k'\in K'\),
  \(q'' \in Q''\), \(k \in K\).
\end{lemma}

The following proposition is an analogue of Lemma~\ref{lem:rel=abs}
for completions, describing a situation when a partial completion
relative to a submodule is equal to the completion.

\begin{proposition}
  \label{completion:pushout}
  Let \(K\overset{i}\into E\overset{p}\onto Q\) be an extension of
  bornological \(\dvr\)\nb-modules.  Assume~\(Q\) to be complete and
  bornologically torsion-free.  Form the pushout diagram
  \[
    \begin{tikzcd}
      K \arrow[r,rightarrowtail, "i"] \arrow[d, "\mathrm{can}_K"] &
      E \arrow[r, twoheadrightarrow, "p"] \arrow[d, "\gamma"] &
      Q \arrow[d, equal] \\
      \comb{K} \arrow[r,rightarrowtail, "i'"] &
      E' \arrow[r, twoheadrightarrow, "p'"] &
      Q.
    \end{tikzcd}
  \]
  There is a unique isomorphism
  \(\varphi\colon E' \congto \comb{E}\) such that
  \(\varphi\circ \gamma\) is the canonical map \(E \to \comb{E}\).
\end{proposition}

\begin{proof}
  The bottom row is an extension by Lemma
  \ref{lem:pull-back_push-out_bornological_extension}.  Then~\(E'\)
  is complete by
  \cite{Meyer-Mukherjee:Bornological_tf}*{Theorem~2.3}.  The maps
  \(\mathrm{can}_E\colon E \to \comb{E}\) and
  \(\comb{i}\colon \comb{K} \to \comb{E}\) induce a bounded
  \(\dvr\)\nb-module map \(\varphi\colon E' \to \comb{E}\) by the
  universal property of pushouts.  Since~\(E'\) is complete, the
  universal property of~\(\comb{E}\) gives a unique map
  \(\psi \colon \comb{E} \to E'\) with
  \(\psi \circ \mathrm{can}_E = \gamma\).  Then
  \(\varphi \circ \psi \circ \mathrm{can}_E = \varphi \circ \gamma =
  \mathrm{can}_E\).  This implies
  \(\varphi \circ \psi = \id_{\comb{E}}\).  Next,
  \(\psi\circ \comb{i}\circ \mathrm{can}_K = \gamma \circ i = i'
  \circ \mathrm{can}_K\) implies \(\psi \circ \comb{i} = i'\), and
  then \(\psi\circ \varphi \circ i' = \psi \circ \comb{i} = i'\) and
  \(\psi \circ \varphi \circ \gamma = \psi \circ \mathrm{can}_E =
  \gamma\) imply \(\psi \circ \varphi = \id_{E'}\).  So~\(\varphi\)
  is an isomorphism.
\end{proof}

\subsection{Injective maps between completions}
\label{sec:injective_completions}

Unlike in the Archimedean case, all Banach spaces over~\(\dvf\) have
a simple structure.  This implies that they all satisfy a variant of
Grothendieck's Approximation Property.  This is
Proposition~\ref{pro:tensor_injective}, and it will be useful to
describe completions of tensor products.

\begin{definition}
  \label{def:Cont0}
  Let~\(D\) be a set.  Let \(\Cont_0(D,\dvr)\) be the set of all
  functions \(f\colon D \to \dvr\) such that for each \(\delta>0\)
  there is a finite subset \(S\subseteq D\) with
  \(\abs{f(x)}<\delta\) for all \(x\in D\setminus S\).  Define
  \(\Cont_0(D,\dvf)\) similarly.  Equip both with the supremum norm.
\end{definition}

\begin{theorem}
  \label{the:structure}
  Let~\(W\) be a complete, torsion-free bornological
  \(\dvr\)\nb-\hspace{0pt}module.  Any \(\dvgen\)\nb-adically
  complete bounded \(\dvr\)\nb-submodule~\(M\) of~\(W\) is
  isomorphic to \(\Cont_0(D,\dvr)\) for some set~\(D\).
\end{theorem}

\begin{proof}
  The map \(W \to W\otimes \dvf\) is a bornological embedding by
  Proposition~\ref{pro:torf_embed_to_fraction_field}.  The subset
  \(\dvf \cdot M \subseteq W\otimes \dvf\) is an \(\dvf\)\nb-vector
  subspace.  Define the \emph{gauge norm} on \(\dvf\cdot M\) by
  \[
    \norm{x} \defeq
    \inf \setgiven{\abs{\dvgen}^j}{\dvgen^{-j}\cdot x \in M}.
  \]
  It is a non-Archimedean norm and makes \(\dvf\cdot M\) a Banach
  \(\dvf\)\nb-vector space with unit ball~\(M\).  It takes values in
  \(\setgiven{\abs{\dvgen}^n}{n\in\Z}\cup \{0\}\) by construction.
  Hence there is a set~\(D\) and an isometric isomorphism
  \(\dvf M \cong \Cont_0(D,\dvf)\) (see
  \cite{Schneider:Nonarchimedean}*{Remark~10.2}).  It maps~\(M\)
  isomorphically onto the the unit ball of \(\Cont_0(D,\dvf)\),
  which is \(\Cont_0(D,\dvr)\).
\end{proof}

\begin{corollary}
  \label{cor:nice_inductive_system}
  Let~\(W\) be a complete, torsion-free bornological
  \(\dvr\)\nb-\hspace{0pt}module.  Then~\(W\) is isomorphic to the
  colimit of an inductive system of complete \(\dvr\)\nb-modules of
  the form \((\Cont_0(D_n,\dvr), f_{n,m})_{n,m\in S}\) with a
  directed set \((S,\le)\), sets~\(D_n\) for \(n\in S\), and
  injective, bounded \(\dvr\)\nb-linear maps
  \(f_{n,m}\colon \Cont_0(D_m,\dvr) \hookrightarrow
  \Cont_0(D_n,\dvr)\) for \(n,m\in S\), \(n \ge m\).
\end{corollary}

\begin{proof}
  The complete \(\dvr\)\nb-submodules of~\(W\) form a directed set
  under inclusion, and this defines an inductive system with
  injective structure maps and with colimit~\(W\) by
  \cite{Cortinas-Cuntz-Meyer-Tamme:Nonarchimedean}*{Proposition~2.10}.
  Theorem~\ref{the:structure} identifies the entries in this
  inductive system with \(\Cont_0(D,\dvr)\) for suitable sets~\(D\).
\end{proof}

\begin{lemma}
  \label{lem:tensor_injective}
  Let
  \(f\colon \Cont_0(D_1,\dvr) \hookrightarrow \Cont_0(D_2,\dvr)\)
  and
  \(g\colon \Cont_0(D_3,\dvr) \hookrightarrow \Cont_0(D_4,\dvr)\) be
  injective, bounded \(\dvr\)\nb-linear maps.  Then the induced
  bounded map
  \[
    f\haotimes g\colon \Cont_0(D_1,\dvr) \haotimes \Cont_0(D_3,\dvr)
    \to \Cont_0(D_2,\dvr) \haotimes \Cont_0(D_4,\dvr)
  \]
  is injective as well.  And here
  \(\Cont_0(D_m,\dvr) \haotimes \Cont_0(D_n,\dvr) \cong \Cont_0(D_m
  \times D_n,\dvr)\).
\end{lemma}

\begin{proof}
  The universal property of the complete bornological tensor product
  implies that
  \(\Cont_0(D_1,\dvr) \haotimes \Cont_0(D_3,\dvr) \cong \Cont_0(D_1
  \times D_3, \dvr)\) for all sets \(D_1\) and~\(D_2\).  There is a
  canonical isomorphism
  \[
    \Cont_0(D_1 \times D_3, \dvr) \overset{\cong}\to
    \Cont_0(D_1,\Cont_0(D_3,\dvr)),
    \qquad f \mapsto (s \mapsto f(s,{\cdot}))
  \]
  Similarly,
  \(\Cont_0(D_1\times D_3, \dvr) \cong \Cont_0(D_3,
  \Cont_0(D_1,\dvr))\).  Now we factorise the map~\(f\hot g\) as
  \begin{multline*}
    \Cont_0(D_1,\dvr) \haotimes \Cont_0(D_3,\dvr)
    \cong \Cont_0(D_1 \times D_3, \dvr)
    \cong \Cont_0(D_1, \Cont_0(D_3,\dvr))
    \\ \overset{g_*}{\hookrightarrow}
    \Cont_0(D_1,\Cont_0(D_4,\dvr))
    \cong  \Cont_0(D_4, \Cont_0(D_1,\dvr))
    \\ \overset{f_*}{\hookrightarrow}
    \Cont_0(D_4,\Cont_0(D_2,\dvr))
    \cong \Cont_0(D_2 \times D_4,\dvr)
    \cong \Cont_0(D_2,\dvr) \haotimes \Cont_0(D_4, \dvr);
  \end{multline*}
  here the maps \(f_*\) and~\(g_*\) are injective because \(f\)
  and~\(g\) are injective.
\end{proof}

\begin{proposition}
  \label{pro:tensor_injective}
  Let \(M_1\), \(W_1\), \(M_2\) and~\(W_2\) be complete,
  torsion-free bornological \(\dvr\)\nb-modules and let
  \(\varphi_j\colon M_j \hookrightarrow W_J\) for \(j=1,2\) be
  injective bounded \(\dvr\)\nb-module maps.  Then
  \(\varphi_1 \hot \varphi_2\colon M_1 \hot M_2 \to W_1 \hot W_2\)
  is injective.
\end{proposition}

\begin{proof}
  Write \(W_1\) and~\(W_2\) as inductive limits as in
  Corollary~\ref{cor:nice_inductive_system}.  Then
  \(W_1 \otimes W_2\) is naturally isomorphic to the inductive limit
  of the inductive system defined by the maps
  \(f_{1,n_1,m_1}\otimes f_{2,n_2,m_2}\colon \Cont_0(D_{n_1},\dvr)
  \otimes \Cont_0(J_{n_2},\dvr) \to \Cont_0(D_{m_1},\dvr) \otimes
  \Cont_0(J_{m_2},\dvr)\), and \(W_1 \hot W_2\) is naturally
  isomorphic to the inductive limit of the inductive system defined
  by the maps
  \(f_{1,n_1,m_1}\haotimes f_{2,n_2,m_2}\colon \Cont_0(D_{n_1},\dvr)
  \haotimes \Cont_0(J_{n_2},\dvr) \to \Cont_0(D_{m_1},\dvr)
  \haotimes \Cont_0(J_{m_2},\dvr)\).  All these bounded maps are
  injective by Lemma~\ref{lem:tensor_injective}.  Therefore, the
  tensor product is isomorphic to an ordinary union of these
  \(\dvr\)\nb-modules, equipped with the bornology cofinally
  generated by these \(\dvr\)\nb-submodules.  The tensor products
  \(M_1 \otimes M_2\) and \(M_1 \hot M_2\) are described similarly,
  and the maps \(\varphi_1\) and~\(\varphi_2\) are described by
  injective maps between the entries of the appropriate inductive
  systems.  Then Lemma~\ref{lem:tensor_injective} shows that
  \(\varphi_1 \hot \varphi_2\) is injective.
\end{proof}

\subsection{The bimodule of differential 1-forms}
\label{sec:Omega_1}

We are going to describe the (complete) bimodule
\(\comb{\Omega}{}^1(A)\) of \emph{noncommutative differential
  \(1\)\nb-forms} over a complete bornological
\(\dvr\)\nb-algebra~\(A\).  It is defined succinctly as the kernel
of the multiplication map \(A^+\hot A^+ \to A^+\).  This is a direct
summand as a bornological \(\dvr\)\nb-module.  Then it is
bornologically closed and a complete bornological \(A\)\nb-bimodule.
The map
\[
  \diff\colon A \to \comb{\Omega}{}^1(A),\qquad
  \diff(x)\defeq 1\otimes x - x \otimes 1,
\]
is the universal bounded derivation into a complete \(A\)\nb-bimodule, that
is, any bounded derivation \(\partial\colon A \to M\) into a complete
\(A\)\nb-bimodule factors uniquely through~\(\diff\).  Namely, there
is a unique bounded bimodule homomorphism \(\comb{\Omega}{}^1(A) \to M\),
\(a_0 \,\diff a_1 \mapsto a_0\cdot \partial(a_1)\).  This
factorisation exists because there
are bornological isomorphisms
\begin{alignat*}{2}
  A^+ \hot A &\to \comb{\Omega}{}^1(A),&\qquad x\otimes y&\mapsto x\,\diff y,\\
  A \hot A^+ &\to \comb{\Omega}{}^1(A),&\qquad x\otimes y&\mapsto
  (\diff x)\cdot y = \diff(x\cdot y) - x\,\diff y.
\end{alignat*}
The first one is left and the second one right \(A\)\nb-linear.

We now relate \(\comb{\Omega}{}^1(A)\) to sections of semi-split,
square-zero extensions of~\(A\) (see
\cite{Meyer:HLHA}*{Theorem~A.53} or
\cite{Cuntz-Quillen:Algebra_extensions}*{Proposition~3.3}).
Let~\(M\) be a complete bornological \(A\)\nb-bimodule.  Give
\(A\oplus M\) the multiplication
\[
  (a_1,m_1) \cdot (a_2,m_2) \defeq
  (a_1\cdot a_2, a_1\cdot m_2 + m_1\cdot a_2).
\]
The inclusion \(M \into A\oplus M\) and the projection
\(A\oplus M\onto A\) form a square-zero extension that splits by the
inclusion homomorphism \(A \hookrightarrow A\oplus M\).

\begin{lemma}
  \label{lem:Omega1_vs_sections}
  Let~\(A\) be a complete bornological algebra and let~\(M\) be a
  complete bornological \(A\)\nb-bimodule.  There is a natural
  bijection between bounded bimodule homomorphisms
  \(\comb{\Omega}{}^1(A) \to M\) and bounded \(\dvr\)\nb-algebra
  homomorphisms \(A \to A\oplus M\) that split the extension
  \(M \into A\oplus M \onto A\).
\end{lemma}

\begin{proof}
  Any bounded linear section \(s\colon A\to A\oplus M\) has the form
  \(a\mapsto (a,\partial(m))\) for a bounded linear map
  \(\partial\colon A\to M\).  And~\(s\) is multiplicative if and
  only if~\(\partial\) is a derivation.  Bounded bimodule maps
  \(\comb{\Omega}{}^1(A) \to M\) are in bijection with bounded
  derivations.
\end{proof}

We shall also apply the definition and the lemma above to incomplete
bornological algebras, where we define \(\Omega^1(A)\) by leaving
out the completions in the construction above.  And we shall use a
variant of \(\Omega^1(A)\) for projective systems of algebras.  In
general, the definition and the lemma above carry over to algebras
in any additive monoidal category.

\subsection{Tensor algebras and noncommutative differential forms}
\label{sec:tensor_nc-forms}

We describe the tensor algebra of a bornological \(\dvr\)\nb-module
and the algebra of differential forms over a bornological algebra
and relate the two.  All this goes back to Cuntz and Quillen.  Their
constructions make sense in any additive monoidal category with
countable direct sums, and we specialise this generalisation of
their constructions to bornological \(\dvr\)\nb-modules and to
complete bornological \(\dvr\)\nb-modules.  We shall mainly use the
uncomplete versions below because we are going to modify tensor
algebras further before completing them.

Let~\(W\) be a bornological \(\dvr\)\nb-module.  Equip
\(W^{\otimes n}\) for \(n\ge 1\) with the tensor product bornology
and \(\tens W\defeq \bigoplus_{n\ge 1} W^{\otimes n}\) with the
direct sum bornology; that is, a subset~\(M\) of~\(\tens W\) is
bounded if and only if it is contained in the image of
\(\bigoplus_{j=1}^n N^{\otimes j}\) for some \(n\ge 1\) and some
bounded submodule \(N\subseteq W\).  The multiplication
\(\tens W \times \tens W \to \tens W\) defined by
\[
  (x_1 \otimes \dotsb \otimes x_n)\cdot (x_{n+1} \otimes \dotsb
  \otimes x_{n+m}) \defeq x_1 \otimes \dotsb \otimes x_{n+m}
\]
makes \(\tens W\) a bornological algebra, called the \emph{tensor
  algebra} of~\(W\).  Let \(\sigma_W\colon W\to \tens W\) be the
inclusion of the first summand.  It is a bounded \(\dvr\)\nb-module
homomorphism, but not an algebra homomorphism.

\begin{lemma}
  \label{lem:tensor_algebra}
  The map \(\sigma_W\colon W\to \tens W\) is the universal bounded
  \(\dvr\)\nb-module map from~\(W\) to a bornological algebra.  That
  is, \(\tens W\) is a bornological \(\dvr\)\nb-algebra and if
  \(f\colon W\to S\) is a bounded \(\dvr\)\nb-module map to a
  bornological \(\dvr\)\nb-algebra~\(S\), then there is a unique
  bounded algebra homomorphism \(f^\#\colon \tens W \to S\) with
  \(f^\#\circ \sigma_W = f\).
\end{lemma}

\begin{proof}
  The multiplication above is well defined and bounded by the
  universal property of the bornological tensor product.  Let
  \(f\colon W\to S\) be a bounded \(\dvr\)\nb-module map.  Then
  there is a unique bounded \(\dvr\)\nb-module map
  \(f^\#\colon \tens W \to S\) with
  \[
    f^\#(x_1 \otimes \dotsb \otimes x_n) \defeq f(x_1) \dotsm f(x_n)
  \]
  for all \(x_1,\dotsc,x_n \in W\).  This is a bounded algebra
  homomorphism.  And it is the unique one with
  \(f^\# \circ \sigma_W = f\).
\end{proof}

Let~\(W\) be a complete bornological \(\dvr\)\nb-module.  The
completion of~\(\tens W\) is
\[
  \comb{\tens} W \defeq \bigoplus_{n\ge 1} W^{\hot n},
\]
the bornological direct sum of the completed tensor products.  By
the universal property of completions, the canonical arrow
\(\comb{\sigma_W}\colon W \to \comb{\tens} W\) is the universal
bounded \(\dvr\)\nb-module map from~\(W\) to a complete bornological
algebra.  That is, \(\comb{\tens} W\) is a complete bornological
\(\dvr\)\nb-algebra and if \(f\colon W\to S\) is a bounded
\(\dvr\)\nb-module map to a complete bornological
\(\dvr\)\nb-algebra~\(S\), then there is a unique bounded algebra
homomorphism \(f^\#\colon \comb{\tens} W \to S\) with
\(f^\#\circ \comb{\sigma_W} = f\).

\begin{remark}
  \label{rem:tensor_tf}
  If~\(W\) is torsion-free, then so is~\(\tens W\).  If~\(W\) is
  complete and torsion-free, then so is~\(\comb{\tens} W\).  This
  uses \cite{Meyer-Mukherjee:Bornological_tf}*{Theorem~4.6 and
    Proposition~4.12} and that completeness and torsion-freeness are
  hereditary for direct sums.
\end{remark}

Let~\(R\) be a bornological \(\dvr\)\nb-algebra.  Then so
is~\(\tens R\).  The identity map on~\(R\) induces a bounded
homomorphism \(p \defeq \id_R^\#\colon \tens R \to R\) by
Lemma~\ref{lem:tensor_algebra}.  Let
\begin{equation}
  \label{eq:def_jens}
  \jens R \defeq \ker (p\colon \tens R \onto R).
\end{equation}
This is a closed two-sided ideal in~\(\tens R\).  The inclusion
\(\jens R \into \tens R\) and the projection
\(p\colon \tens R \onto R\) form an extension of bornological
\(\dvr\)\nb-algebras, which splits by the bounded \(\dvr\)\nb-module
map \(\sigma_R\colon R \to \tens R\).  Similarly, if~\(R\) is a
complete bornological \(\dvr\)\nb-algebra, then there is an
extension of complete bornological \(\dvr\)\nb-algebras
\[
  \comb{\jens} R \into \comb{\tens} R \onto R
\]
that splits by the bounded \(\dvr\)\nb-module
map~\(\comb{\sigma_R}\).

The \emph{unitalisation} of~\(R\) is \(\AU{R}\defeq R\oplus \dvr\) with
the multiplication
\[
  (x,\lambda)\cdot (y,\mu) \defeq (xy+\mu x+\lambda y,\lambda\mu)
\]
for \(x,y\in R\), \(\lambda,\mu\in \dvr\).  So \((0,1)\) is the unit
element in~\(R^+\), which we denote simply by~\(1\).  The inclusion
map \(R \to \AU{R}\) is the universal bounded homomorphism
from~\(R\) to a unital bornological algebra.

Let \(\Omega^0R\defeq R\) and, for \(n\ge 1\), let
\(\Omega^nR \defeq \AU{R} \otimes R^{\otimes n}\), equipped with the
tensor product bornology.  That is, a submodule
\(N\subseteq \Omega^nR\) is bounded if and only if there is a
bounded submodule \(M\subseteq R\) such that~\(N\) is contained in
the image of \(\Omega^n M = \AU{M} \otimes M^{\otimes n}\).  Let
\(\Omega R \defeq \bigoplus_{n\ge 0} \Omega^n R\), equipped with the
direct sum bornology.  We interpret an element
\(x_0\otimes x_1 \otimes \dotsb \otimes x_n \in \Omega^n R\) as a
noncommutative differential form
\(x_0 \, \diff x_1 \dotsc \diff x_n\).  There is a unique structure
of differential graded algebra on \(\Omega R\) whose multiplication
restricts to the given multiplication on \(R=\Omega^0 R\) and whose
differential satisfies
\[
  \diff(x_0 \, \diff x_1 \dotsc \diff x_n) \defeq
  1\cdot \diff x_0 \, \diff x_1 \dotsc \diff x_n.
\]
Namely, the (graded) Leibniz rule dictates that
\[
  x_0 \, \diff x_1 \dotsc \diff x_n \cdot
  x_{n+1} \,\diff x_{n+2} \dotsc \,\diff x_{n+m}
  \defeq
  \sum_{j=0}^{n} (-1)^{n-j}
  x_0\,\diff x_1 \dotsc \diff (x_j\cdot x_{j+1}) \dotsc \diff x_{n+m}.
\]
The \emph{Fedosov product} on a differential graded algebra such as
\(\Omega R\) is defined by
\begin{equation}
  \label{prod:fedosov}
  \xi\Fed\eta \defeq
  \xi\eta - (-1)^{i \cdot j} d(\xi) d(\eta)
  \qquad \text{for }\xi\in \Omega^iR,\ \eta\in\Omega^jR.
\end{equation}
Recall the notation \(M^{(n)} \defeq \sum_{i=1}^n M^i\).  If
\(p,q\ge 0\) and \(M,N\subseteq R\) are bounded
\(\dvr\)\nb-submodules, then
\begin{equation}
  \label{prod:fedosov2}
  \Omega^p M\Fed \Omega^q N\subseteq
  \Omega^{p+q}((M+N)^{(2)})\oplus \Omega^{p+q+2}((M+N)).
\end{equation}
Hence \((\Omega R,\Fed)\) is a bornological algebra.  Its completion
\(\comb{\Omega} R\) is the bornological direct sum
\(\bigoplus_{n\ge 0} \comb{\Omega}^n R\) of the completed
differential forms.  Let \(\Omega^\even R\subseteq\Omega R\) be the
bornological subalgebra of differential forms of even degree.  In
the following, we always equip \(\Omega^\even R\) with the Fedosov
product.

The inclusion map \(R = \Omega^0 R \hookrightarrow \Omega^\even R\)
induces a bounded homomorphism
\begin{equation}
  \label{map:fediso1}
  \tens R\to \Omega^\even R,\qquad
  x_1\otimes\dotsb\otimes x_n\mapsto x_1\Fed\dotsb\Fed x_n,
\end{equation}
by Lemma~\ref{lem:tensor_algebra}, which is, in fact, a bornological
isomorphism.  To understand why, let \(f\colon R\to S\) be a
\(\dvr\)\nb-module map.  Its \emph{curvature} is the
\(\dvr\)\nb-module map
\[
  \omega_f\colon R\otimes R\to S,\qquad
  \omega_f(x,y)=f(x\cdot y)- f(x)\cdot f(y).
\]
It is bounded if~\(f\) is.  The composite of the induced
homomorphism \(f^\#\colon \tens R \to S\) with the inverse of the
map in~\eqref{map:fediso1} must be given by the formula
\begin{equation}
  \label{eq:fsharp_curvature}
  f^\#(x_0 \, \diff x_1 \dotsc \diff x_{2 n})
  = f(x_0) \cdot \omega_f(x_1,x_2) \dotsm \omega_f(x_{2 n-1},x_{2 n})
\end{equation}
because the inclusion map \(R \to \Omega^\even R\) has the curvature
\((x,y) \mapsto x\cdot y - x\Fed y = \diff x\,\diff y\).  Indeed,
this defines a bounded homomorphism
\(f^\#\colon \Omega^\even R \to S\).  So \(\Omega^\even R\) enjoys
the same universal property as \(\tens R\).  Then the map
in~\eqref{map:fediso1} is a bornological isomorphism.

The map \(p\colon \tens R \to R\) corresponds to the map
\(p\colon \Omega^\even R \to R\) that vanishes on \(\Omega^{2 n} R\)
for \(n\ge1\) and is the identity on \(\Omega^0 R = R\).  Therefore,
the isomorphism \(\tens R \cong \Omega^\even R\) maps \(\jens R\)
onto \(\bigoplus_{n\ge 1} \Omega^{2 n} R\).  Then it follows by
induction that the isomorphism maps the ideal~\(\jens R^m\) onto
\(\bigoplus_{n\ge m} \Omega^{2 n} R\).  This simple description of
all the powers \(\jens R^m\) is the main point of rewriting the
tensor algebra using the Fedosov product on the even-degree
differential forms.

\begin{remark}
  \label{rem:jquot}
  The map \(\jens R^{\otimes m}\to \jens R^m\) splits by the bounded
  \(\dvr\)\nb-module map given by
  \[
    a_0 \,\diff a_1\dotsc \diff a_{2(m+n)}\mapsto
    (a_0 \,\diff a_1 \,\diff a_2)\otimes
    \diff a_{2m-3}\,\diff a_{2m-2} \otimes \diff a_{2m-1} \dotsc
    \diff a_{2n}.
  \]
  Thus \(\jens R^{\otimes m}\to \jens R^m\) is a quotient map, and the same is
  true upon completion.
\end{remark}

\subsection{The X-complex}
\label{sec:X-complex}

The \(X\)\nb-complex is another ingredient in the Cuntz--Quillen
approach to cyclic homology theories.  It is defined for algebras in
an additive monoidal category, and we shall specialise its
definition to the additive monoidal category of complete
bornological algebras over \(\dvf\) or~\(\dvr\).

Let \(\comb{\Omega}{}^1(S)/[,]\) be the commutator quotient
of~\(\comb{\Omega}{}^1(S)\), that is, the quotient
of~\(\comb{\Omega}{}^1(S)\) by the closure of the image of
\[
  S\hot \comb{\Omega}{}^1(S) \to \comb{\Omega}{}^1(S),\qquad
  x\otimes \omega \mapsto x\cdot \omega - \omega\cdot x.
\]
With the quotient bornology, this is a complete bornological
\(\dvr\)\nb-module (see
\cite{Meyer-Mukherjee:Bornological_tf}*{Theorem~2.3}).  (The closure
comes in because we take a cokernel in the category of complete
bornological \(\dvr\)\nb-modules, which forces us to make the
quotient separated.)

Let \(q\colon \comb{\Omega}{}^1(S) \to \comb{\Omega}{}^1(S)/[,]\) be
the quotient map.  There is a unique bounded linear map
\(b\colon \comb{\Omega}{}^1(S) \to S\) that satisfies
\(b(x\,\diff y) = x\cdot y - y\cdot x\).  It descends to a bounded
linear map \(\tilde{b}\colon \comb{\Omega}{}^1(S)/[,] \to S\).  The
\(X\)\nb-complex of~\(S\) is the following \(\Z/2\)\nb-graded chain
complex of complete bornological \(\dvr\)\nb-modules:
\[
  X(S)\defeq
  \Bigl(
    \begin{tikzcd}
      S \arrow[r, shift left, "q\circ \diff"]&
      \comb{\Omega}{}^1(S)/[\cdot,\cdot] \arrow[l, shift left, "\tilde{b}"]
    \end{tikzcd}
  \Bigr).
\]
We briefly call \(\Z/2\)\nb-graded chain complexes
\emph{supercomplexes}.  If~\(S\) is a complete bornological
\(\dvf\)\nb-algebra, then \(X(S)\) is even a supercomplex of
complete bornological \(\dvf\)\nb-vector spaces.

\section{Definition of analytic cyclic homology}
\label{sec:definition_theory}
\numberwithin{equation}{section}

Let~\(A\) be a torsion-free, complete bornological
\(\dvr\)\nb-algebra.  We are going to define the analytic cyclic
homology of~\(A\) by a sequence of small steps.  First, let
\[
  R \defeq \tens A,\qquad
  I \defeq \jens A,
\]
be the tensor algebra over~\(A\) and the kernel of the canonical
homomorphism \(\tens A \onto A\).

The second step enlarges~\(R\) to a projective system of tube
algebras relative to powers of the ideal~\(I\):

\begin{definition}
  \label{def:tube_algebra}
  Let~\(R\) be a torsion-free bornological \(\dvr\)\nb-algebra
  and~\(I\) an ideal in~\(R\).  Let~\(I^j\) for \(j\in\N^*\) denote
  the \(\dvr\)\nb-linear span of products \(x_1\dotsm x_j\) with
  \(x_1,\dotsc,x_j\in I\).  The \emph{tube algebra of
    \(I^l\idealin R\)} for \(l\in\N^*\) is
  \[
    \tub{R}{I^l}\defeq \sum_{j=0}^{\infty} \dvgen^{-j}I^{l\cdot j}
    \subseteq R \otimes \dvf
  \]
  with the subspace bornology; this is indeed a
  \(\dvr\)\nb-subalgebra of \(R \otimes \dvf\).  If \(l \ge j\),
  then \(\tub{R}{I^l} \subseteq \tub{R}{I^j}\) is a bornological
  subalgebra.  Let \(\tub{R}{I^\infty}\) be the projective system of
  bornological \(\dvr\)\nb-algebras \((\tub{R}{I^l})_{l\in\N^*}\).
\end{definition}

Since \(\tub{R}{I^l}\) is defined as a bornological submodule of an
\(\dvf\)\nb-vector space, it is bornologically torsion-free.  And
the inclusion \(R \hookrightarrow \tub{R}{I^l}\) induces a
bornological isomorphism
\(\tub{R}{I^l} \otimes \dvf \cong R \otimes \dvf\).

\begin{remark}
  \label{rem:bornology_on_tube}
  In
  \cite{Cortinas-Cuntz-Meyer-Tamme:Nonarchimedean}*{Definition~3.1.19},
  the tube algebra \(\tub{R}{I^l}\) of a bornological
  \(\dvr\)\nb-algebra is equipped with a different bornology,
  namely, the bornology that is generated by subsets bounded
  in~\(R\) and subsets of the form~\(\dvgen^{-1} M^l\) for bounded
  subsets \(M\subseteq I\).  This makes no difference if~\(R\)
  carries the fine bornology.  For general~\(R\), however, the two
  bornologies on the tube algebra need not be the same.  It is easy
  to check that both bornologies induce the same bornology on
  \(\tub{R}{I^l} \otimes \dvf \cong R \otimes \dvf\).  Thus the two
  bornologies coincide if and only if the bornology defined
  in~\cite{Cortinas-Cuntz-Meyer-Tamme:Nonarchimedean} is
  torsion-free.  This concept is introduced only later
  in~\cite{Meyer-Mukherjee:Bornological_tf}.  The more complicated
  bornology defined
  in~\cite{Cortinas-Cuntz-Meyer-Tamme:Nonarchimedean} gives the tube
  algebra the expected universal property for bornological algebras
  that are torsion-free as algebras, but not bornologically
  torsion-free.
\end{remark}

The third step equips \(\tub{R}{I^l}\) for \(l\in\N^*\) with the
linear growth bornology relative to the ideal \(\tub{I}{I^l}\).
This gives a projective system of bornological algebras
\[
  \relling{\tub{R}{I^\infty}}{\tub{I}{I^\infty}}
  = \bigl(\relling{\tub{R}{I^l}}{\tub{I}{I^l}}\bigr)_{l\in\N^*}
\]
because the inclusion homomorphism
\(\tub{R}{I^{l+1}} \hookrightarrow \tub{R}{I^l}\) maps
\(\tub{I}{I^{l+1}}\) to \(\tub{I}{I^l}\).  All these bornological
algebras are torsion-free by Lemma~\ref{lem:relgtf}.

The fourth step applies the completion functor.  By
\cite{Meyer-Mukherjee:Bornological_tf}*{Theorem~4.6}, this gives a
projective system of complete, torsion-free bornological
\(\dvr\)\nb-algebras
\[
  (\tub{R}{I^\infty},\tub{I}{I^\infty})^\updagger
  = \bigl((\tub{R}{I^l},\tub{I}{I^l})^\updagger\bigr)_{l\in \N^*}.
\]

The fifth step is to tensor with~\(\dvf\).  This gives a projective
system of complete bornological \(\dvf\)\nb-algebras
\[
  (\tub{R}{I^\infty},\tub{I}{I^\infty})^\updagger \otimes \dvf
  \defeq ((\tub{R}{I^l},\tub{I}{I^l})^\updagger
  \otimes \dvf)_{l\in\N^*}.
\]

The sixth step is to take the \(X\)\nb-complex.  Being natural, it
extends to a functor from projective systems of complete
bornological algebras to projective systems of supercomplexes.  In
particular, the canonical maps \(\tub{R}{I^{l+1}} \to \tub{R}{I^l}\)
induce bounded chain maps
\[
  \sigma_l\colon X\bigl((\tub{R}{I^{l+1}},\tub{I}{I^{l+1}})^\updagger
  \otimes \dvf\bigr)
  \to X\bigl((\tub{R}{I^l},\tub{I}{I^l})^\updagger \otimes \dvf\bigr).
\]
These define a projective system of supercomplexes of complete
bornological \(\dvf\)\nb-vector spaces, which we denote by
\[
  \HAC(A) \defeq
  X\bigl((\tub{R}{I^\infty},\tub{I}{I^\infty})^\updagger \otimes \dvf\bigr).
\]

The seventh step takes the \emph{homotopy projective limit}
\(\holim \HAC(A)\).  More explicitly, this is the mapping cone of
the chain map
\begin{align*}
  \prod_{l\in \N^*} X\bigl((\tub{R}{I^l},\tub{I}{I^l})^\updagger
  \otimes \dvf\bigr) &\to \prod_{l\in \N^*}
  X\bigl((\tub{R}{I^l},\tub{I}{I^l})^\updagger \otimes \dvf\bigr),\\
  (x_l) &\mapsto \bigl(x_l - \sigma_l(x_{l+1})\bigr)_{l\in\N^*}.
\end{align*}
It is a supercomplex of complete bornological \(\dvf\)\nb-vector
spaces.  The final, eighth step takes its homology:

\begin{definition}
  \label{def:HA_k-algebra}
  The \emph{analytic cyclic homology} \(\HA_*(A)\) of a complete,
  torsion-free bornological \nb-algebra~\(A\) for \(*\in\Z/2\) is
  the homology of \(\holim \HAC(A)\), that is, the quotient of the
  kernel of the differential by the image of the differential.
\end{definition}

\subsection{Bivariant analytic cyclic homology}
\label{sec:bivariant_theory}
\numberwithin{equation}{subsection}

Besides the analytic cyclic homology functor~\(\HA_*\), we also have
the functor \(\HAC\) taking values in suitable homotopy categories
of chain complexes of projective systems of bornological
\(\dvr\)\nb-algebras.  This functor contains more information.  In
particular, it yields a bivariant analytic cyclic homology theory by
letting \(\HA_*(A_1,A_2)\) be the set of morphisms
\(\HAC(A_1) \to \HAC(A_2)\).  This depends on the choice of the
target category, and there is a certain flexibility here.  We do not
pick any choice in this article, but only point out two natural
options.

The analytic cyclic homology computations in this paper often prove
a chain homotopy equivalence \(\HAC(A) \simeq \HAC(B)\), as
supercomplexes of projective systems of bornological \(\dvr\)\nb-modules.
These are equivalences in the homotopy category of supercomplexes,
where homotopy is understood simply as chain homotopy.  In all cases
where we compute \(\HA_*(A)\) in this paper, we actually prove that
\(\HAC(A)\) is chain homotopy equivalent to a supercomplex with
zero boundary map, so that it contains no more information than the
bornological \(\dvf\)\nb-vector space \(\HA_*(A)\).  Homotopy
projective limits are sufficiently compatible with chain homotopies
to preserve chain homotopy equivalence; and this implies an
isomorphism on homology.

A larger class of weak equivalences is used
in~\cite{Cortinas-Valqui:Excision} to define a homotopy category of
chain complexes of projective systems.  A good aspect of this
construction is that it clarifies the role of the homotopy
projective limit: this just replaces a given complex by one that is
weakly equivalent to it and fibrant in a suitable sense, so that the
arrows to it in the homotopy category are the same as chain homotopy
classes of chain maps.  Thus \(\HA_*(A)\) is isomorphic to the space
of arrows from the trivial supercomplex~\(\dvr\) to \(\HAC(A)\) in
the homotopy category of~\cite{Cortinas-Valqui:Excision}.  We will
see later that \(\HAC(\dvr)\) is chain homotopy equivalent to the
trivial supercomplex~\(\dvr\) (see Corollary~\ref{cor:HA_V}).  So
the homotopy category of~\cite{Cortinas-Valqui:Excision} is such
that the bivariant analytic cyclic homology group \(\HA_*(\dvr,A)\)
simplifies to \(\HA_*(A)\).

\section{Analytic nilpotence and analytically quasi-free resolutions}
\label{sec:aqf}
\numberwithin{equation}{section}

Cuntz and Quillen described the periodic cyclic homology of an
algebra~\(A\) as the homology of the \(X\)\nb-complex of a certain
projective system built from the tensor algebra~\(\tens A\)
of~\(A\).  This approach to periodic cyclic homology is the key to
proving that it satisfies excision.  The Cuntz--Quillen approach is
carried over to more analytic versions of periodic cyclic homology
in~\cite{Meyer:HLHA}.  Our proof of excision for~\(\HA_*\) in
Section~\ref{sec:excision} will follow the pattern
in~\cite{Meyer:HLHA}.  In this section, we explain how \(\HA_*\) as
defined above fits into this framework.

\subsection{Pro-Algebras}
\label{sec:pro-algebras}
\numberwithin{equation}{subsection}

An important idea in~\cite{Meyer:HLHA} is that an analytic variant
of periodic cyclic homology is defined by a suitable notion of
``analytic nilpotence''.  This leads to an analytic tensor algebra
of an algebra~\(A\), which is universal among analytically nilpotent
extensions of~\(A\).  It also leads to the concept of analytically
quasi-free algebras.  The theory is set up so that any two
analytically quasi-free, analytically nilpotent extensions of a
given algebra are homotopy equivalent.  In characteristic~\(0\),
this implies that their \(X\)\nb-complexes are chain homotopy
equivalent.  Thus the \(X\)\nb-complex of the analytic tensor
algebra is chain homotopy equivalent to the \(X\)\nb-complex of any
analytically quasi-free resolution of~\(A\).  In this discussion, ``algebras''
are always more complex objects -- such as projective systems of
algebras or bornological algebras -- because there is no suitable
concept of analytic nilpotence for mere algebras without extra
structure.  For the analytic cyclic homology defined above, the
appropriate type of algebra is a projective system of torsion-free,
complete bornological \(\dvr\)\nb-algebras.  For brevity, we call
torsion-free, complete bornological \(\dvr\)\nb-algebras
\emph{algebras} and projective systems of them \emph{pro-algebras}.

A pro-algebra is given by a directed set \((N,\le)\),
algebras~\(A_n\) for \(n\in N\), and bounded algebra homomorphisms
\(\alpha_{m,n}\colon A_n \to A_m\) for \(m,n\in N\) with \(n\ge m\)
that satisfy \(\alpha_{m,m} = \id_{A_m}\) for all \(m\in N\) and
\(\alpha_{m,n}\circ \alpha_{n,p} = \alpha_{m,p}\) for all
\(m,n,p\in N\) with \(p \ge n \ge m\).  The morphism set between
two pro-algebras is
\[
  \Hom\bigl((A_l)_{l\in L},(B_n)_{n\in N}\bigr)
  \defeq \varprojlim_{n} \varinjlim_{l} \Hom(A_l,B_n).
\]
We shall only need pro-algebras \((A_n)_{n\in N}\) where~\(N\) is
countable.  Restricting to a cofinal increasing sequence in~\(N\)
gives an isomorphic pro-algebra with \(N=\N\).  Then the
maps~\(\alpha_{m,n}\) are uniquely determined by
\(\alpha_{n,n+1}\colon A_{n+1} \to A_n\) for \(n\in\N\).

An algebra~\(A\) is also a pro-algebra by taking \(A_n = A\) and
\(\alpha_{n,n+1} \defeq \id_A\) for all \(n\in\N\).  Such projective
systems are called \emph{constant}.  For a pro-algebra
\(A=(A_n,\alpha_{m,n})\), there are canonical morphisms
\(A \to \mathrm{const}(A_n)\) for all \(n\in N\).

The analytic tensor algebra of a torsion-free algebra ~\(A\) is the
torsion-free pro-algebra
\((\tub{\tens A}{\jens A^\infty},\tub{\jens A}{\jens
  A^\infty})^\updagger\) in the above definition of analytic cyclic
homology.  This comes with a canonical homomorphism to~\(A\), whose
kernel is the pro-algebra
\((\tub{\jens A}{\jens A^\infty})^\updagger\).  This projective
system of complete, torsion-free bornological algebras has two
important extra properties: it is semi-dagger -- hence dagger -- and
nilpotent mod~\(\dvgen\) -- this concept will be defined below.  A
pro-algebra with these two properties is called analytically
nilpotent.
The tube algebra construction and the relative dagger completion in
the construction of the analytic tensor algebra are the universal
way to make a pro-algebra extension of~\(A\) have an analytically
nilpotent kernel.

Any functor from algebras to algebras extends canonically to an
endofunctor on the category of pro-algebras by applying it
entrywise.  The definition of analytic cyclic homology already used
this extension to pro-algebras for completions and tensor products
with~\(\dvf\).  The constructions of \(\tens A\) and~\(\jens A\) for
algebras are also functors and thus extend to pro-algebras.  So is
the tensor product bifunctor \({-} \hot {-}\), which extends to
pro-algebras by
\begin{multline*}
  (A_n,\alpha_{m,n})_{m,n\in N} \hot
  (B_n,\beta_{m,n})_{m,n\in N'} \\\defeq
  (A_{n_1} \hot B_{n_2}, \alpha_{m_1,n_1}\hot
  \beta_{m_2,n_2})_{m_1,n_1\in N, m_2,n_2\in N'}.
\end{multline*}
In particular, we may tensor a pro-algebra with an algebra such
as~\(\dvr[t]^\updagger\), viewed as a constant pro-algebra.

\begin{definition}
  \label{def:daghomo}
  An \emph{elementary dagger homotopy} between two morphisms of
  pro-algebras \(f_0,f_1\colon A \rightrightarrows B\) is a morphism
  of pro-algebras \(f\colon A \to B\hot \dvr[t]^\updagger\) that
  satisfies \((\id_A \otimes \ev_t)\circ f = f_t\) for \(t=0,1\).
  We call \(f_0,f_1\) \emph{elementary dagger homotopic} if there is
  such a homotopy.  \emph{Dagger homotopy} is the equivalence
  relation generated by elementary dagger homotopy.
\end{definition}

\subsection{The universal property of the tube algebra construction}
\label{subsec:nilpotence}
\numberwithin{equation}{subsection}

First, we generalise the construction of tube algebras to
pro-algebras.  Actually, in this subsection, we drop the
completeness assumption for algebras because tube algebras are
usually incomplete.  So ``algebras'' are torsion-free bornological
algebras and pro-algebras are projective systems of such algebras
until the end of this subsection.

An \emph{ideal in a pro-algebra} \(A=(A_n, \alpha_{m,n})_{i\in N}\)
is a family of ideals \(I_n\idealin A_n\) with
\(\alpha_{m,n}(I_n)\subseteq I_m\) for all \(n,m\in N\) with
\(n\ge m\); then~\(\alpha_{m,n}\) induces homomorphisms
\(\tub{A_n}{I_n^l} \to \tub{A_m}{I_m^l}\) for all \(l\in\N^*\),
which intertwine the inclusion maps
\(\tub{A_n}{I_n^l} \hookrightarrow \tub{A_n}{I_n^j}\) for
\(l \ge j\).  These homomorphisms form a pro-algebra
\[
  \tub{A}{I^\infty} \defeq \bigl(\tub{A_n}{I_n^l}\bigr)_{n\in N, l\in\N^*}.
\]
If \(l\in\N^*\), then
\(\tub{A}{I^l} \defeq \bigl(\tub{A_n}{I_n^l}\bigr)_{n\in N}\) is a
pro-algebra.  The pro-algebra \(\tub{A}{I^l}\) for
\(l\in\N^*\cup \{\infty\}\) contains \(\tub{I}{I^l}\) as an ideal.
Since \(A_n \subseteq \tub{A_n}{I_n^l}\) for all \(n\in N\),
\(l\in\N^*\), the inclusion maps define a pro-algebra homomorphism
\(\iota_{A,I}\colon A \to \tub{A}{I^\infty}\).

\begin{definition}
  \label{def:nilpotent_mod_pi}
  A pro\nb-algebra \((A_n,\alpha_{m,n})_{n\in N}\) is
  \emph{nilpotent mod~\(\dvgen\)} if, for each \(m\in N\), there are
  \(n\in N_{\ge m}\) and \(l\in\N^*\) such that
  \(\alpha_{m,n}(A_n^l) \subseteq \dvgen A_m\); here~\(A_n^l\)
  denotes the \(\dvr\)\nb-submodule generated by all products
  \(x_1\dotsm x_l\) of \(l\)~factors in~\(A_n\).
\end{definition}

\begin{remark}
  \label{rem:pro-nilpotent}
  Let \(A=(A_n,\alpha_{m,n})_{m,n\in N}\) be a pro-algebra.  Let
  \(A/(\dvgen)\) be the projective system of \(\resf\)\nb-algebras
  formed by the quotients \(A_n/(\dvgen)\) with the homomorphisms
  induced by~\(\alpha_{m,n}\).  By definition, \(A\) is nilpotent
  mod~\(\dvgen\) if and only if~\(A/(\dvgen)\) has the following
  property: for each \(n\in N\) there are \(m\in N\) and
  \(l\in\N^*\) such that the \(l\)\nb-fold multiplication map
  \((A_m/(\dvgen))^{\otimes l} \to A_n/(\dvgen)\) is zero.  This is
  equivalent to the definition that a projective system of
  \(\resf\)\nb-algebras is pro-nilpotent in
  \cite{Meyer:HLHA}*{Definition~4.3}.
\end{remark}

\begin{proposition}
  \label{pro:make-pro-nilpotent}
  Let \(A\) and~\(B\) be pro-algebras and let \(I\) and~\(J\) be
  ideals in \(A\) and~\(B\), respectively.  Let
  \(\varphi\colon A\to B\) be a pro-algebra morphism that restricts
  to a pro-algebra morphism \(I\to J\).  Let
  \(\iota_{A,I}\colon A \to \tub{A}{I^\infty}\) denote the canonical
  pro-algebra morphism.
  \begin{enumerate}
  \item \label{pro:make-pro-nilpotent1}%
    The pro-algebra \(\tub{I}{I^\infty}\) is nilpotent mod~\(\dvgen\).
  \item \label{pro:make-pro-nilpotent2}%
    If~\(J\) is nilpotent mod~\(\dvgen\), then there is a unique
    morphism \(\bar\varphi\colon \tub{A}{I^\infty} \to B\) with
    \(\bar\varphi\circ \iota_{A,I} = \varphi\).  It restricts to a
    morphism \(\tub{I}{I^\infty}\to J\).
  \item \label{pro:make-pro-nilpotent3}%
    There is a unique morphism
    \(\varphi_*\colon \tub{A}{I^\infty} \to \tub{B}{J^\infty}\) with
    \(\varphi_*\circ \iota_{A,I} = \iota_{B,J} \circ \varphi\).  It
    restricts to a morphism
    \(\tub{I}{I^\infty}\to \tub{J}{J^\infty}\).
  \end{enumerate}
\end{proposition}

\begin{proof}
  Write \(A=(A_n,\alpha_{m,n})_{n\in N}\), \(I=(I_n)_{n\in N}\) with
  ideals~\(I_n\) in~\(A_n\) with \(\alpha_{m,n}(I_n) \subseteq I_m\)
  and \(B=(B_n,\beta_{m,n})_{n\in N'}\), \(J=(J_n)_{n\in N'}\) with
  ideals~\(J_n\) in~\(B_n\) with \(\beta_{m,n}(J_n) \subseteq J_m\).
  The tube algebra \(\tub{A}{I^\infty}\) is the projective limit of
  the tube algebras \(\tub{A_n}{I_n^\infty}\) in the category of
  pro-algebras.

  Being nilpotent mod~\(\dvgen\) is hereditary for
  projective limits.  So it suffices to
  prove~\ref{pro:make-pro-nilpotent1} when~\(A\) is a constant
  pro-algebra.  Fix \(n\in\N^*\) and let \(m = 2n\), \(l=n\).  Then
  \begin{equation}
    \label{eq:tube_power}
    \tub{I}{I^m}^l
    = \tub{I}{I^{2n}}^n
    = \biggl( I + \sum_{j=1}^\infty \dvgen^{-j} I^{2 n j}\biggr)^n
    \subseteq I^n + \sum_{j=1}^\infty \dvgen^{-j} I^{2 n j}
  \end{equation}
  because \(\sum_{j=1}^\infty \dvgen^{-j} I^{2 n j}\) is an ideal in
  \(\tub{A}{I^{2 n}}\).  Since \(\dvgen^{-1} I^n\) and
  \(\dvgen^{-2 j} I^{2 n j}\) are contained in \(\tub{I}{I^n}\), all
  summands on the right hand side of~\eqref{eq:tube_power} are
  contained in \(\dvgen\cdot \tub{I}{I^n}\).  Thus
  \(\tub{I}{I^\infty}\) is nilpotent mod~\(\dvgen\).

  We prove statement~\ref{pro:make-pro-nilpotent2}.  The morphism
  \(\varphi\colon A\to B\) is described by a coherent family of
  \(\dvr\)\nb-algebra homomorphisms
  \(\varphi_n\colon A_{\psi(n)}\to B_n\) for all \(n\in N'\).
  Each~\(B_n\) is torsion-free by our definition of ``algebra''.
  Then the homomorphism~\(\varphi_n\) is determined by
  \(\varphi_n\otimes \id_\dvf \colon A_{\psi(n)} \otimes \dvf \to
  B_n\otimes \dvf\).  By construction,
  \(\tub{A_\nu}{I^m} \otimes \dvf = A_\nu \otimes \dvf\) for all
  \(\nu\in N\), \(m\in\N^*\).  Thus a factorisation of~\(\varphi\)
  through \(\tub{A}{I^\infty}\) is unique if it exists.

  Fix \(n\in N'\).  Since~\(J\) is nilpotent mod~\(\dvgen\), there
  are \(m\in N'_{\ge n}\) and \(l\in\N^*\) with
  \(\beta_{n,m}(J_m^l) \subseteq \dvgen \cdot J_n\).
  Since~\(\varphi\) is coherent, there is \(\nu \in N_{\ge\psi(m)}\)
  with
  \(\beta_{n,m}\circ \varphi_m \circ \alpha_{\psi(m),\nu} =
  \varphi_n\circ \alpha_{n,\nu}\).  Since~\(\varphi\) restricts to a
  morphism \(I\to J\), we may also arrange that
  \(\varphi_m\circ\alpha_{\psi(m),\nu}(I_\nu) \subseteq J_m\) by
  increasing~\(\nu\) if necessary.  Hence
  \[
    \varphi_n\circ \alpha_{n,\nu}(I_\nu^l)
    = \beta_{n,m}\circ \varphi_m\circ\alpha_{\psi(m),\nu}(I_\nu^l)
    \subseteq \beta_{n,m}(J_m^l) \subseteq \dvgen\cdot J_n.
  \]
  Thus the homomorphism
  \((\varphi_n\circ \alpha_{n,\nu}) \otimes \id_\dvf\colon A_\nu
  \otimes \dvf \to B_n \otimes \dvf\) maps the tube algebra
  \(\tub{A_\nu}{I_\nu^l} \subseteq A_\nu \otimes \dvf\) into
  \(B_n \subseteq B_n \otimes \dvf\) and
  \(\tub{I_\nu}{I_\nu^l} \subseteq A_\nu \otimes \dvf\) into
  \(J_n \subseteq B_n \otimes \dvf\).  This gives a homomorphism
  \(\bar\varphi_n\colon \tub{A_\nu}{I_\nu^l} \to B_n\) with
  \(\bar\varphi_n\circ \iota_{A_\nu,I_\nu^l} = \varphi_n\circ
  \alpha_{n,\nu}\).  Since
  \(\tub{A_\nu}{I^m} \subseteq A_\nu \otimes \dvf\), the
  homomorphisms~\(\bar\varphi_n\) inherit the coherence property of
  a pro-algebra morphism from the maps~\(\varphi_n\).

  We prove statement~\ref{pro:make-pro-nilpotent3} of the
  proposition.  We compose \(\varphi\colon A\to B\) with the
  canonical map \(B\to\tub{B}{J^\infty}\) to get a morphism
  \(A\to \tub{B}{J^\infty}\).  It restricts to a morphism
  \(I\to J \to \tub{J}{J^\infty}\).  The ideal \(\tub{J}{J^\infty}\)
  in \(\tub{B}{J^\infty}\) is nilpotent mod~\(\dvgen\)
  by~\ref{pro:make-pro-nilpotent1}.
  So~\ref{pro:make-pro-nilpotent2} shows that our morphism extends
  uniquely to a morphism \(\tub{A}{I^\infty} \to \tub{B}{J^\infty}\)
  that maps \(\tub{I}{I^\infty}\) to \(\tub{J}{J^\infty}\).
\end{proof}

We summarise the tube algebra construction in category-theoretic
language.  Let~\(\mathfrak{Pro}\) be the category whose objects are
pairs \((A,I)\), where~\(A\) is a pro-algebra and~\(I\) is an ideal
in~\(A\) and whose morphisms are pro-algebra morphisms that
restrict to a morphism between the ideals.  The pairs~\((A,I)\)
where~\(I\) is nilpotent mod~\(\dvgen\) form a
subcategory~\(\mathfrak{Pro}_\mathrm{nil}\) in~\(\mathfrak{Pro}\).
The first two statements in Proposition~\ref{pro:make-pro-nilpotent}
say that the canonical arrow
\((A,I)\to (\tub{A}{I^\infty},\tub{I}{I^\infty})\) is a universal
arrow from~\((A,I)\) to an object
in~\(\mathfrak{Pro}_\mathrm{nil}\).  Thus
\(\mathfrak{Pro}_\mathrm{nil}\) is a reflective subcategory
in~\(\mathfrak{Pro}\) and the reflector acts on objects by
\((A,I)\mapsto (\tub{A}{I^\infty},\tub{I}{I^\infty})\).  Its
functoriality is Proposition
\ref{pro:make-pro-nilpotent}.\ref{pro:make-pro-nilpotent3}.
If~\(I\) is already nilpotent mod~\(\dvgen\), then it follows that
the identity map on~\(A\) extends uniquely to an isomorphism of
pro-algebras \(\tub{A}{I^\infty} \cong A\).

The heredity properties of nilpotence mod~\(\dvgen\) proven in the
following proposition are needed by the analytic cyclic homology
machinery in~\cite{Meyer:HLHA}.

\begin{proposition}
  \label{pro:extension_nilpotent_mod_pi}
  The class of nilpotent mod~\(\dvgen\) pro-algebras is closed under
  the following operations:
  \begin{itemize}
  \item Let \(A\overset{i}\into B \overset{p}\onto C\) be an
    extension of pro-algebras.  If \(A\) and~\(C\) are nilpotent
    mod~\(\dvgen\), then so is~\(B\), and vice versa.
  \item A pro-subalgebra \(D\subseteq B\) is nilpotent
    mod~\(\dvgen\) if~\(B\) is so and~\(B/D\) is isomorphic to a
    projective system of torsion-free bornological
    \(\dvr\)\nb-modules.
  \item Being nilpotent mod~\(\dvgen\) is hereditary for projective
    limits.
  \item A tensor product \(A\hot B\) is nilpotent mod~\(\dvgen\) if
    \(A\) or~\(B\) is nilpotent mod~\(\dvgen\).
  \end{itemize}
\end{proposition}

\begin{proof}
  Remark~\ref{rem:pro-nilpotent} translates all these statements to
  statements about the class of pro-nilpotent projective systems of
  \(\resf\)\nb-algebras.  In this way, the statements follow from
  \cite{Meyer:HLHA}*{Theorem~4.4}.  We briefly explain direct proofs
  for the first two claims.  The claims about projective limits and
  tensor products are easy and left to the reader.

  As in~\cite{Meyer:HLHA}, we may write any extension of
  pro-algebras \(A\overset{i}\into B \overset{p}\onto C\) as a
  projective system of extensions
  \(A_n\overset{i_n}\into B_n \overset{p_n}\onto C_n\), with
  morphisms of extensions
  \[
  \begin{tikzcd}
    A_n \arrow[r,rightarrowtail, "i_n"] \arrow[d,"\alpha_{m,n}"] &
    B_n \arrow[r, twoheadrightarrow, "p_n"] \arrow[d,"\beta_{m,n}"] &
    C_n  \arrow[d,"\gamma_{m,n}"] \\
    A_m \arrow[r,rightarrowtail, "i_m"] &
    B_m \arrow[r, twoheadrightarrow, "p_m"] &
    C_m
  \end{tikzcd}
  \]
  for \(n\ge m\) as structure maps (this construction is also
  explained during the proof of Proposition~\ref{pro:extension_an}
  below).  Assume that \(A\) and~\(C\) are nilpotent mod~\(\dvgen\).
  Pick \(m\in N\).  There are \(n_1\in N_{\ge m}\) and
  \(j_1\in\N^*\) so that
  \(\alpha_{m,n_1}(A_{n_1}^{j_1})\subseteq \dvgen\cdot A_m\).  And
  there are \(n_2\in N_{\ge n_1}\) and \(j_2\in\N^*\) so that
  \(\gamma_{n_1,n_2}(C_{n_2}^{j_2})\subseteq \dvgen\cdot C_{n_1}\).
  Then
  \(p_{n_1}(\beta_{n_1,n_2}(B_{n_2}^{j_2})) \subseteq \dvgen\cdot
  C_{n_1}\).  This implies
  \(\beta_{n_1,n_2}(B_{n_2}^{j_2}) \subseteq \dvgen\cdot B_{n_1} +
  i_{n_1}(A_{n_1})\).  Then
  \begin{multline*}
    \beta_{m,n_2}(B_{n_2}^{j_1\cdot j_2})
    \subseteq \beta_{m,n_1}(\dvgen\cdot B_{n_1} + i_{n_1}(A_{n_1}))^{j_1}
    \subseteq \dvgen \cdot B_m + i_m(\alpha_{m,n_1}(A_{n_1}^{j_1}))
    \\\subseteq \dvgen \cdot B_m + i_m(\dvgen A_m)
    \subseteq \dvgen\cdot B_m.
  \end{multline*}
  So~\(B\) is nilpotent mod~\(\dvgen\).  Conversely, if~\(B\) is
  nilpotent mod~\(\dvgen\), then~\(C\) is nilpotent mod~\(\dvgen\)
  because \(p_m(B_m) = C_m\) and
  \(p_m(\dvgen \cdot B_m) = \dvgen \cdot C_m\).  The claim
  that~\(A\) is nilpotent mod~\(\dvgen\) if~\(B\) is follows from
  the claim about pro-subalgebras.

  Given a pro-subalgebra \(D\subseteq B\), we may write
  \(B= (B_n,\beta_{m,n})_{n\in N}\) and
  \(D=(D_n,\delta_{m,n})_{n\in N}\) so that \(D_n\subseteq B_n\) for
  all \(n\in N\) and
  \(\delta_{m,n} = \beta_{m,n}|_{D_n}\colon D_n \to D_m\) for all
  \(m,n\in N\) with \(m \le n\).  Let \(m\in N\).  Since~\(B/D\) is
  isomorphic to a projective system of torsion-free bornological
  \(\dvr\)\nb-modules, there is \(n\in N_{\ge m}\) so that the
  structure map \(B_n/D_n \to B_m/D_m\) kills all elements
  \(x\in B_n/D_n\) with \(\dvgen \cdot x=0\).  Equivalently, if
  \(x\in B_n\) satisfies \(\dvgen\cdot x\in D_n\), then
  \(\beta_{m,n}(x) \in D_m\).  Thus
  \(\beta_{m,n}(\dvgen\cdot B_n\cap D_n) \subseteq \dvgen \cdot
  D_m\).  If~\(B\) is nilpotent mod~\(\dvgen\), then there are
  \(l\in N_{\ge n}\) and \(j\in\N^*\) with
  \(\beta_{n,l}(B_l^j) \subseteq \dvgen\cdot B_n\).  Hence
  \[
  \delta_{m,l}(D_l^j)
  \subseteq \delta_{m,n} (\delta_{n,l}(D_l^j))
  \subseteq \beta_{m,n} (\dvgen\cdot B_n \cap D_n)
  \subseteq \dvgen\cdot D_m.
  \]
  Thus~\(D\) is nilpotent mod~\(\dvgen\).
\end{proof}

\subsection{Analytically nilpotent pro-algebras}
\label{sec:analytic_nilpotence}

From now on, ``algebra'' means a complete, torsion-free bornological
algebra.

\begin{definition}
  \label{def:an_nilpotent}
  A pro-algebra~\(J\) is \emph{analytically nilpotent} if it is
  isomorphic to a pro-dagger algebra and nilpotent mod~\(\dvgen\).
  It is \emph{square-zero} if its multiplication map is~\(0\).  An
  extension of pro-algebras \(J\into E \onto A\) is
  \emph{analytically nilpotent} or \emph{square-zero} if~\(J\) is
  analytically nilpotent or square-zero, respectively.
\end{definition}

\begin{definition}
  A \emph{pro-linear map} between two pro-algebras is a morphism of
  projective systems of bornological \(\dvr\)\nb-modules between
  them; so pro-linear maps need not be multiplicative.  An extension
  of pro-algebras \(J\into E \onto A\) is \emph{semi-split} if it
  splits by a pro-linear map.
\end{definition}

\begin{definition}
  \label{def:an_quasi-free}
  A pro-algebra~\(A\) is \emph{analytically quasi-free} if any
  semi-split analytically nilpotent extension \(J \into E \onto A\)
  splits by a pro-algebra homomorphism \(A \to E\).  It is
  \emph{quasi-free} if any semi-split square-zero extension
  \(J \into E \onto A\) splits by a pro-algebra homomorphism
  \(A \to E\).
\end{definition}

The following lemma gives an equivalent reformulation of the last
definition:

\begin{lemma}
  \label{lem:quasi-free_with_lifts}
  A pro-algebra~\(A\) is analytically quasi-free if and only if, for
  any semi-split analytically nilpotent extension
  \(J \into E \onto B\), any homomorphism \(f\colon A\to B\) lifts
  to a homomorphism \(A\to E\).  A pro-algebra~\(A\) is quasi-free
  if and only if, for any semi-split square-zero extension
  \(J \into E \onto B\), any homomorphism \(f\colon A\to B\) lifts
  to a homomorphism \(A\to E\).
\end{lemma}

\begin{proof}
  We may pull the given extension back to a semi-split extension
  \(J \into \hat{E} \onto A\), such that a section \(A \to \hat{E}\)
  is equivalent to a lifting of~\(f\).
\end{proof}

\begin{remark}
  \label{rem:an_quasi-free_is_quasi-free}
  A pro-algebra is square-zero if and only if it is isomorphic to a
  projective system of torsion-free complete bornological
  \(\dvr\)\nb-modules, each equipped with the zero map as
  multiplication.  Then it is analytically nilpotent.  As a
  consequence, analytically quasi-free algebras are quasi-free.
\end{remark}

\begin{proposition}
  \label{pro:V_analytically_quasi_free}
  The base ring~\(\dvr\) viewed as a constant pro-algebra is
  analytically quasi-free.
\end{proposition}

\begin{proof}
  The proof follows
  \cite{Cuntz-Quillen:Cyclic_nonsingularity}*{Section~12}.  Let
  \(J\into E \onto Q\) be a semi-split, analytically nilpotent
  extension of pro-algebras.  Analytic quasi-freeness of~\(\dvr\) is
  equivalent to the assertion that any idempotent in~\(Q\) lifts to
  an idempotent in~\(E\).  Here by an idempotent in a pro-algebra
  \(A=(A_n)_n\), we mean a collection \(a = (a_n)_n\) of idempotents
  \(a_n \in A_n\).  Each \(a_n \in A_n\) is equivalent to a
  homomorphism \(\dvr \to A_n\).

  Let \(\dot{e}=(\dot{e}_n)_n \in Q\) be an idempotent and let
  \(e \in E\) be the image of~\(\dot{e}\) under a pro-linear section
  for \(E \overset{p}\onto Q\).  Let \(x\defeq e - e^2 \in J\).  We
  use an Ansatz by Cuntz and Quillen to find an idempotent
  \(\hat{e} \in E\) with \(e - \hat{e} \in J\).  Namely, we assume
  \(\hat{e} = e+ (2e - 1) \varphi(x)\) for some power series
  \(\varphi \in t\Z[[t]]\).  As~\(J\) is nilpotent mod~\(\dvgen\),
  for every \(l \in N\), there are \(m(l)\ge l\) and
  \(j(l) \in N^*\) with \(x_{m(l)}^{j(l)} = \dvgen y_l\).  To
  simplify notation, we simply write this as \(x^j = \dvgen y\) for
  some \(y \in J\) and \(j \in N^*\).  Finally, since~\(J\) is also
  a pro-dagger algebra, \(\varphi(x) \in J\) for all
  \(\varphi \in t\Z[[t]]\).  We compute
  \[
    \hat{e}^2 - \hat{e} = (\varphi(x)^2 + \varphi(x))(1-4x) - x.
  \]
  So \(\hat{e}^2 = \hat{e}\) if and only if
  \(\varphi(x)^2 + \varphi(x) = \frac{x}{1-4x}\).  This is solved by
  \(\varphi(x) \defeq \sum_{n=1}^\infty \binom{2n-1}{n}x^n \).  This
  defines an element of~\(J\).  Then~\(\hat{e}\) is the desired
  idempotent lifting.
\end{proof}

\begin{proposition}
  An algebra~\(A\) is analytically quasi-free if and only if its
  unitalisation~\(A^+\) is analytically quasi-free.
\end{proposition}

\begin{proof}
  Proposition~\ref{pro:V_analytically_quasi_free} implies this as in
  the proof of \cite{Meyer:HLHA}*{Proposition~5.53}.

\end{proof}

\begin{proposition}
  \label{pro:sum_quasi-free}
  Let \((A_n)_{n\in \N}\) be a sequence of unital, analytically
  quasi-free pro-algebras.  Then \(\bigoplus_{n\in \N} A_n\) is
  analytically quasi-free.
\end{proposition}

\begin{proof}
  The proof of \cite{Meyer:HLHA}*{Proposition~5.53} carries over to
  this context.
\end{proof}

\begin{corollary}
  \label{cor:sum_V_quasi-free}
  The direct sum \(\bigoplus_{n\in \N} \dvr\) is analytically quasi-free.
\end{corollary}

\begin{proposition}
  \label{pro:an_quasi-free_extension_universal}
  Let \(J_i \into E_i \onto A_i\) for \(i=1,2\) be semi-split,
  analytically nilpotent extensions of pro-algebras.  Assume
  that~\(E_1\) is analytically quasi-free.
  \begin{enumerate}
  \item \label{pro:an_quasi-free_extension_universal_1}%
    Any pro-algebra morphism \(f\colon A_1 \to A_2\) lifts to a
    morphism of extensions
    \[
      \begin{tikzcd}
        J_1 \arrow[r,rightarrowtail] \arrow[d] &
        E_1 \arrow[r, twoheadrightarrow, "q_1"] \arrow[d, "\hat{f}"] &
        A_1 \arrow[d, "f"] \\
        J_2 \arrow[r,rightarrowtail] &
        E_2 \arrow[r, twoheadrightarrow, "q_2"] &
        A_2
      \end{tikzcd}
    \]
    This lifting is unique up to dagger homotopy.

  \item \label{pro:an_quasi-free_extension_universal_2}%
    Let \(\hat{f},\hat{g}\colon E_1 \rightrightarrows E_2\) be
    pro-algebra homomorphisms that lift homomorphisms
    \(f,g\colon A_1 \rightrightarrows A_2\).  Then an
    elementary dagger homotopy
    \(h \colon A_1 \to A_2 \hot \dvr[t]^\updagger\) between \(f\)
    and~\(g\) lifts to an elementary dagger homotopy
    \(\hat{h}\colon E_1 \to E_2 \hot \dvr[t]^\updagger\) between
    \(\hat{f}\) and~\(\hat{g}\).
  \item \label{pro:an_quasi-free_extension_universal_3}%
    Any elementary dagger homotopy
    \(A_1 \to A_2 \hot \dvr[t]^\updagger\) lifts to an elementary
    dagger homotopy \(E_1 \to E_2 \hot \dvr[t]^\updagger\).
  \end{enumerate}
\end{proposition}

\begin{proof}
  Let \(f\colon A_1 \to A_2\) be a pro-algebra homomorphism.
  Since~\(E_1\) is analytically quasi-free and the extension
  \(J_2 \into E_2 \onto A_2\) is semi-split and analytically
  nilpotent, the homomorphism \(f\circ q_1\) lifts to a homomorphism
  \(\hat{f}\colon E_1 \to E_2\).  Since
  \(q_2\circ \hat{f} = f\circ q_1\) vanishes on~\(J_1\), \(\hat{f}\)
  restricts to a homomorphism \(J_1 \to J_2\).  Thus~\(\hat{f}\)
  gives a morphism of extensions.

  The uniqueness claim
  in~\ref{pro:an_quasi-free_extension_universal_1} follows
  from~\ref{pro:an_quasi-free_extension_universal_2} by taking
  \(f=g\).  And~\ref{pro:an_quasi-free_extension_universal_3}
  follows from \ref{pro:an_quasi-free_extension_universal_1}
  and~\ref{pro:an_quasi-free_extension_universal_2}.  So it remains
  to prove~\ref{pro:an_quasi-free_extension_universal_2}.  Assume
  that we are in the situation
  of~\ref{pro:an_quasi-free_extension_universal_2}.  Let
  \(\ev_0,\ev_1\colon A_2 \hot \dvr[t]^\updagger \rightrightarrows
  A_2\) and
  \(\ev_0,\ev_1\colon E_2 \hot \dvr[t]^\updagger \rightrightarrows
  E_2\) denote the evaluation homomorphisms.  Form the pull-back
  pro-algebra
  \[
    \begin{tikzcd}[column sep=1.5cm]
      \mathcal{E} \arrow[r] \arrow[d] &
      E_2 \oplus E_2 \arrow[d,"q_2 \oplus q_2"] \\
      A_2 \hot \dvr[t]^\updagger \arrow[r, "(\ev_0\ \ev_1)" ] &
      A_2 \oplus A_2.
    \end{tikzcd}
  \]
  The universal property of the pull back gives pro-algebra
  homomorphisms
  \begin{gather*}
    q\defeq (\ev_0,\ev_1,q_2 \hot \id_{\dvr[t]^\updagger})_*\colon
    E_2 \hot \dvr[t]^\updagger \to \mathcal{E},\\
    (\hat{f},\hat{g},h\circ q_1)_*\colon
    E_1 \to \mathcal{E},
  \end{gather*}
  because \(\hat{f}\) and~\(\hat{g}\) lift \(\ev_t\circ h\) for
  \(t=0,1\), respectively.  Let
  \[
    \dvr[t]^\updagger_0 \defeq
    \setgiven{\varphi\in \dvr[t]^\updagger}{\varphi(0)=0,\ \varphi(1)=0}.
  \]
  We claim that~\(q\) is part of a semi-split, analytically
  nilpotent extension of pro-algebras
  \begin{equation}
    \label{eq:extension_for_lifting_homotopies}
    J_2 \hot \dvr[t]^\updagger_0\into
    E_2 \hot \dvr[t]^\updagger \onto \mathcal{E}.
  \end{equation}
  To see this, we forget multiplications and treat everything as a
  projective system of bornological \(\dvr\)\nb-modules.  In this
  category, a pro-linear section \(s\colon A_2 \to E_2\) for the
  semi-split extension \(J_2 \to E_2 \to A_2\) gives a direct sum
  decomposition \(E_2 \cong J_2 \oplus A_2\).  And
  \(\dvr[t]^\updagger \cong \dvr[t]^\updagger_0 \oplus \dvr \oplus
  \dvr\), where the latter two summands are, say, spanned by the
  functions \(1-t\) and~\(t\).  This induces decompositions
  \[
    A_2 \hot \dvr[t]^\updagger \cong
    (A_2 \hot \dvr[t]^\updagger_0) \oplus A_2 \oplus A_2,\qquad
    E_2 \hot \dvr[t]^\updagger \cong
    (E_2 \hot \dvr[t]^\updagger_0) \oplus E_2 \oplus E_2,
  \]
  such that \((\ev_0,\ev_1)\) is the projection to the second and
  third summand both for \(A_2\) and~\(E_2\).  These direct sum
  decompositions imply
  \[
    E_2 \hot \dvr[t]^\updagger
    \cong (J_2 \hot \dvr[t]^\updagger_0)
    \oplus (A_2 \hot \dvr[t]^\updagger_0)
    \oplus E_2 \oplus E_2
    \cong (J_2 \hot \dvr[t]^\updagger_0)
    \oplus  \mathcal{E}.
  \]
  And this proves the claim.

  Corollary~\ref{cor:tensor_dagger} and
  Proposition~\ref{pro:extension_nilpotent_mod_pi} imply that the
  tensor product \(J_2 \hot \dvr[t]^\updagger_0\) is analytically
  nilpotent.  Since~\(E_1\) is analytically quasi-free, the
  homomorphism \((\hat{f},\hat{g},h\circ q_1)\) lifts to a
  homomorphism \(\hat{h}\colon E_1 \to E_2 \hot \dvr[t]^\updagger\)
  in the extension~\eqref{eq:extension_for_lifting_homotopies}.
  This finishes the proof
  of~\ref{pro:an_quasi-free_extension_universal_2}.
\end{proof}

\begin{corollary}
  \label{cor:an_quasi-free_extension_universal}
  Any two analytically quasi-free, analytically nilpotent extensions
  of a pro-algebra are dagger homotopy equivalent.
\end{corollary}

\begin{proof}
  By Proposition~\ref{pro:an_quasi-free_extension_universal}, there
  are morphisms of extensions in both directions which lift the
  identity map on~\(A\) and whose composite maps are
  dagger homotopic to the identity maps.
\end{proof}

\begin{proposition}
  \label{pro:extension_an}
  Let \(A \into E \onto B\) be an extension of pro-algebras.  If
  \(A\) and~\(B\) are isomorphic to projective systems of dagger
  algebras, then so is~\(E\).  If \(A\) and~\(B\) are analytically
  nilpotent, then so is~\(E\).
\end{proposition}

\begin{proof}
  Being nilpotent mod~\(\dvgen\) is hereditary for pro-algebra
  extensions by Proposition~\ref{pro:extension_nilpotent_mod_pi}.
  Hence the second statement follows from the first one.  Its proof
  has several steps.  First, we rewrite the given extension of
  pro-algebras as a projective limit of a projective system of
  algebra extensions.  Similar ideas in a less specialised setting
  also appear in~\cite{Artin-Mazur:Etale}*{Appendix}.

  Write \(E\) and~\(B\) as projective systems
  of (torsion-free, complete bornological) algebras
  \((E_n,\gamma_{n,m})\) and \((B_n,\beta_{n,m})\) that are indexed
  by directed sets \(N_E\) and~\(N_B\), respectively.  By
  assumption, \(B\) is isomorphic to a projective system of dagger
  algebras.  We assume that we have picked this representative
  above, that is, each~\(B_n\) is a dagger algebra.  We describe the
  pro-algebra morphism \(E\to B\) by a coherent family of bounded
  homomorphisms \(\varphi_n\colon E_{m(n)} \to B_n\) for all
  \(n\in N_B\).  Let
  \(N\defeq \setgiven{(m,n)\in N_E \times N_B}{m \ge m(n)}\).
  Define a partial order on~\(N\) by \((m_1,n_1) \ge (m_2,n_2)\) if
  \(m_1\ge m_2\), \(n_1\ge n_2\), \(m_1\ge m(n_2)\), and
  \(\beta_{n_2,n_1} \circ \varphi_{n_1} \circ \gamma_{m(n_1),m_1} =
  \varphi_{n_2} \circ \gamma_{m(n_2),m_1}\).  This partially ordered
  set is directed because \(N_B\) and~\(N_E\) are directed and the
  maps~\(\varphi_n\) for \(n\in N\) form a morphism of projective
  systems.  The objects \(E_m\) and~\(B_n\) for \((m,n)\in N\) and
  the maps \(\gamma_{m_1,m_2}\) and \(\beta_{n_1,n_2}\) for
  \(m_1\ge m_2\) and \(n_1\ge n_2\) form projective systems \(E'\)
  and~\(B'\) of bornological algebras.  They are isomorphic to \(E\)
  and~\(B\), respectively.  The homomorphisms
  \[
    \varphi'_{(m,n)} \defeq \varphi_n\circ \gamma_{m(n),m}\colon
    E'_{(m,n)} = E_m \to B_n = B'_{(m,n)}
  \]
  for \((m,n)\in N\) are coherent in the strong sense that
  \[
    \beta'_{(m_1,n_1),(m_2,n_2)} \circ \varphi'_{(m_2,n_2)}
    = \varphi'_{(m_1,n_1)} \circ \gamma'_{(m_1,n_1),(m_2,n_2)}
  \]
  for all \((m_1,n_1),(m_2,n_2)\in N\) with
  \((m_1,n_1) \le (m_2,n_2)\).  Here \(\gamma'\) and~\(\beta'\)
  denote the structure maps of the projective systems \(E'\)
  and~\(B'\), respectively.  By construction, each~\(B'_n\) is a
  dagger algebra.

  By assumption, the inclusion \(A\to E\) is the kernel of the
  morphism \(E \to B\).  This is isomorphic to the kernel of
  \(\varphi'\colon E' \to B'\).  So~\(A\) is isomorphic to the
  projective system~\(A'\) formed by the closed ideals
  \(A'_n \defeq \ker \varphi_n \subseteq E_n'\) for \(n\in N\) with
  the structure maps
  \(\alpha'_{n_1,n_2} = \gamma'_{n_1,n_2}|_{A_{n_2}}\) for
  \(n_1,n_2\in N\) with \(n_1 \le n_2\); and the canonical morphism
  \(A'\to E'\) is the strongly coherent family of inclusion maps
  \(A'_n \hookrightarrow E'_n\) for \(n\in N\).  Each~\(A'_n\) is
  complete and torsion-free because \(E'_n\) and~\(B'_n\) are (see
  \cite{Meyer-Mukherjee:Bornological_tf}*{Theorem~2.3 and
    Lemma~4.2}).

  The quotients \(E'_n/A'_n\) with the structure maps
  \(\dot\gamma'_{n,m}\) induced by~\(\gamma'_{n,m}\) form a
  projective system of complete bornological algebras, which is the
  cokernel for the inclusion \(A' \hookrightarrow E'\).  The
  map~\(\varphi'_n\) for \(n\in N\) descends to an injective,
  bounded homomorphism \(\varrho_n\colon E'_n/A'_n \to B'_n\).  The
  pro-algebra morphism \(\varrho= (\varrho_n)_{n\in N}\) is an
  isomorphism because \(E\to B\) is assumed to be another cokernel
  for the map \(A\to E\).  Next, we modify our projective systems so
  that these become equalities; this replaces the quotients
  \(E'_n/A'_n\) by dagger algebras.  The inverse of~\(\varrho\) is
  given by a choice of \(m(n)\in N\) for \(n\in N\) and bounded
  homomorphisms \(\psi_n\colon B'_{m(n)} \to E'_n/A'_n\).
  Increasing \(m(n)\) if necesessary, we may arrange that
  \(\varrho_n\circ\psi_n = \beta'_{n,m(n)}\colon B'_{m(n)} \to
  B'_n\) and
  \(\psi_n\circ\varrho_{m(n)} = \dot\gamma'_{n,m(n)}\colon
  E'_{m(n)}/A'_{m(n)} \to E'_n/A'_n\).  Let
  \(N' \defeq \setgiven{(m,n)\in N\times N}{m\ge m(n)}\).  For
  \((m,n)\in N'\), pull the extension \(A'_n \into E'_n \onto
  E'_n/A'_n\) back along~\(\psi_n\) as in
  Lemma~\ref{lem:pull-back_push-out_bornological_extension}.  This
  gives a diagram of extensions of bornological \(\dvr\)\nb-modules
  \[
    \begin{tikzcd}
      A''_{(m,n)} \arrow[r,rightarrowtail ] \arrow[d, equal] &
      E''_{(m,n)} \arrow[r, twoheadrightarrow] \arrow[d] &
      B''_{(m,n)} \arrow[d, "\psi_n"] \\
      A'_n \arrow[r,rightarrowtail] &
      E'_n \arrow[r, twoheadrightarrow] &
      E'_n/A'_n
    \end{tikzcd}
  \]
  with \(A''_{(m,n)} = A'_n\) and \(B''_{(m,n)} = B'_m\).  The
  latter is a dagger algebra because it is equal to~\(B_m\) for
  suitable \(m \in N_B\) depending on \(n\in N'\).  There is a
  unique bornological algebra structure on~\(E''_{(m,n)}\) for which
  all maps in this diagram are homomorphisms.  We claim
  that~\(E''_{(m,n)}\) is complete.  First, \(A'_n\) is closed
  in~\(E'_n\) because~\(B'_n\) is separated.  Then \(E'_n/A'_n\) is
  separated (see
  \cite{Meyer-Mukherjee:Bornological_tf}*{Lemma~2.1}).
  Then~\(E''_{(m,n)}\) is closed in \(B'_m \oplus E'_n\).  And
  then~\(E''_{(m,n)}\) is complete by
  \cite{Meyer-Mukherjee:Bornological_tf}*{Theorem~2.3}.  As above,
  there is a partial order on~\(N'\) that makes it a directed set
  and such that \(A''_n \into E''_n \onto B''_n\) becomes a
  projective system of algebra extensions.  This projective system
  is isomorphic to \(A' \into E' \onto E'/A'\) because it is the
  pullback along the pro-algebra isomorphism \(B' \congto E'/A'\).
  Thus it is isomorphic to the original extension
  \(A \into E \onto B\).  We have now replaced this pro-algebra
  extension by a projective system of algebras extensions where the
  quotients~\(B''_n\) are dagger algebras.

  To simplify notation, we remove the primes now and assume that our
  pro-algebra extension already comes to us as a projective system
  of algebra extensions \(A_n \into E_n \onto B_n\), where \(A_n\)
  and~\(E_n\) are torsion-free, complete bornological algebras
  and~\(B_n\) are dagger algebras for all \(n\in N\).  The dagger
  completions~\(E_n^\updagger\) for \(n\in N\) form a projective
  system of dagger algebras, and the canonical maps
  \(E_n \to E_n^\updagger\) form a pro-algebra morphism.  We claim
  that this pro-algebra morphism is an isomorphism.  Equivalently,
  for each \(n\in N\) there are \(m \in N\) with \(m \ge n\) and a
  bounded homomorphism
  \(\tilde\gamma_{n,m}\colon E_m^\updagger \to E_n\) such that the
  composite map \(E_m \to E_m^\updagger \to E_n\)
  is~\(\gamma_{n,m}\); then the other composite map
  \(E_m^\updagger \to E_n \to E_n^\updagger\) is the map on the
  dagger completions induced by~\(\gamma_{n,m}\), and these two
  equalities of compositions say that we are dealing with morphisms
  of pro-algebras inverse to each other.  Fix \(n\in N\).  We are
  going to build the following commuting diagram, where the dashed
  arrow is the desired map~\(\tilde\gamma_{n,m}\):
  \[
    \begin{tikzcd}
      A_m \arrow[rr,rightarrowtail] \arrow[dd, "f"] \arrow[swap, ddd, bend right=30, "\alpha_{n,m}"] &&
      E_m \arrow[rrr, twoheadrightarrow] \arrow[dd] \arrow[dl] \arrow[ddd, bend left=30, "\gamma_{n,m}"]&&&
      B_m \arrow[dd, equal] \\
      & E_m^\updagger \arrow[rd]
      \arrow[swap,rdd, near end, dotted, "\tilde{\gamma}_{n,m}"] &&&&\\
      \tilde{A}_{n'} \arrow[rr,rightarrowtail] \arrow[d, "g"] &&
      \tilde{E}_{n} \arrow[rrr, twoheadrightarrow] \arrow[d]&&&
      B_m \arrow[d, "\beta_{n,m}"]\\
      A_n \arrow[rr,rightarrowtail] &&
      E_{n} \arrow[rrr, twoheadrightarrow]&&&
      B_n
    \end{tikzcd}
  \]
  By assumption, \(A\) is isomorphic to a projective system of
  dagger algebras~\((\tilde{A}_{n'})_{n'\in N'}\).  Therefore, there are
  \(m\in N\), \(n'\in N'\), and maps
  \(f \colon A_m \to \tilde{A}_{n'}\) and
  \(g \colon \tilde{A}_{n'} \to A_n\) such that \(m \ge n\) and
  \(g \circ f = \alpha_{n,m}\colon A_m \to A_n\).
  Let~\(\tilde{E}_n\) be the pushout bornological
  \(\dvr\)\nb-module of the maps \(A_m \to E_m\) and
  \(A_m \to \tilde{A}_{n'}\).  This fits in an extension of
  bornological \(\dvr\)\nb-modules
  \(\tilde{A}_{n'} \into \tilde{E}_n \onto B_m\) by
  Lemma~\ref{lem:pull-back_push-out_bornological_extension}.  Since
  \(\tilde{A}_{n'}\) and~\(B_m\) are torsion-free and complete,
  \(\tilde{E}_n\) is complete by
  \cite{Meyer-Mukherjee:Bornological_tf}*{Theorem~2.3}.
  Since~\(\tilde{A}_{n'}\) is semi-dagger, the canonical map
  \(E_m\to \tilde{E}_n\) remains bounded when we give~\(E_m\) the
  linear growth bornology relative to the ideal~\(A'_m\).  This
  bornology is equal to the absolute linear growth bornology
  on~\(E_m\) by Lemma~\ref{lem:rel=abs} because \(B_m=E_m/A_m\) is a
  dagger algebra.  Since~\(\tilde{E}_n\) is complete, the map
  \(E_m\to \tilde{E}_n\) extends to a bounded \(\dvr\)\nb-module
  homomorphism \(E_m^\updagger \to \tilde{E}_n\).  By construction,
  the map \(\gamma_{n,m}\colon E_m \to E_n\) agrees on~\(A_m\) with
  the composite map
  \[
    A_m \overset{f}\to \tilde{A}_{n'} \overset{g}\to A_n \to
    E_n.
  \]
  Then the universal property of pushouts gives an induced bounded
  \(\dvr\)\nb-module homomorphism
  \(\tau\colon \tilde{E}_n \to E_n\).  Let
  \(\tilde\gamma_{n,m}\colon E_m^\updagger \to E_n\) be the
  composite of the bounded \(\dvr\)\nb-module homomorphisms
  \(E_m^\updagger \to \tilde{E}_n\) and \(\tilde{E}_n \to E_n\)
  defined above.  The composite map
  \(E_m \to E_m^\updagger \to E_n\) is~\(\gamma_{n,m}\) by
  construction.  This finishes the proof that~\(E_n\) is isomorphic
  to a projective system of dagger algebras.
\end{proof}

\subsection{The analytic tensor algebra}
\label{sec:analytic_tensor}

Let~\(R\) be a constant pro-algebra.  The definitions of \(\HAC(R)\)
and \(\HA_*(R)\) use a certain pro-algebra~\(\tans R\) defined by
completing the tensor algebra~\(\tens R\).  We call~\(\tans R\) the
analytic tensor algebra of~\(R\).  We show that there is a
semi-split analytically nilpotent extension
\(\jans R \into \tans R \onto R\) and that~\(\tans R\) is
analytically quasi-free.  Since it is not more difficult, we extend
the construction of the analytic tensor algebra to pro-algebras
right away.

\begin{definition}
  \label{def:tans}
  Let \(R= (R_n,\alpha_{m,n})_{m,n\in N}\) be a pro-algebra.
  Extending the tensor algebra construction to pro-algebras gives a
  natural semi-split pro-algebra extension
  \(\jens R \into \tens R \onto R\) with
  \(\tens R = (\tens R_n)_{n\in N}\) and
  \(\jens R = (\jens R_n)_{n\in N}\).  For each \(n\in N\), we form
  the tube algebras \(\tub{\tens R}{(\jens R)^l}\) with the ideals
  \(\tub{\jens R}{(\jens R)^l}\), and their relative dagger
  completions
  \((\tub{\tens R}{(\jens R)^l}, \tub{\jens R}{(\jens
    R)^l})^\updagger\).  These form a pro-algebra indexed by the
  product set \(N\times \N\), which we call the \emph{analytic
    tensor algebra} of~\(R\) and denote by \(\tans R\).
\end{definition}

\begin{lemma}
  \label{lem:tans_properties}
  The canonical homomorphism \(p\colon \tens R \to R\) extends
  uniquely to a pro-algebra homomorphism
  \(\tilde{p}\colon \tans R \to R\).  The composite~\(\sigma_{\an}\)
  of the pro-linear map \(\sigma_R\colon R \to \tens R\) and the
  canonical homomorphism \(\tens R \to \tans R\) is a section
  for~\(\tilde{p}\).
\end{lemma}

\begin{proof}
  Fix \(n\in N\) and \(l\in\N^*\).  The canonical homomorphism
  \(\tens R_n \to R_n\) vanishes on~\(\jens R_n\).  Then it extends
  uniquely to the tube algebra \(\tub{\tens R_n}{(\jens R_n)^l}\) by
  Proposition~\ref{pro:make-pro-nilpotent}.  This extension vanishes
  on \(\tub{\jens R_n}{(\jens R_n)^l}\).  Then it remains bounded
  for the linear growth bornology relative to this ideal and extends
  uniquely to a homomorphism on the relative dagger completion.
  These maps for all \(n\) and~\(l\) form a morphism of pro-algebras
  \(\tilde{p}\colon \tans R \to R\).  The canonical maps
  \(\sigma_{R_n}\colon R_n \to \tens R_n\) form a pro-linear section
  for \(p\colon \tens R \to R\).  Composing with the canonical map
  \(\tens R \to \tans R\) gives a section for~\(\tilde{p}\).
\end{proof}

\begin{definition}
  \label{def:jans}
  Let~\(\jans R\) be the kernel of
  \(\tilde{p}\colon \tans R \onto R\).
\end{definition}

Lemma~\ref{lem:tans_properties} implies that there is a semi-split
extension of pro-algebras
\[
  \begin{tikzcd}
    \jans R \arrow[r, rightarrowtail] &
    \tans R \arrow[r, twoheadrightarrow, "\tilde{p}"] &
    R.  \arrow[l, dashed, bend left, "\sigma_R"]
  \end{tikzcd}
\]

\begin{proposition}
  \label{pro:jans_nilpotent}
  The pro-algebra~\(\jans R\) is analytically nilpotent.
\end{proposition}

\begin{proof}
  Let \(m\in\N^*\).  The linear growth bornology on
  \(\tub{\tens R}{(\jens R)^m}\) relative to
  \(\tub{\jens R}{(\jens R)^m}\) restricts to the ``absolute''
  linear growth bornology on \(\tub{\jens R}{(\jens R)^m}\) by
  Lemma~\ref{lem:relling}.  The tensor algebra is bornologically
  torsion-free by Remark~\ref{rem:tensor_tf}.  Then so is
  \(\tub{\tens R}{(\jens R)^m}\) by the definition of the bornology
  on the tube algebra.  Then the relative linear growth bornology on
  it is torsion-free by Lemma~\ref{lem:relgtf}, and this property is
  preserved by completions (see
  \cite{Meyer-Mukherjee:Bornological_tf}*{Theorem~4.6}).  Therefore,
  the completion of \(\tub{\jens R}{(\jens R)^m}\) in the linear
  growth bornology is a dagger algebra.  Then \(\jans R\) is a
  pro-dagger algebra.  And \(\tub{\jens R}{(\jens R)^\infty}\) is
  nilpotent mod~\(\dvgen\) by
  Proposition~\ref{pro:make-pro-nilpotent}.  This remains unaffected
  when we equip the tube algebras with the linear growth bornology
  and complete.
\end{proof}

\begin{remark}
  \label{rem:tans_dagger}
  Let \(R= (R_n,\alpha_{m,n})_{m,n\in N}\) be a projective system of
  dagger algebras.  Since
  \(\tub{\tens R}{(\jens R)^l} \bigm/ \tub{\jens R}{(\jens R)^l}
  \cong R\) is semi-dagger, the linear growth bornology on
  \(\tub{\tens R}{(\jens R)^l}\) is equal to the linear growth
  bornology relative to \(\tub{\jens R}{(\jens R)^l}\) by
  Lemma~\ref{lem:rel=abs}.  Hence~\(\tans R\) is also equal to the
  ``absolute'' dagger completion,
  \[
    \tans R \cong  \tub{\tens R}{(\jens R)^\infty}^\updagger.
  \]
\end{remark}

\begin{proposition}
  \label{pro:tans_quasi-free}
  The analytic tensor algebra~\(\tans R\) is analytically quasi-free
  and quasi-free.  The bimodule \(\comb{\Omega}{}^1(\tans R)\) is isomorphic
  to the free bimodule on~\(R\), that is,
  \begin{equation}
    \label{eq:Omega1_tansD}
    (\tans R)^+ \hot R \hot (\tans R)^+ \cong \comb{\Omega}{}^1(\tans R);
  \end{equation}
  the isomorphism is the map
  \(\omega \otimes x\otimes \eta \mapsto \omega\cdot (\diff
  \sigma_R(x))\cdot \eta\).  And the following maps are isomorphisms
  of left or right \(\tans R\)\nb-modules, respectively:
  \begin{alignat*}{2}
    (\tans R)^+ \hot R &\congto \tans R,&\qquad
    \omega \otimes x\mapsto \omega \Fed \sigma_R(x),\\
    R \hot (\tans R)^+ &\congto \tans R,&\qquad
    x \otimes \omega\mapsto \sigma_R(x) \Fed \omega.
  \end{alignat*}
\end{proposition}

\begin{proof}
  Let \(J \into E \overset{q}\onto \tans R\) be a semi-split,
  analytically nilpotent pro-algebra extension.  Pull it back along
  the inclusion \(\jans R \hookrightarrow \tans R\) to a pro-algebra
  extension \(J \into K \onto \jans R\) and identify~\(K\) with an
  ideal in~\(E\).  Since \(J\) and~\(\jans R\) are analytically
  nilpotent, so is~\(K\) by Proposition~\ref{pro:extension_an}.  Let
  \(s\colon \tans R\to E\) be a pro-linear section and let
  \(\sigma_R\colon R \to \tans R\) be the canonical pro-linear
  section.  The pro-linear map \(s\circ \sigma_R\) induces a
  pro-algebra homomorphism
  \((s\circ \sigma_R)^\#\colon \tens R \to E\) by
  Lemma~\ref{lem:tensor_algebra}.  It satisfies
  \(q\circ (s\circ \sigma_R)^\# = \sigma_R^\#\), and
  \(\sigma_R^\#\colon \tens R \to \tans R\) is the canonical
  homomorphism because \(\sigma_R^\#\) and the inclusion map agree
  on the image of~\(R\) in~\(\tens R\).  In particular,
  \((s\circ \sigma_R)^\#\) maps \(\jens R\) into \(K\idealin E\).
  Since~\(K\) is nilpotent mod~\(\dvgen\),
  Proposition~\ref{pro:make-pro-nilpotent} shows that
  \((s\circ \sigma_R)^\#\) extends to the tube algebra
  \(\tub{\tens R}{(\jens R)^\infty}\), in such a way that
  \(\tub{\jens R}{(\jens R)^\infty}\) is mapped to~\(K\).  And
  since~\(K\) is a pro-dagger algebra, the criterion in
  Proposition~\ref{pro:relative_dagger} shows that the morphism
  \(\tub{\tens R}{(\jens R)^\infty} \to E\) extends uniquely to the
  dagger completion relative to \(\tub{\jens R}{(\jens R)^\infty}\).
  This gives a pro-algebra morphism \(\tans R \to E\) that is a
  section for the extension \(J \into E \overset{q}\onto \tans R\).
  So~\(\tans R\) is analytically quasi-free.

  If \(h\colon R\to E\) is any pro-linear map with
  \(q\circ h = \sigma_R\), then the argument above shows that
  \(h^\#\colon \tens R \to E\) extends uniquely to a pro-algebra
  morphism \(\tans R \to E\) that is a section for the extension.
  Conversely, any multiplicative section \(g\colon \tans R \to E\)
  is of this form for \(h \defeq g\circ \sigma_R\).  Thus the
  multiplicative sections for the extension
  \(J \into E \overset{q}\onto \tans R\) are in bijection with
  pro-linear maps \(R\to E\) with \(q\circ h = \sigma_R\).  Any such
  pro-linear map is equal to \(s\circ \sigma_R+ h_0\) for a unique
  pro-linear map \(h_0\colon R\to J\).  So multiplicative sections
  for our extension are in bijection with pro-linear maps
  \(R\to J\).  Combined with Lemma~\ref{lem:Omega1_vs_sections}, we
  get a natural bijection for all \(\tans R\)-bimodules~\(M\)
  between pro-bimodule homomorphisms \(\comb{\Omega}{}^1(\tans R) \to M\) and
  pro-linear maps \(R\to M\).  Thus \(\comb{\Omega}{}^1(\tans R)\) is
  isomorphic to the free bimodule on~\(R\), which is
  \((\tans R)^+ \hot R \hot (\tans R)^+\).  And this isomorphism is
  indeed induced by the map
  \(\omega \otimes x\otimes \eta \mapsto \omega\cdot (\diff
  \sigma_R(x))\cdot \eta\).

  Now let~\(M\) be a left \(\tans R\)\nb-module.  Turn~\(M\) into a
  \(\tans R\)\nb-bimodule by taking the zero map as right module
  structure.  Then a bimodule derivation \(\tans R \to M\) is just a
  left module map.  Therefore, left module homomorphisms
  \(\tans R \to M\) are in bijection with pro-linear maps
  \(R \to M\).  Thus the map
  \[
    (\tans R)^+ \hot R,\qquad
    \omega \otimes x\mapsto \omega \Fed \sigma_R(x),
  \]
  is an isomorphism of left \(\tans R\)\nb-modules.  Here we have
  written~\(\Fed\) for the multiplication in~\(\tans R\) because we
  will later use these formulas when~\(\tans R\) is identified
  with~\(\Omega^\even R\) with the Fedosov product.  A similar
  argument works for right modules.
\end{proof}

We now describe the analytic tensor algebra and its bornology more
concretely.  For this, we assume that~\(R\) is a torsion-free,
complete bornological algebra.  A projective system
\((R_n)_{n\in N}\) is treated by applying the following discussion
to~\(R_n\) for each \(n\in N\).  We identify~\(\tens R\)
with~\(\Omega^\even R\) with the Fedosov product as in
Section~\ref{sec:tensor_nc-forms}.  Recall that the isomorphism
\(\tens R \cong \Omega^\even R\) maps the ideal~\(\jens R^m\) onto
\(\bigoplus_{n\ge m} \Omega^{2 n} R\).  Thus
\(\tub{\tens R}{(\jens R)^m}\) is spanned by
\(\dvgen^{-j} \Omega^{2 n} R\) with \(n \ge m\cdot j\).  And
\(\tub{\jens R}{(\jens R)^m}\) is spanned by
\(\dvgen^{-j} \Omega^{2 n} R\) with \(n \ge m\cdot j\) and
\(n\ge1\).  Equivalently,
\begin{equation}
  \label{eq:tub_tens}
  \tub{\tens R}{(\jens R)^m}
  = \sum_{n=0}^\infty \dvgen^{-\floor{n/m}} \Omega^{2 n} R,
  \qquad
  \tub{\jens R}{(\jens R)^m}
  = \sum_{n=1}^\infty \dvgen^{-\floor{n/m}} \Omega^{2 n} R.
\end{equation}

The following lemma estimates the growth of Fedosov products
in~\(\Omega R\):

\begin{lemma}
  \label{lem:fedo}
  Let~\(R\) be an algebra and let \(M\subseteq R\) be a submodule.
  Let \(i_0,\dotsc,i_n\ge 1\) and \(i\defeq i_0+\dotsb+i_n\).  Then
  \[
    \Omega^{i_0} M \Fed\dotsb \Fed \Omega^{i_n}M
    \subseteq \bigoplus_{j=0}^n\Omega^{i+2j}(M^{(3)}).
  \]
\end{lemma}

\begin{proof}
  As in the proof of \cite{Meyer:HLHA}*{Theorem~5.11}, we show the
  more precise estimate
  \begin{equation}
    \label{inc:fedo}
    \Omega^{i_0}M \Fed \dotsb \Fed \Omega^{i_n} M
    \subseteq \bigoplus_{j=0}^n {}\AU{(M^{(2)})}
    \,\diff(M^{(3)})^{i+2 j}
  \end{equation}
  by induction on~\(n\).  This is trivial for \(n=0\).  The
  induction step uses~\eqref{prod:fedosov2} and
  \[
    \Omega^iM \Fed \AU{(M^{(2)})}
    \subseteq \AU{(M^{(2)})} \diff (M^{(3)})^i
    + (\diff M)^{i+1}\, \diff (M^{(2)}).\qedhere
  \]
\end{proof}

\begin{proposition}
  \label{pro:anborno}
  Let~\(R\) be a torsion-free bornological algebra and \(m\ge 1\).
  If \(M\subseteq R\) is bounded, \(\alpha\in \Q\cap (0,1/m)\), and
  \(f\in\N_0\), then define
  \begin{equation}
    \label{bded:anborno}
    D_m(M,\alpha,f) \defeq \bigoplus_{n=0}^\infty
    \dvgen^{-\floor*{\min \{n/m, \alpha\cdot n + f\}}}
    \Omega^{2 n}M.
  \end{equation}
  These are \(\dvr\)\nb-submodules of \(\tub{\tens R}{(\jens R)^m}\)
  that cofinally generate its linear growth bornology relative to
  the ideal \(\tub{\jens R}{(\jens R)^m}\).
\end{proposition}

\begin{proof}
  Let \(M\subseteq R\) be bounded, \(\alpha\in \Q\cap (0,1/m)\), and
  \(f\in\N_0\).  Equation~\eqref{eq:tub_tens} implies
  \(D_m(M,\alpha,f)\subseteq \tub{\tens R}{\jens R^m}\).  Our first
  goal is to show that \(D_m(M,\alpha,f)\) has linear growth
  relative to \(\tub{\jens R}{\jens R^m}\).  Let \(N\subseteq R\) be
  a submodule and \(e\ge 1\).  We claim that
  \begin{equation}
    \label{bded:anaborno1}
    \AU{N}\cdot
    \bigg(\sum_{n=1}^{e m} \dvgen^{-\floor*{n/m}}
    (\diff N\,\diff N)^n \biggr)^\diamond
    = \bigoplus_{n=1}^\infty \dvgen^{-\floor*{n/m}
      + \ceil*{\frac{n}{e m}}-1} \Omega^{2 n}N.
  \end{equation}
  By definition, the left hand side is spanned by Fedosov products
  \begin{multline*}
    \dvgen^{j-1- \floor{i_1/m} - \dotsb - \floor{i_j/m}} \AU{N}\Fed
    (\diff N,\diff N)^{i_1}\Fed \dotsb \Fed (\diff N,\diff N)^{i_j}
    \\= \dvgen^{j-1- \floor{i_0/m} - \dotsb - \floor{i_j/m}}
    \Omega^{2(i_1 + \dotsb + i_j)}(N)
  \end{multline*}
  for \(j\ge 1\) and \(1 \le i_1,\dotsc,i_j \le e m\).  These
  contribute to \(\Omega^{2 n} N\) if \(i_1 + \dotsb + i_j = n\).
  For fixed \(n\) and~\(j\), the sum of floors
  \(\floor{i_1/m} + \dotsb + \floor{i_j/m}\) is maximal if all but
  one of the~\(i_j\) are divisible by~\(m\), and then it
  becomes~\(\floor{n/m}\).  For fixed~\(n\), the term
  \(j-1-\floor{n/m}\) becomes minimal if~\(j\) is minimal.
  Equivalently, we choose \(i_j = e m\) for all but one~\(j\), and
  then \(j = \ceil{n/e m}\).  This finishes the proof
  of~\eqref{bded:anaborno1}.

  The right hand side in~\eqref{bded:anaborno1} is one of the
  generators of the linear growth bornology relative to
  \(\tub{\jens R}{(\jens R)^m}\).  For fixed \(\alpha<1/m\)
  and~\(f\) as above, there is \(e\in\N^*\) with
  \(1/m - 1/(e m)>\alpha\).  Then there is \(k\in\N\) with
  \[
    \floor*{n/m} - \ceil*{\frac{n}{e m}} + 1
    \ge \floor*{\min \{n/m, \alpha\cdot n + f\}}
  \]
  for \(n>k\).  Then
  \[
    D_m(M,\alpha,f)
    \subseteq \sum_{n=0}^k
    \dvgen^{-\floor*{\min \{n/m, \alpha\cdot n + f\}}}
    \Omega^{2 n}M
    + \AU{N}\cdot
    \bigg(\sum_{n=1}^{e m} \dvgen^{-\floor*{n/m}} (\diff
    N\,\diff N)^n \biggr)^\diamond.
  \]
  The first, finite sum is already bounded in
  \(\tub{\tens R}{\jens R^m}\).  Therefore, \(D_m(M,\alpha,f)\)
  has linear growth relative to \(\tub{\jens R}{(\jens R)^m}\).




  Now let~\(S\) be any \(\dvr\)\nb-submodule of
  \(\tub{\tens R}{\jens R^m}\) that has linear growth relative to
  \(\tub{\jens R}{(\jens R)^m}\).  We claim that~\(S\) is contained
  in \(D_m(M,\alpha,f)\) for suitable \(M,\alpha,f\).  By definition
  of the relative linear growth bornology, there are \(k, e\in\N\)
  and a bounded submodule \(M\subseteq R\) such that~\(S\) is
  contained in the sum of
  \(\sum_{n=0}^k \dvgen^{-\floor*{n/m}} \Omega^{2n}M\) and
  \(\bigl(\sum_{i=1}^{e m} \dvgen^{-\floor*{\frac{i}{m}}} \Omega^{2
    i} M\bigr)^\diamond\).  The latter is spanned by Fedosov
  products
  \[
    \dvgen^{j-1 -\floor*{\frac{i_1}{m}} - \dotsb -
      \floor*{\frac{i_j}{m}}}
    \Omega^{2 i_1} M \Fed \dotsb \Fed \Omega^{2 i_j} M
  \]
  with \(j\in\N^*\), \(1\le i_1,\dotsc, i_j \le e m\).  By
  Lemma~\ref{lem:fedo},
  \(\Omega^{2 i_1} M \Fed \dotsb \Omega^{2 i_j} M\) is contained in
  the sum of \(\Omega^{2 n} (M^{(3)})\), where~\(n\) lies between
  \(i\defeq \sum_{k=1}^j i_j\) and \(i+j\).  As above, the sum of
  the floors \(\floor{i_j/m}\) for fixed~\(i\) is maximal if all but
  one~\(i_j\) are divisible by~\(m\), and then it
  is~\(\floor{i/m}\).  The constraints \(i_k \le e m\) are
  equivalent to the constraint \(i \le j\cdot e m\).  So~\(S\) is
  contained in the sum of
  \(\dvgen^{j-1- \floor{i/m}} \Omega^{2 n}(M^{(3)})\) with
  \(i \le n \le i+j\) and \(i \le j\cdot e m\).  For fixed \(n,j\),
  the exponent \(j-1 - \floor{i/m}\) is minimal if~\(i\) is maximal,
  so we may assume that~\(i\) is the minimum of \(n\) and~\(j e m\).
  Then the optimal choice for~\(j\) is the minimal one, which is
  \(\ceil{n/(e m)}\) if \(i=n\) and \(j = \ceil{n/(e m+1)}\) if
  \(i=j e m\).  The resulting exponents of~\(\dvgen\) become
  \(\ceil{n/(e m)} - 1 - \floor{n/m}\) in the first case and
  \(\ceil{n/(e m+1)} - 1 - \ceil{n/(e m+1)}\cdot e\) in the second.
  If \(\alpha > 1/m - 1/(e m)\) and~\(n\) is large enough, then both
  terms are greater or equal \(-\floor{\alpha n}\).  Choosing~\(f\)
  big enough, we may arrange that both are greater or equal
  \(-\floor{\min \{n/m,\alpha n+f\}}\) for all \(n\in \N\).  Then
  \(S \subseteq D_m(M^{(3)},\alpha,f)\).
\end{proof}

\begin{corollary}
  \label{cor:anborno}
  For \(m\in\N^*\), let~\(\mathcal{B}_m\) be the bornology on
  \(\tub{\tens R}{\jens R^m}\) that contains a subset if and only if
  it is contained in
  \(\bigoplus_{n=0}^\infty \dvgen^{-\floor*{\frac{n}{m}}}
  \Omega^{2 n} M\) for some bounded \(\dvr\)\nb-submodule
  \(M\subseteq R\).  This bornology makes
  \(\tub{\tens R}{\jens R^m}\) a torsion-free bornological algebra.
  The projective system of bornological algebras
  \((\tub{\tens R}{\jens R^m},\mathcal{B}_m)_{m\in\N^*}\) is
  isomorphic to the projective system formed by
  \(\tub{\tens R}{\jens R^m}\) with the linear growth bornology
  relative to \(\tub{\jens R}{\jens R^m}\).
\end{corollary}

\begin{proof}
  By Lemma~\ref{lem:fedo}, the Fedosov product is bounded for the
  bornology~\(\mathcal{B}_m\).  The subsets \(D_m(M,\alpha,f)\)
  in~\eqref{bded:anborno} are clearly in~\(\mathcal{B}_m\).
  Conversely,
  \[
    \bigoplus_{n=0}^\infty \dvgen^{-\floor*{\frac{n}{m+1}}}
    \Omega^{2 n} M
    = D_m\bigl(M,\textstyle\frac1{m+1},0\bigr).
  \]
  Thus any subset in~\(\mathcal{B}_{m+1}\) is mapped to a subset of
  \(\tub{\tens R}{\jens R^m}\) with linear growth relative to
  \(\tub{\jens R}{\jens R^m}\).  The asserted isomorphism of
  projective systems follows.
\end{proof}

Now we can describe the completion~\(\tans R\).  Recall
that~\(\comb{\Omega}^n R\) denotes the completion \(R^+ \hot R^{\hot n}\)
of \(\Omega^n R = R^+ \otimes R^{\otimes n}\).  For \(m\in\N^*\) and
a bounded \(\dvr\)\nb-submodule \(M\subseteq R\), the canonical map
\(\comb{\Omega}^{2 n} M \to \comb{\Omega}^{2 n} R\) is injective by
Proposition~\ref{pro:tensor_injective}.  Then we may view
\(\prod_{n=0}^\infty \dvgen^{-\floor*{\frac{n}{m}}} \comb{\Omega}^{2
  n} M\) as a \(\dvr\)\nb-submodule of
\(\prod_{n=0}^\infty \comb{\Omega}^{2 n} R \otimes \dvf\).  Let
\(\comb{\Omega}^\ev(R)_m\) be the union of
\(\prod_{n=0}^\infty \dvgen^{-\floor*{\frac{n}{m}}} \comb{\Omega}^{2
  n} M\) for all bounded \(\dvr\)\nb-submodules \(M\subseteq R\),
with the bornology where a subset is bounded if and only if it is
contained in
\(\prod_{n=0}^\infty \dvgen^{-\floor*{\frac{n}{m}}} \comb{\Omega}^{2
  n} M\) for some bounded \(\dvr\)\nb-submodules \(M\subseteq R\).
These form a decreasing sequence of subalgebras with bounded
inclusion maps
\(\comb{\Omega}^\ev(R)_{m+1} \hookrightarrow \comb{\Omega}^\ev(R)_m\).

\begin{proposition}
  \label{pro:completed_analytic_tensor}
  If~\(R\) is a torsion-free, complete bornological algebra,
  then~\(\tans R\) is naturally isomorphic to the projective system
  of complete bornological algebras
  \((\comb{\Omega}^\ev(R)_m)_{m\in\N^*}\).
\end{proposition}

\begin{proof}
  We shall use the explicit description of the relative linear
  growth bornology in Proposition~\ref{pro:anborno}.  Each
  \(\dvgen^{-\floor{n/m}} \Omega^{2 n}R\) is a direct summand of
  \(\tub{\tens R}{\jens R^m}\), and the projection is bounded in the
  linear growth bornology relative to \(\tub{\jens R}{\jens R^m}\).
  This gives us maps from the completed tube to
  \(\dvgen^{-\floor{n/m}} \comb{\Omega}^{2 n}R\) for all \(n\in\N\).
  It is easy to see that the \(\dvgen\)\nb-adic completion of
  \(D_m(M,\alpha,f)\) is isomorphic to the subspace of
  \(\prod_{n=0}^\infty \dvgen^{-\floor{n/m}} \comb{\Omega}^{2 n} M\)
  consisting of all \((\omega_n)_{n\in\N}\) for which there is a
  sequence \((h_j)_{j\in\N}\) in~\(\N\) with \(\lim h_j = \infty\)
  and
  \(\omega_n \in \dvgen^{-\floor*{\min \{n/m, \alpha\cdot n + f\}} +
    h_n} \comb{\Omega}^{2 n}M\) for all \(n\in\N\).  Any such subset
  is bounded in \(\comb{\Omega}^\ev(R)_m\).  Conversely, any bounded subset
  in \(\comb{\Omega}^\ev(R)_{m+1} \) is contained in a subset of this form
  with \(f=0\) and \(m< 1/\alpha <m+1\).  Therefore, the projective
  system formed by the relative dagger completions
  \(\bigl(\tub{\tens R}{\jens R^m}, \tub{\jens R}{\jens
    R^m}\bigr)^\updagger\) is isomorphic to the projective system
  \((\comb{\Omega}^\ev(R)_m)_{m\in\N^*}\).
\end{proof}

\subsection{Pro-Linear maps with nilpotent curvature}
\label{sec:nilpotent_curvature}

Let \(R\) and~\(S\) be pro-algebras.  We are going to describe
pro-algebra homomorphisms \(\tans R \to S\) through a certain class
of pro-linear maps \(R\to S\), namely, those with analytically
nilpotent curvature.  This follows rather easily from the concrete
description of the relative linear growth bornology above.  The main
issue is to define analytically nilpotent curvature.  We begin with
the analogue of nilpotent curvature mod~\(\dvgen\).

\begin{definition}
  Let \(X= (X_{n'})_{n'\in N'}\) be a bornological pro-module and
  \(S = (S_n)_{n\in N}\) a pro-algebra and let
  \(\omega\colon X \to S\) be a pro-linear map.  We call~\(\omega\)
  \emph{nilpotent mod~\(\dvgen\)} if, for every \(n\in N\), there is
  \(m \in \N^*\) such that the following composite is zero:
  \begin{equation}
    \label{map:pronilpo}
    X^{\otimes m} \xrightarrow{\omega^{\otimes m}}
    S_m^{\otimes m} \xrightarrow{\mathrm{mult}}
    S_m \to S_m/\dvgen S_m;
  \end{equation}
  here \(\mathrm{mult}\) denotes the \(m\)\nb-fold multiplication
  map of~\(S\).
\end{definition}

Let \(\omega\colon X\to S\) be nilpotent mod~\(\dvgen\) and
represent~\(\omega\) by a coherent family of bounded
\(\dvr\)\nb-module maps \(\omega_n\colon X_{r(n)} \to S_n\) with
\(r(n) \in N'\) for \(n\in N\).  For \(n\in N\) and \(n' \in N'\)
with \(n' \ge r(n)\), let \(\omega_{n,n'}\colon X_{n'} \to S_n\) be
the composite map \(X_{n'} \to X_{r(n)} \to S_n\).  Let \(n\in N\)
and choose~\(m\) so that the map in~\eqref{map:pronilpo} vanishes.
Then there is \(n' \in N'\) with \(n' \ge r(n)\) such that the
composite map
\(X_{n'}^{\otimes m} \to S_n^{\otimes m}\to S_n \to S_n /\dvgen
S_n\) vanishes.  That is,
\(\omega_{n,n'}(x_1) \dotsm \omega_{n,n'}(x_m) \in \dvgen\cdot S_n\)
for all \(x_1,\dotsc,x_m \in X_{n'}\).  Let \(M\subseteq X_{n'}\) be
bounded.  Since~\(\omega_{n,n'}\) is bounded and~\(S_n\) is
torsion-free, it follows that
\(\omega_{n,n'}(M)^m \subseteq \dvgen S_n\) and that
\(\dvgen^{-1}\cdot \omega_{n,n'}(M)^m \subseteq S_n\) is bounded.
Then
\begin{equation}
  \label{bded:lgcurva}
  \omega_{n,n'}(M)_e \defeq
  \sum_{j=1}^{e m} \dvgen^{-\floor{j/m}} \omega_{n,n'}(M)^j
\end{equation}
is bounded for every \(e\ge 1\).

\begin{definition}
  Let \(X= (X_m)_{m\in N'}\) be a bornological pro-module and
  \(S = (S_n)_{n\in N}\) a pro-algebra and let
  \(\omega\colon X \to S\) be a pro-linear map.
  Represent~\(\omega\) by a coherent family of bounded
  \(\dvr\)\nb-module maps \(\omega_{n,n'}\colon X_{n'} \to S_n\) as
  above.  The map~\(\omega\) is called \emph{analytically nilpotent}
  if, for every~\(n\), there are \(m\in \N^*\) and \(n'\in N'\) with
  \(n' \ge r(n)\) such that for any bounded subset
  \(M\subseteq X_{n'}\), the subset
  \[
    \sum_{j=0}^\infty \dvgen^{-\floor*{j/m}} \omega_{n,n'}(M)^j
    \subseteq S_n \otimes \dvf
  \]
  is bounded in~\(S_n\).
\end{definition}

\begin{proposition}
  \label{pro:puanatube}
  Let \(R\) and~\(S\) be pro-algebras and let \(f\colon R\to S\) be
  a pro-linear map.  Let \(\omega\colon R\otimes R\to S\),
  \(x\otimes y\mapsto f(x\cdot y) - f(x)\cdot f(y)\), be its
  curvature.  There is a pro-algebra homomorphism
  \(f^\#\colon \tans R \to S\) with \(f=f^\# \sigma_R = f\) if and
  only if~\(\omega\) is analytically nilpotent.
\end{proposition}

\begin{proof}
  Write \(R = (R_{n'})_{n' \in N'}\) and \(S = (S_n)_{n \in N}\) as
  projective systems of algebras.  Then~\(\tans R\) is the
  completion of the projective system of bornological algebras
  \(T\defeq (\tub{\tens R_{n'}}{\jens R_{n'}^m},
  \mathcal{B}_m)_{n'\in N', m\in\N^*}\) with the
  bornologies~\(\mathcal{B}_m\) in Corollary~\ref{cor:anborno}.
  Since~\(S\) is complete, any homomorphism of projective systems of
  bornological algebras \(T\to S\) extends uniquely to~\(\tans R\).
  Since~\(S\) is torsion-free, such a homomorphism \(T\to S\) is
  determined by its restriction to~\(\tens R\).  Then there is a
  unique pro-linear map \(f\colon R\to S\) such that the
  homomorphism is \(f^\#\colon \tens R \to S\) as
  in~\eqref{eq:fsharp_curvature}.  Corollary~\ref{cor:anborno} shows
  that~\(f^\#\) extends to a homomorphism \(T\to S\) if and only
  if~\(f\) has analytically nilpotent curvature.
\end{proof}

\begin{corollary}
  \label{coro:composmallcurvature}
  Let \(f\colon R\to S\), \(g\colon S\to T\) be pro-linear maps and
  let \(U\) be a projective system of dagger algebras.  If \(f\)
  and~\(g\) have analytically nilpotent curvature, then so do
  \(g\circ f\) and \(f\hot U\colon R\hot U \to S\hot U\).
\end{corollary}

\begin{proof}
  The assertion about \(g\circ f\) follows as in the proof of
  \cite{Meyer:HLHA}*{Theorem~5.23}, using
  \cite{Meyer-Mukherjee:Bornological_tf}*{Theorems 3.7 and~4.5}.
  Since~\(f\) has analytically nilpotent curvature, there is a
  homomorphism \(f^\#\colon \tans R \to S\) with
  \(f^\#\circ\sigma_R = f\).  The extension
  \[
    (\jans R) \hot U \into (\tans R)\hot U \onto R\hot U
  \]
  is analytically nilpotent because \((\jans R) \hot U\) is
  nilpotent mod~\(\dvgen\) by
  Proposition~\ref{pro:extension_nilpotent_mod_pi} and a pro-dagger
  algebra by the extension of Corollary~\ref{cor:tensor_dagger} to
  projective systems.  The pro-linear section \(\sigma_R \hot U\)
  induces a homomorphism \(\tans(R\hot U) \to (\tans R)\hot U\)
  which, when composed with \(f^\#\), gives a homomorphism
  \(\tans (R\hot U) \to S\) that extends~\(f\hot U\).
  Thus~\(f\hot U\) has analytically nilpotent curvature.
\end{proof}


\subsection{Homotopy invariance of the X-complex}
\label{sec:homotopy_X}

In this section, we assume that the field~\(\dvf\) has
characteristic~\(0\).  This is needed to prove that homotopic
homomorphisms defined on a quasi-free algebra induce chain homotopic
maps between the \(X\)\nb-complexes.  If we understand homotopy to
mean ``polynomial homotopy'', then this is already shown by Cuntz
and Quillen (see
\cite{Cuntz-Quillen:Cyclic_nonsingularity}*{Sections 7--8}).  In the
context of complete bornological \(\dvr\)\nb-algebras, the proof for
polynomial homotopies still works for dagger homotopies.  The
corresponding statement for the \(B,b\)-bicomplexes is
\cite{Cortinas-Cuntz-Meyer-Tamme:Nonarchimedean}*{Proposition~4.3.3}.
For quasi-free algebras, the canonical projection from the
\(B,b\)-bicomplex to the \(X\)\nb-complex is a chain homotopy
equivalence.  This implies the following:

\begin{proposition}
  \label{pro:X-homotopy}
  Let \(R\) and~\(S\) be projective systems of complete bornological
  \(\dvf\)\nb-\hspace{0pt}algebras.  Let
  \(f,g\colon R \rightrightarrows S\) be two homomorphisms that are
  dagger homotopic.  Assume that~\(\dvf\) has characteristic~\(0\)
  and that~\(R\) is quasi-free.  Then the induced chain maps
  \(X(f),X(g)\colon X(R) \rightrightarrows X(S)\) are chain
  homotopic.
\end{proposition}

\begin{proof}
  It suffices to treat an elementary dagger homotopy.  Define
  \begin{align*}
    \eta_n \colon \Omega^n(S \otimes \dvr[t])\otimes \dvf
    &\to \Omega^{n-1}(S) \otimes \dvf,\\
    a_0\,\diff a_1 \dotsc \diff a_n
    &\mapsto \int_0^1 a_0(t) \frac{\partial a_1(t)}{\partial t}
      \diff a_2(t) \dotsc \diff a_n(t) \,\diff t,
  \end{align*}
  for \(n = 1, 2\).  Here integration and differentiation are
  defined formally by rescaling the coefficients of polynomials
  \(a_i \in S \otimes \dvf[t]\).  We claim that~\(\eta_n\) extends
  to a bounded linear map
  \(\eta_n \colon \Omega^n(S \otimes \dvr[t]^\updagger ) \otimes
  \dvf \to \Omega^{n-1}(S )\otimes \dvf\).  To see this, let
  \(T \defeq S \otimes \ling{\dvr[t]}\).  Then
  \(\Omega^n(T) \cong T^+ \otimes T^{\otimes n} \cong T^{\otimes
    n}\oplus T^{\otimes n+1}\).  So it suffices to show
  that~\(\eta_n\) is bounded on
  \(T^{\otimes n} \otimes \dvf \cong S^{\otimes n} \otimes
  \ling{\dvr[t]}^{\otimes n} \otimes \dvf\).  This follows if the map
  \begin{align*}
    \ling{\dvr[t]}^{\otimes n} \otimes \dvf
    &\to \dvf,\\
    a_0 \otimes a_1 \otimes \dotsc \otimes a_n
    &\mapsto \int_0^1 a_0(t) \frac{\partial a_1(t)}{\partial t}
      \cdot a_2(t) \dotsm a_n(t)\,\diff t
  \end{align*}
  is bounded.  The formal differentiation on~\(\ling{\dvr[t]}\) is
  clearly bounded.  And~\(\ling{\dvr[t]}\) is a bornological
  algebra.  So this happens if and only if the integration map
  \[
    \ling{\dvr[t]} \otimes \dvf
    \to \dvf,\qquad
    a(t) = \sum_{l=0}^\infty c_l
    t^l \mapsto \sum_{l=0}^\infty \frac{c_l}{l+1}
  \]
  is bounded.  If~\(\resf\) has characteristic~\(0\), then \(l+1\)
  is invertible in~\(\dvr\) for all \(l\in\N\).  If~\(\resf\) has
  finite characteristic~\(p\), then the valuation of \(l+1\) grows
  at most logarithmically.  In any case, this is dominated by the
  linear growth of the exponents of~\(\dvgen\) for a subset of
  linear growth in~\(\dvr[t]\).  Thus the integration map above is
  bounded, and then so are the maps~\(\eta_n\).  We still
  write~\(\eta_n\) for their unique bounded extensions to the
  completions.
  
  Let \(\eta_0 = 0\).  Then \([\eta,b] = 0\).  Therefore,
  \(\eta_2(b(\comb{\Omega}^3(S \hot \dvr[t]^\updagger))) \subseteq
  b(\comb{\Omega}^2(S))\).  So~\(\eta\) defines a map
  \(X^{(2)}(S \hot \dvr[t]^\updagger) \to X(S)\), where~\(X^{(2)}\)
  is the truncated \(B-b\)-complex defined in
  \cite{Meyer:HLHA}*{Definition A.122}.

  Let
  \(\xi_2 \colon X^{(2)}(S \hot \dvr[t]^\updagger) \to X(S \hot
  \dvr[t]^\updagger)\) be the canonical projection.  Then
  \[
    [\eta, B+b] = (X(\ev_1) - X(\ev_0))\circ \xi \colon X^{(2)}(S
    \hot \dvr[t]^\updagger) \to X(S).
  \]
  Now let \(H \colon R \to S \otimes \dvr[t]^\updagger\) be an
  elementary dagger homotopy between \(f\) and~\(g\).  Then
  \(\eta \circ X^{(2)}(H) \colon X^{(2)}(R) \to X(S)\) is a chain
  homotopy between \(X(f)\circ \xi_2\) and \(X(g) \circ \xi_2\),
  where \(\xi_2 \colon X^{(2)}(R) \to X(R)\) is the canonical
  projection.  Since~\(R\) is analytically quasi-free, it is in
  particular quasi-free, so that \(\xi_2\) is a chain homotopy
  equivalence.  Let \(\alpha \colon X(R) \to X^{(2)}(R)\) be the
  homotopy inverse of \(\xi_2\).  Then \(\eta \circ \alpha\) is the
  desired chain homotopy between \(X(f)\) and~\(X(g)\).
\end{proof}

\begin{theorem}
  \label{the:homotopy_invariance_X}
  Let \(A\) and~\(B\) be pro-algebras.  If two homomorphisms
  \(f_0,f_1\colon A \rightrightarrows B\) are dagger homotopic, then
  they induce homotopic chain maps \(\HAC(A) \to \HAC(B)\).  And
  then \(\HA_*(f_0) = \HA_*(f_1)\).
\end{theorem}

\begin{proof}
  The homomorphisms
  \(\tans f_0,\tans f_1\colon \tans A \rightrightarrows \tans B\)
  lift \(f_0\) and~\(f_1\).  Since~\(\tans A\) is analytically
  quasi-free and~\(\jans B\) is analytically nilpotent,
  Proposition~\ref{pro:an_quasi-free_extension_universal} provides a
  dagger homotopy between \(\tans f_0\) and~\(\tans f_1\).  Then the
  chain maps
  \(X(\tans A\otimes\dvf) \rightrightarrows X(\tans B\otimes\dvf)\)
  induced by \(f_0\) and~\(f_1\) are chain homotopic by
  Proposition~\ref{pro:X-homotopy}.  This remains so on the homotopy
  projective limits.  And then \(f_0\) and~\(f_1\) induce the same
  map on the homology of the homotopy projective limits.  That is,
  \(\HA_*(f_0)=\HA_*(f_1)\).
\end{proof}

\subsection{Invariance under analytically nilpotent extensions}
\label{sec:homotopy_invariance_applications}

We continue to assume that~\(\dvf\) has characteristic~\(0\).

\begin{theorem}
  \label{Goodwillie}
  Let \(J \into E \overset{p}\onto A\) be a semi-split, analytically
  nilpotent extension of pro-algebras.  Then~\(p\) induces a chain
  homotopy equivalence \(\HAC(E) \simeq \HAC(A)\) and \(\HAC(J)\) is
  contractible.  So \(\HA_*(E) \cong \HA_*(A)\) and \(\HA_*(J)=0\).
  If~\(E\) is analytically quasi-free, then \(\HAC(A)\) is chain
  homotopy equivalent to \(X(E \otimes \dvf)\) and \(\HA_*(A)\) is
  isomorphic to the homology of the homotopy projective limit of
  \(X(E \otimes \dvf)\).
\end{theorem}

\begin{proof}
  The composite map \(\tans E \onto E \onto A\) is pro-algebra
  homomorphism with a pro-linear section.  Its kernel~\(K\) is an
  extension of \(\jans E\) by~\(J\) and hence analytically nilpotent
  by Proposition~\ref{pro:extension_an}.  Both \(\tans E\)
  and~\(\tans A\) are analytically quasi-free by
  Proposition~\ref{pro:tans_quasi-free}.
  Proposition~\ref{pro:an_quasi-free_extension_universal} applied to
  the extensions \(K \into \tans E \onto A\) and
  \(\jens A \into \tans A \onto A\) shows that \(\tans A\)
  and~\(\tans E\) are dagger homotopy equivalent.  This together
  with Proposition~\ref{pro:X-homotopy} implies that
  \(\HAC(A) = X(\tans A\otimes \dvf)\) and
  \(\HAC(E) = X(\tans E\otimes \dvf)\) are homotopy equivalent.
  This remains so for their homotopy projective limits.
  So \(\HA_*(E) \cong \HA_*(A)\).  More
  precisely, the isomorphism is the map induced by the quotient map
  \(E\onto A\).

  Since \(J\) and~\(\jans J\) are analytically nilpotent, so
  is~\(\tans J\) by Proposition~\ref{pro:extension_an}.
  Since~\(\tans J\) is analytically quasi-free,
  Proposition~\ref{pro:an_quasi-free_extension_universal} may be
  applied to the extensions \(\tans J = \tans J \to 0\) and
  \(0 = 0 = 0\) of~\(0\).  Thus~\(\tans J\) is dagger homotopy
  equivalent to~\(0\).  Then \(\HAC(J)\simeq0\) and
  \(\HA_*(J) \cong 0\).

  Now assume~\(E\) to be analytically quasi-free.  Then
  Proposition~\ref{pro:an_quasi-free_extension_universal} shows that
  the extensions of~\(A\) by \(\tans A\) and~\(E\) are
  dagger homotopy equivalent.  Then
  \(X(E)\otimes \dvf\) is homotopy equivalent to
  \(X(\tans A)\otimes \dvf\).  Then \(\HAC(A)\) is homotopy
  equivalent to the homotopy projective limit of the projective
  system of chain complexes \(X(E)\otimes \dvf\).
\end{proof}

\begin{corollary}
  \label{cor:analytically_quasi_free_computation}
  Let~\(A\) be an analytically quasi-free algebra.  Then \(\HAC(A)\)
  is chain homotopy equivalent to \(X(A \otimes \dvf)\) and
  \(\HA_*(A)\) is isomorphic to the homology of
  \(X(A \otimes \dvf)\).
\end{corollary}

\begin{proof}
  Theorem~\ref{Goodwillie} shows that \(\HAC(A)\) is homotopy
  equivalent to \(X(A \otimes \dvf)\).  Then \(\HA_*(A)\) is
  isomorphic to the homology of \(\holim X(A \otimes \dvf)\).  Since
  \(X(A \otimes \dvf)\) is a constant projective system, it is chain
  homotopy equivalent to its homotopy projective limit.  So we
  simply get the ordinary homology of \(X(A \otimes \dvf)\).
\end{proof}

\begin{corollary}
  \label{cor:HA_V}
  \(\HAC(\dvr)\) is homotopy equivalent to~\(\dvr\) with zero boundary
  map.
\end{corollary}

\begin{proof}
  The algebra~\(\dvr\) is analytically quasi-free by
  Proposition~\ref{pro:V_analytically_quasi_free}.  Then
  \(\HAC(\dvr) \simeq X(\dvr)\) by
  Corollary~\ref{cor:analytically_quasi_free_computation}.  A small
  calculation shows that any element of \(\Omega^1(\dvr)\) is a
  commutator.  So \(X(\dvr)\) is~\(\dvr\) with zero boundary map.
\end{proof}

\section{Excision}
\label{sec:excision}
\numberwithin{equation}{section}

The goal of this section is to prove the following excision theorem
for analytic cyclic homology:

\begin{theorem}
  \label{the:excision}
  Let \(K \overset{i}\into E \overset{p}\onto Q\) be a semi-split
  extension of pro-algebras with a pro-linear section
  \(s \colon Q \to E\).  Then there is a natural exact triangle
  \[
    \HAC(K) \xrightarrow{i_*}
    \HAC(E) \xrightarrow{p_*}
    \HAC(Q) \xrightarrow{\delta}
    \HAC(K)[-1]
  \]
  in the homotopy category of chain complexes of projective systems
  of bornological \(\dvr\)\nb-modules.  Thus there is a natural long
  exact sequence
  \[
    \begin{tikzcd}
      \HA_0(K) \arrow[r, "i_*"] &
      \HA_0(E) \arrow[r, "p_*"] &
      \HA_0(Q) \arrow[d, "\delta"] \\
      \HA_1(Q) \arrow[u, "\delta"] &
      \HA_1(E) \arrow[l, "p_*"] &
      \HA_1(K).  \arrow[l, "i_*"]
    \end{tikzcd}
  \]
\end{theorem}

The proof will take up the rest of this section.  It follows
\cites{Meyer:Excision, Meyer:HLHA}.  We use the left
ideal~\(\mathcal{L}\) in~\(\tans E\) generated by~\(K\) and prove
chain homotopy equivalences \(X(\tans K) \simeq X(\mathcal{L})\) and
\(X(\mathcal{L}) \simeq X(\tans E : \tans Q)\) as chain complexes in
the additive category of projective systems of bornological
\(\dvr\)\nb-modules.  First, the pro-linear section~\(s\) yields two
bounded maps
\(s_L,s_R\colon \Omega^\even Q\rightrightarrows \Omega^\even E\)
defined by
\begin{align*}
  s_L(q_0\,\diff q_1\dotsc \diff q_{2n})
  &\defeq s(q_0)\,\diff s(q_1)\dotsc \,\diff s(q_{2n}),\\
  s_R(\diff q_1\dotsc \diff q_{2n}\, q_{2n+1})
  &\defeq \diff s(q_1)\dotsc \diff s(q_{2n})\, s(q_{2n+1})
\end{align*}
for all \(q_0,q_{2n+1}\in Q^+\) and \(q_i\in Q\) for
\(1 \le i \le 2 n\).  Let \(m\in\N^*\).  Both \(s_L\) and~\(s_R\)
map~\(\jens Q^{m j}\) to~\(\jens E^{m j}\) for all \(j\in\N\)
by~\eqref{eq:tub_tens}.  Thus
they induce bounded linear maps on the tubes, from
\(\tub{\tens Q}{\jens Q^m}\) to \(\tub{\tens E}{\jens E^m}\).  Both
are sections for the canonical projection
\(\tub{\tens E}{\jens E^m} \to \tub{\tens Q}{\jens Q^m}\).  These
sections remain bounded for the linear growth bornologies relative
to \(\tub{\jens E}{\jens E^m}\) and \(\tub{\jens Q}{\jens Q^m}\) by
Proposition~\ref{pro:anborno}.  Thus they extend to bounded
\(\dvr\)\nb-module maps on the completions.  These maps for all
\(m\in\N^*\) form two pro-linear sections for
\(\tans p \colon \tans E \to \tans Q\).  They induce two sections
for the canonical chain map
\(X(\tans p) \colon X(\tans E) \to X(\tans Q)\).  Let
\[
  X(\tans E : \tans Q) \defeq
  \ker\bigl(X(\tans p) \colon X(\tans E) \to X(\tans Q)\bigr).
\]
There is a semi-split extension of chain complexes
\[
  X(\tans E: \tans Q) \into\tans E \to \tans Q.
\]
Since \(X(\tans p) \circ X(\tans i) = X(\tans (p\circ i)) = 0\), the
chain map \(X(\tans i)\) factors through \(X(\tans E : \tans Q)\).
We are going to prove that this chain map
\(X(\tans K) \to X(\tans E : \tans Q)\) is a chain homotopy
equivalence.  Then the homotopy projective limit of \(X(\tans K)\)
is homotopy equivalent to that of \(X(\tans E : \tans Q)\), and the
latter fits into a semi-split extension of chain complexes with the
homotopy projective limits of \(X(\tans E)\) and \(X(\tans Q)\).  As
a result, Theorem~\ref{the:excision} follows if the inclusion map
\(X(\tans K) \to X(\tans E : \tans Q)\) is a chain homotopy
equivalence.

Our construction of the chain homotopy equivalence will, in
principle, be explicit and natural, using only the multiplication
maps in our pro-algebras and the pro-linear sections \(s_L\)
and~\(s_R\) above.  Therefore, we assume for simplicity from now on
that we are dealing with an extension of (complete, torsion-free
bornological) algebras \(K \into E \onto Q\).  In general, we may
rewrite the semi-split extension above as a projective system of
semi-split algebra extensions \(K_n \into E_n \onto Q_n\) with
compatible bounded linear sections; this uses arguments as in the
proof of Proposition~\ref{pro:extension_an}.  To simplify notation,
we write down the proof below only for a semi-split algebra
extension.  The chain maps and homotopies that we are going to build
for the extensions \(K_n \into E_n \onto Q_n\) form morphisms of
projective systems.  So the same proof works for a semi-split
extension of pro-algebras.

\subsection{The pro-algebra \texorpdfstring{$\mathcal{L}$}{L}}
\numberwithin{equation}{subsection}

In the following, we identify \(\tens E\) with \(\Omega^\even E\)
and~\(E\) with \(\Omega^0(E) \subseteq \Omega^\even E\).  So the map
\(\sigma_E\colon E \to \tans E\) disappears from our notation.
Proposition~\ref{pro:tans_quasi-free} gives an isomorphism of left
\(\tans E\)\nb-modules
\begin{equation}
  \label{eq:tansEplus_E_is_tansE}
  (\tans E)^+ \hot E \congto \tans E,\qquad
  \omega \otimes x\mapsto \omega \Fed x.
\end{equation}
Explicitly, the inverse of this isomorphism is given by
\begin{equation}
  \label{eq:tansE_decompose_inverse}
  \omega \,\diff e_{2n-1}\,\diff e_{2n} \mapsto
  \omega \otimes (e_{2n-1}\cdot e_{2n})
  - (\omega \Fed e_{2n-1}) \otimes e_{2n}.
\end{equation}
These two maps also define an isomorphism for the purely algebraic
tensor algebras:
\begin{equation}
  \label{eq:tensEplus_E_is_tensE}
  (\tens E)^+\otimes E\congto \tens E,\qquad
  \omega \otimes e\mapsto \omega \Fed e.
\end{equation}
Variants of this isomorphism and the following ones were proven
already in \cite{Meyer:HLHA}*{Section~4.3.2}.  Let
\(L\subseteq \tens E\) be the left ideal generated by~\(K\).  The
bounded linear section \(s\colon Q\to E\) yields an isomorphism of
bornological \(\dvr\)\nb-modules \(E \cong K \oplus Q\).
Then~\eqref{eq:tensEplus_E_is_tensE} implies an isomorphism
\begin{equation}
  \label{eq:tensEplus_K_is_L}
  (\tens E)^+\otimes K\congto L,\qquad
  \omega \otimes k\mapsto \omega \Fed k.
\end{equation}
The explicit formula for the isomorphism
in~\eqref{eq:tansE_decompose_inverse} and its inverse imply
\begin{equation*}
  \label{L:concrete}
  L = K \oplus \bigoplus_{n\ge 1} \Omega^{2n-1}(E)\,\diff K
\end{equation*}
as in the proof of \cite{Meyer:HLHA}*{Lemma~4.55}.  Let
\(I\defeq \ker(\tens p\colon \tens E \onto \tens Q)\).  This is part
of semi-split extensions
\begin{equation}
  \label{eq:I_tensE_tensQ_extension}
  \begin{tikzcd}
    I \arrow[r, rightarrowtail] &
    \tens E \arrow[r, twoheadrightarrow, "\tens p"] &
    \tens Q \arrow[l, dashed, bend left, "s_L"]
  \end{tikzcd}\qquad
  \begin{tikzcd}
    I \arrow[r, rightarrowtail] &
    (\tens E)^+ \arrow[r, twoheadrightarrow, "(\tens p)^+"] &
    (\tens Q)^+.  \arrow[l, dashed, bend left, "s_L"]
  \end{tikzcd}
\end{equation}

\begin{lemma}
  \label{lem:isomorphisms_tensE_L_I}
  The following maps are isomorphisms:
  \begin{alignat}{2}
    \label{eq:Lplus_tensQplus_is_tensEplus}
    \Psi\colon L^+ \otimes (\tens Q)^+&\congto (\tens E)^+,\qquad&
    l\otimes \eta &\mapsto l\Fed s_L(\eta),\\
    \label{eq:L_tensQplus_is_I}
    L \otimes (\tens Q)^+&\congto I,\qquad&
    l\otimes \eta &\mapsto l\Fed s_L(\eta),\\
    \label{eq:tensEplus_K_tensQplus_is_I}
    (\tens E)^+ \otimes K \otimes (\tens Q)^+&\congto I,\qquad&
    \omega\otimes k\otimes \eta &\mapsto \omega\Fed k\Fed s_L(\eta),\\
    \label{eq:tensQplus_K_tensEplus_is_I}
    (\tens Q)^+ \otimes K \otimes (\tens E)^+&\congto I,\qquad&
    \eta\otimes k\otimes \omega &\mapsto s_R(\eta) \Fed k\Fed \omega,\\
    \label{eq:tensQplus_K_Lplus_is_L}
    (\tens Q)^+ \otimes K \otimes L^+&\congto L,\qquad&
    \eta \otimes k \otimes l&\mapsto s_R(\eta) \Fed k\Fed l.
  \end{alignat}
\end{lemma}

\begin{proof}
  The computations in \cite{Meyer:HLHA}*{Section 4.3.1} show this.
  We briefly sketch them.  The isomorphisms
  \eqref{eq:Lplus_tensQplus_is_tensEplus}
  and~\eqref{eq:L_tensQplus_is_I} are equivalent because of the
  semi-split extension~\eqref{eq:I_tensE_tensQ_extension}.  And
  \eqref{eq:L_tensQplus_is_I}
  and~\eqref{eq:tensEplus_K_tensQplus_is_I} are equivalent because
  of the isomorphism~\eqref{eq:tensEplus_K_is_L}.  The isomorphisms
  \eqref{eq:tensEplus_K_tensQplus_is_I}
  and~\eqref{eq:tensQplus_K_tensEplus_is_I} imply each other by
  taking opposite algebras because this reverses the order of
  multiplication and exchanges \(s_L\) and~\(s_R\).
  And~\eqref{eq:tensQplus_K_tensEplus_is_I}
  implies~\eqref{eq:tensQplus_K_Lplus_is_L} by substituting
  \((\tens E)^+ \cong L^+ \otimes (\tens Q)^+\) and
  \(I \cong L \otimes (\tens Q)^+\)
  in~\eqref{eq:tensQplus_K_Lplus_is_L} and then cancelling the
  factor \((\tens Q)^+\) on both sides.

  So it suffices to prove that~\(\Psi\) is an isomorphism.  We
  describe its inverse~\(\Psi^{-1}\).  Split a differential form
  \(e_0 \,\diff e_1 \dotsc \diff e_{2 n}\in \Omega^{2 n}E\) so that
  each coefficient \(e_j\) belongs either to \(K\) or~\(s(Q)\), or
  is~\(1\) in case of~\(e_0\); this is possible because of the
  direct sum decomposition \(E\cong K \oplus s(Q)\); write
  \(k_i \defeq e_i\) or \(q_i \defeq s^{-1}(e_i)\) accordingly.  If
  no~\(e_i\) belongs to~\(K\), then
  \[
    \Psi^{-1}\bigl( s(q_0)\,\diff s(q_1)\dotsc \diff s(q_{2n})\bigr)
    = 1 \otimes q_0\,\diff q_1\dotsc \diff q_{2n}.
  \]
  Otherwise, there is a largest~\(i\le 2 n\) with \(e_i \in K\).
  If \(i=0\), then
  \[
    \Psi^{-1}\bigl( k_0\,\diff s(q_1)\dotsc \diff s(q_{2 n})\bigr)
    = k_0 \otimes \diff q_1\dotsc \diff q_{2 n}.
  \]
  If~\(i\) is even and non-zero, then
  \[
    \Psi^{-1}\bigl( e_0\,\diff e_1\dotsc \diff e_{i-1}\,
    \diff k_i \,\diff s(q_{i+1})\dotsc \diff s(q_{2 n})\bigr)
    = e_0\,\diff e_1\dotsc \diff e_{i-1}\, \diff k_i \otimes
    \diff q_{i+1}\dotsc \diff q_{2 n}.
  \]
  If~\(i\) is odd, then
  \begin{multline*}
    \Psi^{-1}\bigl( e_0\,\diff e_1\dotsc \diff e_{i-1}\,
    \diff k_i \,\diff s(q_{i+1})\dotsc \diff s(q_{2 n})\bigr)
    \\
    \begin{aligned}
      &= e_0\,\diff e_1\dotsc \diff e_{i-1} \Fed (k_i\cdot
      s(q_{i+1})) \otimes \diff q_{i+2}\dotsc \,\diff q_{2 n}
      \\&\qquad- e_0\,\diff e_1\dotsc \diff e_{i-1} \Fed k_i
      \otimes q_{i+1} \, \diff q_{i+2}\dotsc \diff q_{2 n}.
    \end{aligned}
  \end{multline*}
  A direct computation using
  \(\diff k_i\,\diff s(q_{i+1}) = k_i\cdot s(q_{i+1}) - k_i\Fed
  s(q_{i+1})\) shows that
  \[
    \Psi\circ\Psi^{-1}(e_0\,\diff e_1\dotsc \,\diff e_{2n}) =
    e_0\,\diff e_1\dotsc \,\diff e_{2n}
  \]
  for all \(e_0\in \{1\}\cup K\cup s(Q)\),
  \(e_1,\dotsc,e_n\in K\cup s(Q)\).  Then one shows that the
  map~\(\Psi^{-1}\) is surjective: its image contains all elements
  of the form \(1\otimes \eta\) for \(\eta\in (\tens Q)^+\) and
  \(\omega \otimes \diff q_1\dotsc \diff q_{2 n}\) with
  \(\omega \in L^+\) by the first cases where there is no~\(i\)
  or~\(i\) is even, respectively.  And modulo a term of this form,
  the image of~\(\Psi^{-1}\) contains all
  \(\omega \odot k \otimes q_0\,\diff q_1\dotsc \diff q_{2 n}\) with
  \(\omega \in (\tens E)^+\), \(k\in K\) because of the formula in
  the case where~\(i\) is odd.  This exhausts
  \(L^+ \otimes (\tens Q)^+\) because of the
  isomorphism~\eqref{eq:tensEplus_K_is_L}.
\end{proof}

We are going to pass to the analytic tensor algebras and describe
``analytic'' analogues of \(L,I\subseteq \tens E\) and of the
isomorphisms and semi-split extensions above.  For \(m\in\N^*\), let
\begin{align*}
  I_{(m)} &\defeq
            \ker \bigl(\tub{\tens E}{\jens E^m} \to
            \tub{\tens Q}{\jens Q^m}\bigr),\\
  L_{(m)} &\defeq K \oplus
            \bigoplus_{n\ge 1} \dvgen^{-\floor*{n/m}}\cdot
            \Omega^{2n-1}(E)\,\diff K.
\end{align*}
It is easy to see that \(I_{(m)}\) is a two-sided and~\(L_{(m)}\) a
left ideal in \(\tub{\tens E}{\jens E^\infty}\).  In particular,
both are \(\dvr\)\nb-algebras in their own right.  Inspection shows
that
\begin{equation}
  \label{eq:ILm_from_tube}
  I_{(m)} = \tub{\tens E}{\jens E^m} \cap (I \otimes \dvf),\qquad
  L_{(m)} = \tub{\tens E}{\jens E^m} \cap (L \otimes \dvf)
\end{equation}
as \(\dvr\)\nb-submodules of \(\tens E \otimes \dvf\).  The maps in
the projective system \(\tub{\tens E}{\jens E^\infty}\) make
\((I_{(m)})_{m\in\N^*}\) and \((L_{(m)})_{m\in\N^*}\) projective
systems by restriction.  We equip each \(\tub{\tens E}{\jens E^m}\)
with the bornology~\(\mathcal{B}_m\) described in
Corollary~\ref{cor:anborno}; using the linear growth bornology
instead would slightly complicate the estimates below.  We give
\(I_{(m)}\) and~\(L_{(m)}\) the subspace bornologies.  So the
bornology on~\(L_{(m)}\) is cofinally generated by
\begin{equation}
  \label{eq:bornology_Lm}
  (M\cap K) \oplus \bigoplus_{n=1}^\infty
  \dvgen^{-\floor*{n/m}} \Omega^{2 n-1}M \,\diff (M\cap K)
\end{equation}
for bounded \(\dvr\)\nb-submodules \(M\subseteq E\).  Let
\(\mathcal{I} \defeq \bigl(\comb{I_{(m)}}\bigr)_{m\in\N^*}\) and
\(\mathcal{L} \defeq \bigl(\comb{L_{(m)}}\bigr)_{m\in\N^*}\) be the
projective systems formed by the completions.

Since \(\tub{\tens E}{\jens E^m}\) is a subalgebra of
\(\tens E \otimes \dvf\) and the maps in
\eqref{eq:tensEplus_E_is_tensE}, \eqref{eq:tensEplus_K_is_L} and
\eqref{eq:Lplus_tensQplus_is_tensEplus}--\eqref{eq:tensQplus_K_Lplus_is_L}
only involve Fedosov products and the maps \(s_L\) and~\(s_R\),
\eqref{eq:ILm_from_tube} implies that these maps
still exist and are bounded if \(\tens E, \tens Q,I, L\) are
replaced by
\(\tub{\tens E}{\jens E^m},\tub{\tens Q}{\jens Q^m},I_{(m)},
L_{(m)}\), respectively, each equipped with the relative linear
growth bornologies specified above.  The inverse maps for these
isomorphisms are slightly more complicated, however: they may shift
the index~\(m\) in the projective system:

\begin{lemma}
  \label{lem:inverses_in_excision_proof}
  The inverses to the isomorphisms above extend to bounded maps
  \begin{align*}
    \tub{\tens E}{\jens E^{m+1}}
    &\to \tub{\tens E}{\jens E^m}^+\otimes E,\\
    L_{(m+1)}
    &\to \tub{\tens E}{\jens E^m}^+\otimes K,\\
    \tub{\tens E}{\jens E^{2 m}}^+
    &\to L_{(m)}^+ \otimes \tub{\tens Q}{\jens Q^m}^+,\\
    I_{(2 m)}
    &\to L_{(m)} \otimes \tub{\tens Q}{\jens Q^m}^+,\\
    I_{(2 m)}
    &\to \tub{\tens E}{\jens E^m}^+ \otimes K \otimes \tub{\tens Q}{\jens Q^m}^+,\\
    I_{(2 m)}
    &\to \tub{\tens Q}{\jens Q^m}^+ \otimes K \otimes \tub{\tens E}{\jens E^m}^+,\\
    L_{(2 m)}
    &\to \tub{\tens Q}{\jens Q^m}^+ \otimes K \otimes L_{(m)}^+.
  \end{align*}
\end{lemma}

\begin{proof}
  Our explicit formula for the first map shows that it reduces the
  total degree of a differential form by at most~\(2\).  Since
  \(\frac{n+1}{m+1} \le \frac{n}{m}\) for all \(n\ge m\) and
  \(\floor*{\frac{n+1}{m+1}} = \floor*{\frac{n}{m}} = 0\) if
  \(n<m\), it follows that it defines a map
  \(\tub{\tens E}{\jens E^{m+1}} \to \tub{\tens E}{\jens E^m}^+
  \otimes E\) that is bounded for the bornologies described in
  Corollary~\ref{cor:anborno}.  The second map is a restriction of
  the first map, so that it is covered by the same argument.

  Our explicit formula for the third map shows that it maps a
  differential form of degree~\(2 n\) to a sum of tensor products
  involving differential forms of degree~\(2 j\) and \(2(n-j-1)\) or
  \(2(n-j)\); in the first case, \(j<n\) and the differential form
  in~\(L\) is already explicitly written as \(\omega \odot k\), so
  that the isomorphism \(L \to (\tens E)^+ \otimes K\) does not
  reduce the degree any further.  This shows that the same degree
  estimate applies to the fourth map in the lemma.  The fifth map
  differs from that only by taking opposite algebras, and the sixth
  map is a restriction of the fifth one.  This is why the following
  estimates cover all these maps at the same time.

  That these maps are well defined between the relevant tube
  algebras amounts to the estimate
  \(\floor{n/2 m} \le \floor{j/m} + \floor{(n-j-1)/m}\) for all
  \(n\in\N\), \(0 \le j< n\).  This is trivial for \(n<2 m\), so
  that we assume \(n\ge 2 m\).  For fixed~\(n\), the right hand side
  is minimal if \(j=m-1\), and then the needed estimate simplifies
  to \(\floor{n/2 m} \le \floor{(n-m)/m}\).  This is true for
  \(2 m \le n<4 m\).  Since adding~\(2 m\) to~\(n\) increases
  \(\floor{n/2 m}\) by~\(1\) and \(\floor{(n-m)/m}\) by~\(2\), the
  inequality follows for all \(n\in\N\).  Now it follows that the
  maps in the lemma are well defined and bounded for the bornologies
  described in Corollary~\ref{cor:anborno}.
\end{proof}

The composite maps
\begin{gather*}
  \tub{\tens E}{\jens E^{m+1}}^+\otimes E
  \to \tub{\tens E}{\jens E^{m+1}}
  \to \tub{\tens E}{\jens E^m}^+\otimes E,\\
  \tub{\tens E}{\jens E^{m+1}}
  \to \tub{\tens E}{\jens E^m}^+\otimes E
  \to \tub{\tens E}{\jens E^m}^+
\end{gather*}
are the structure maps in our projective systems because they extend
the identity maps on \((\tens E)^+\otimes E\) and \(\tens E\),
respectively.  Thus these two families of maps for \(m\in\N^*\) are
isomorphisms of projective systems of bornological
\(\dvr\)\nb-modules that are inverse to each other.  This remains so
when we complete, giving an isomorphism
\((\tans E)^+\hot E\congto \tans E\).  The same argument applies to
the other isomorphisms above.  Summing up, we get the following
isomorphisms of projective systems of bornological
\(\dvr\)\nb-modules:
\begin{alignat}{2}
  \label{eq:tansEplus_E_is_tansE2}
  (\tans E)^+\hot E&\congto \tans E,\qquad&
  \omega \otimes e&\mapsto \omega \Fed e,\\
  \label{eq:tansEplus_K_is_L}
  (\tans E)^+\hot K&\congto \mathcal{L},\qquad&
  \omega \otimes k&\mapsto \omega \Fed k,\\
  \label{eq:Lplus_tansQplus_is_tansEplus}
  \mathcal{L}^+ \hot (\tans Q)^+&\congto (\tans E)^+,\qquad&
  l\otimes \eta &\mapsto l\Fed s_\mathcal{L}(\eta),\\
  \label{eq:L_tansQplus_is_I}
  \mathcal{L} \hot (\tans Q)^+&\congto \mathcal{I},\qquad&
  l\otimes \eta &\mapsto l\Fed s_L(\eta),\\
  \label{eq:tansEplus_K_tansQplus_is_I}
  (\tans E)^+ \hot K \hot (\tans Q)^+&\congto \mathcal{I},\qquad&
  \omega\otimes k\otimes \eta &\mapsto \omega\Fed k\Fed s_L(\eta),\\
  \label{eq:tansQplus_K_tansEplus_is_I}
  (\tans Q)^+ \hot K \hot (\tans E)^+&\congto \mathcal{I},\qquad&
  \eta\otimes k\otimes \omega &\mapsto s_R(\eta) \Fed k\Fed \omega,\\
  \label{eq:tansQplus_K_Lplus_is_L}
  (\tans Q)^+ \hot K \hot \mathcal{L}^+&\congto \mathcal{L},\qquad&
  \eta \otimes k \otimes l&\mapsto s_R(\eta) \Fed k\Fed l.
\end{alignat}
In addition, there are semi-split extensions
\begin{equation}
  \label{eq:I_tansE_tansQ_extension}
  \begin{tikzcd}
    \mathcal{I} \arrow[r, rightarrowtail] &
    \tans E \arrow[r, twoheadrightarrow, "\tans p"] &
    \tans Q \arrow[l, dashed, bend left, "s_L"]
  \end{tikzcd}\qquad
  \begin{tikzcd}
    \mathcal{I} \arrow[r, rightarrowtail] &
    (\tans E)^+ \arrow[r, twoheadrightarrow, "(\tans p)^+"] &
    (\tans Q)^+.  \arrow[l, dashed, bend left, "s_L"]
  \end{tikzcd}
\end{equation}
Here~\eqref{eq:tansEplus_E_is_tansE2} is the same
as~\eqref{eq:tansEplus_E_is_tansE}.  So it follows already from the
analytic nilpotence machinery in Section~\ref{sec:aqf}.
And~\eqref{eq:tansEplus_E_is_tansE2} easily
implies~\eqref{eq:tansEplus_K_is_L}.  The isomorphisms
\eqref{eq:L_tansQplus_is_I}--\eqref{eq:tansQplus_K_Lplus_is_L}
follow from
\eqref{eq:tansEplus_E_is_tansE2}--\eqref{eq:Lplus_tansQplus_is_tansEplus}
and the semi-split extension~\eqref{eq:I_tansE_tansQ_extension} as
in the proof of Lemma~\ref{lem:isomorphisms_tensE_L_I}.  It seems
that the existence of the maps
\(s_L, s_R\colon \tans Q \rightrightarrows \tans E\) and
\eqref{eq:Lplus_tansQplus_is_tansEplus} do not follow from the
machinery in Section~\ref{sec:aqf} and must be checked by hand.

\begin{theorem}
  \label{the:L_to_relative}
  The chain map \(X(\mathcal{L}) \to X(\tans E: \tans Q)\) induced
  by the inclusion \(\mathcal{L} \hookrightarrow \tans E\) is a
  chain homotopy equivalence.
\end{theorem}

\begin{proof}
  The proofs of \cite{Meyer:HLHA}*{Theorems 4.66 and~4.67} carry
  over literallly to our analytic tensor algebras, using the
  isomorphisms
  \eqref{eq:tansEplus_E_is_tansE2}--\eqref{eq:tansQplus_K_Lplus_is_L}
  and the semi-split extension~\eqref{eq:I_tansE_tansQ_extension}.
  We merely have to replace the symbols \(\otimes\),
  \(A\defeq \overleftarrow{\tens} E\), \(\overleftarrow{\tens} Q\),
  \(\overleftarrow{L}\), \(\overleftarrow{I}\)
  and~\(\overleftarrow{G}\) in that proof by \(\hot\), \(\tans E\),
  \(\tans Q\), \(\mathcal{L}\), \(\mathcal{I}\) and
  \((\tans Q)^+ \hot K\), respectively; and
  \(\overleftarrow{\Omega}^\even E\) and
  \(\overleftarrow{\Omega}^\odd E\) in~\cite{Meyer:HLHA} become
  \(\tans E\) and \((\tans E)^+ \hot E\), respectively, with the
  latter identified with differential forms of odd degree.
  \cite{Meyer:HLHA}*{Theorem~5.80} is a similar translation exercise
  for the analytic cyclic homology theory for bornological algebras
  over the complex numbers, and the situation in this article is
  exactly the same.

  We briefly sketch the main idea of the proof.
  Proposition~\ref{pro:tans_quasi-free} and the definition of
  \(\comb{\Omega}{}^1(\tans E)\) imply that there is a semi-split
  free \(\tans E\)\nb-bimodule resolution
  \[
    \comb{\Omega}{}^1(\tans E) \into (\tans E)^+\hot (\tans E)^+
    \onto (\tans E)^+
  \]
  with a natural pro-linear section
  \((\tans E)^+ \to (\tans E)^+\hot (\tans E)^+\),
  \(x\mapsto 1 \otimes x\).  Let
  \begin{align*}
    P_0 &\defeq \mathcal{L}^+ \hot  \mathcal{L}^+
          + (\tans E)^+ \hot \mathcal{L} \subseteq
          (\tans E)^+ \hot (\tans E)^+,\\
    P_1 &\defeq (\tans E)^+ \, D \mathcal{L} \subseteq
          \comb{\Omega}{}^1 (\tans E)^+.
  \end{align*}
  This together with \(\mathcal{L}^+ \subseteq (\tans E)^+\) gives a
  subcomplex of the resolution above, and the standard section above
  yields a contracting homotopy for it, making it a resolution.  The
  bimodules \(P_0\) and~\(P_1\) are free; this is where the
  isomorphisms above enter.  So
  \(P_1 \into P_0 \onto \mathcal{L}^+\) is a free
  \(\mathcal{L}\)\nb-bimodule resolution.  Then~\(\mathcal{L}\) is
  quasi-free, and the \(X\)\nb-complex computes its periodic cyclic
  homology.  And the commutator quotient complex
  \(P_1/[\mathcal{L},P_1] \to P_0/[\mathcal{L},P_0]\) computes the
  Hochschild homology of~\(\mathcal{L}\).  These commutator
  quotients are computed explicitly and shown to compute the
  relative Hochschild homology for the quotient map
  \(\tans E \onto \tans Q\).  And then the isomorphism on Hochschild
  homology implies an isomorphism in cyclic homology and thus
  periodic cyclic homology.
\end{proof}

\subsection{Analytic quasi-freeness of
  \texorpdfstring{$\mathcal{L}$}{L}}
\label{sec:L_aqf}

The proof of the excision theorem is completed by the following
theorem:

\begin{theorem}
  \label{the:L_in_extension}
  There is a semi-split, analytically nilpotent extension
  \(\jans E \cap \mathcal{L} \into \mathcal{L} \onto K\)
  and~\(\mathcal{L}\) is analytically quasi-free.
\end{theorem}

This theorem and Theorem~\ref{Goodwillie} imply that \(\HAC(K)\) is
chain homotopy equivalent to the \(X\)\nb-complex
of~\(\mathcal{L}\).  Theorem~\ref{the:L_to_relative} identifies this
with the homology of the homotopy projective limit of the relative
\(X\)\nb-complex \(X(\tans E: \tans Q)\).  And this yields the
excision theorem.  So it only remains to prove
Theorem~\ref{the:L_in_extension}.

The canonical projection \(\tans E \onto E\) restricts to a
semi-split projection \(\mathcal{L} \onto K\).  Its kernel
\(\jans E \cap \mathcal{L}\subseteq \jans E\) is a projective system
of closed subalgebras.  These are complete and torsion-free by
\cite{Meyer-Mukherjee:Bornological_tf}*{Theorem 2.3 and Lemma 4.2};
and subalgebras also clearly inherit the property of being
semi-dagger.  So \(\jans E\cap \mathcal{L}\) is a projective system
of dagger algebras.
Proposition~\ref{pro:extension_nilpotent_mod_pi} implies that it is
again nilpotent mod~\(\dvgen\) because
\(\jans E/(\jans E\cap \mathcal{L})\) is torsion-free.

The proof of Theorem~\ref{the:L_to_relative} already shows
that~\(\mathcal{L}\) is quasi-free.  We need it to be analytically
quasi-free, however.  This is the main difficulty in
Theorem~\ref{the:L_in_extension}.  The proof of this uses the same
ideas as the proof of the corresponding statement for analytic
cyclic homology for bornological algebras over~\(\C\)
in~\cite{Meyer:HLHA}.  First, we define a homomorphism
\(\upsilon\colon L\to \tens L\) for the purely algebraic
version~\(L\) of~\(\mathcal{L}\).  Then we show that this
homomorphism extends uniquely to a homomorphism of pro-algebras
\(\mathcal{L} \to \tans \mathcal{L}\) that is a section for the
canonical projection \(\tans \mathcal{L} \onto \mathcal{L}\).

We need some notation for elements of~\(\tens L\) and a certain
grading on~\(\tens L\).  Elements of~\(\tens L\) are sums of
differential forms \(l_0 \,\Diff l_1 \dotsc \Diff l_{2 n}\) with
\(l_0\in L^+\), \(l_1,\dotsc,l_{2 n} \in L\).  We write \(\ccirc\)
for the Fedosov product in \(\Omega^\even L\) to distinguish it from
the Fedosov product~\(\odot\) in~\(L\) and the resulting usual
multiplication on \(\Omega L\).  Call an element of~\(\tens L\)
\emph{elementary} if it is of the form
\(l_0 \,\Diff l_1 \dotsc \Diff l_{2 n}\) with
\(l_j = e_{j,0}\,\diff e_{j,1} \dotsc \diff e_{j,2 i_j}\) for
\(0 \le j \le 2 n\), and \(e_{j,k} \in K \cup s(Q)\) for all
occurring indices \(j,k\), except that we allow \(l_0 = 1\) and then
put \(i_0 = 0\); here \(e_{j, 2 i_j} \in K\) because \(l_j\in L\).
Any element of~\(\tens L\) is a finite linear combination of such
elementary elements.  The \emph{entries} of an elementary
element~\(\xi\) are the elements \(e_{j,l} \in E\); its
\emph{internal degree} is \(\deg_i(\xi) = \sum_{j=0}^{2 n} 2 i_j\);
its \emph{external degree} is \(\deg_e(\xi) = 6 n\) if \(l_0\in L\)
and \(\deg_e = 6 n - 4\) if \(l_0=1\), and the total degree
\(\deg_t(\xi)\) is the sum of these two degrees; this particular
total degree already appears in the proof of
\cite{Meyer:HLHA}*{Lemma~5.102}.

The definition of~\(\upsilon\) is based on the isomorphism
\(L \cong (\tens E)^+ \otimes K\) in~\eqref{eq:tensEplus_K_is_L}.
The restriction of~\(\upsilon\) to
\(K = (\Omega^0 \tens E \cap L) \subseteq L\) is the obvious
inclusion of~\(K\) into~\(\tens L\).  We extend this map to~\(L\)
using a homomorphism from~\(\tens E\) to the algebra of
\(\dvr\)\nb-module homomorphisms \(\tens L \to \tens L\).  Such a
homomorphism is equivalent to a linear map
\(E \to \Hom (\tens L,\tens L)\), which is, in turn, equivalent to a
\(\dvr\)\nb-bilinear map \(E\times \tens L \to \tens L\), which we
denote as an operation \((e,\xi) \mapsto e\triangleright \xi\) for
\(e\in E\), \(\xi\in \tens L\).  As in~\cite{Meyer:HLHA}, we first
define the map \(\nabla\colon L \to \Omega^1(L)\) by
\(\nabla (s_R(\xi) \odot k \odot l) \defeq s_R(\xi) \odot k \,\Diff
l\) for all \(\xi \in (\tens Q)^+\), \(k\in K\), \(l\in L^+\), with
the understanding that \(\Diff 1=0\) if~\(l\) is the unit element
of~\(L\); this uses the inverse of the
isomorphism~\eqref{eq:tensQplus_K_Lplus_is_L}.  Then we let
\[
  \begin{aligned}
    e \triangleright x_0\,\Diff x_1\dotsc\Diff x_{2n}
    &= e\Fed x_0\,\Diff x_1\dotsc\Diff x_{2n}
    - \Diff \nabla(e\Fed x_0)\,\Diff x_1\dotsc\Diff x_{2n},
    \\
    e \triangleright \Diff x_1\,\Diff x_2\dotsc\Diff x_{2n}
    &= \nabla(e\Fed x_1)\,\Diff x_2\dotsc\Diff x_{2n}.
  \end{aligned}
\]
The curvature of the corresponding map
\(E\to \Hom(\tens L,\tens L)\) acts by the operation
\(\omega_\triangleright(e_1,e_2) \xi \defeq (e_1 \cdot e_2)
\triangleright \xi - e_1 \triangleright (e_2 \triangleright \xi)\).
It is computed in \cite{Meyer:HLHA}*{Equation~(5.91)}:
\begin{align*}
  \omega_\triangleright(e_1,e_2) l_0\,\Diff l_1 \dotsc \Diff l_{2 n}
  &= (\diff e_1 \diff e_2 \odot l_0)\,\Diff l_1 \dotsc \Diff l_{2 n}
  \\&\qquad+
  \nabla \bigl(e_1 \odot \nabla (e_2 \odot l_0)\bigr)
  \,\Diff l_1 \dotsc \Diff l_{2 n}
  \\&\qquad- \Diff \nabla(\diff e_1 \diff e_2 \odot l_0)
  \,\Diff l_1 \dotsc \Diff l_{2 n},\\
  \omega_\triangleright(e_1,e_2) \Diff l_1 \dotsc \Diff l_{2 n}
  &= \nabla ((e_1 \cdot e_2) \odot l_1)\,\Diff l_2 \dotsc \Diff l_{2 n}
  \\&\qquad-
  e_1 \odot \nabla (e_2 \odot l_1)\,\Diff l_2 \dotsc \Diff l_{2 n}
  \\&\qquad+
  \Diff \nabla \bigl(e_1 \odot \nabla (e_2 \odot l_1)\bigr)\,\Diff l_2
    \dotsc \Diff l_{2 n}.
\end{align*}
Finally, we define
\[
  \upsilon(e_0\,\diff e_1 \dotsc \diff e_{2 n}\odot k)
  \defeq e_0 \triangleright \bigl(\omega_\triangleright(e_1,e_2)
  \circ \dotsb \circ \omega_\triangleright(e_{2 n- 1},e_{2 n})\bigr)(k).
\]

\begin{lemma}
  \label{lem:properties_upsilon}
  The map \(\upsilon\colon L \to \tens L\) is an algebra
  homomorphism and \(p\circ \upsilon = \id_L\) for the canonical
  projection \(p\colon \tens L \to L\).

  If \(l\in \Omega^{2 n-1}(E) \,\diff K\subseteq L\) has
  degree~\(2 n\), then \(\upsilon(l)\) is a sum of elementary
  elements of~\(\tens L\) with total degree at least~\(2 n\).

  Let \(M\subseteq E\) be a bounded \(\dvr\)\nb-submodule.  There is
  a bounded subset \(M'\subseteq E\) such that if
  \(e_0\,\diff e_1 \dotsc \diff e_{2 n} \in \Omega^{2 n} M \cap L\),
  then \(\upsilon(e_0\,\diff e_1 \dotsc \diff e_{2 n})\) is a sum of
  elementary elements of~\(\tens L\) with entries in~\(M'\).
\end{lemma}

\begin{proof}
  As shown in \cite{Meyer:HLHA} or in~\cite{Meyer:Excision}, the
  left action~\(\triangleright\) is by left multipliers, that is,
  \(e\triangleright (\xi\ccirc \tau) = (e\triangleright \xi)\ccirc
  \tau\) for all \(e\in E\), \(\xi,\tau \in \tens L\).  And
  \(k\triangleright \xi = k\ccirc \xi\) for all \(k\in K\).  This
  implies that~\(\upsilon\) is a homomorphism.

  A short computation shows that each summand in the formula for
  \(\omega_\triangleright(e_1,e_2)\) increases the total degree
  defined above by at least~\(2\); this is already shown in the
  proof of \cite{Meyer:HLHA}*{Lemma~5.102}.  By induction on~\(n\),
  it follows that~\(\upsilon\) maps \(\Omega^{2 n} L\) into the
  subgroup spanned by elementary elements of~\(\tens L\) with total
  degree at least~\(2 n\).

  Given a bounded subset \(M\subseteq E\), the proof of
  \cite{Meyer:HLHA}*{Lemma~5.92} provides a bounded subset
  \(M'\subseteq E\) such that
  \(\upsilon(e_0\,\diff e_1 \dotsc \diff e_{2 n}\odot k)\) is a sum
  of elementary elements of~\(\tens L\) with entries in~\(M'\).
\end{proof}

The homomorphism~\(\upsilon\) induces an \(\dvf\)\nb-algebra
homomorphism \(L\otimes \dvf \to \tens L\otimes \dvf\).  Recall that
\[
  L_{(m)} \defeq K \oplus \bigoplus_{n=1}^\infty
  \dvgen^{-\floor{n/m}} \Omega^{2 n-1}(E)\,\diff K
\]
for \(m\in\N^*\).  These are \(\dvr\)\nb-subalgebras of
\(L\otimes \dvf\) that satisfy \(L_{(n)} \subseteq L_{(m)}\) if
\(n \ge m\).  Each~\(L_{(m)}\) is equipped with the bornology
cofinally generated by the submodules in~\eqref{eq:bornology_Lm}.

Let \((\tens L)_{(m)} \subseteq \tens L\otimes \dvf\) be the
subgroup generated by \(\dvgen^{-\floor{d/m}} \xi\) for elementary
elements~\(\xi\) of total degree~\(d\).  These are
\(\dvr\)\nb-subalgebras of \(\tens L\otimes \dvf\) that satisfy
\((\tens L)_{(n)} \subseteq (\tens L)_{(m)}\) if \(n \ge m\).  If
\(M\subseteq E\) is a bounded \(\dvr\)\nb-submodule, then let
\(D^\tens_m(M) \subseteq (\tens L)_{(m)}\) be the subgroup generated
by \(\dvgen^{-\floor*{d/m}} \xi\) for elementary elements~\(\xi\) of
total degree~\(d\).  We give~\((\tens L)_{(m)}\) the bornology that
is cofinally generated by these \(\dvr\)\nb-submodules.  This
bornology is the analogue of the bornology in
Corollary~\ref{cor:anborno}.  It is torsion-free and makes the
multiplication in~\((\tens L)_{(m)}\) and the inclusion maps
\((\tens L)_{(n)} \hookrightarrow (\tens L)_{(m)}\) for \(n \ge m\)
bounded.  So we have turned
\(\bigl((\tens L)_{(m)}\bigr)_{m\in\N^*}\) into a projective system
of torsion-free bornological algebras.

The second paragraph in Lemma~\ref{lem:properties_upsilon} says that
the extension \(L\otimes \dvf \to \tens L\otimes \dvf\)
of~\(\upsilon\) maps \(L_{(m)}\) to~\((\tens L)_{(m)}\) for each
\(m\in\N^*\).  And the third paragraph says that this homomorphism
is bounded.  Thus~\(\upsilon\) is a homomorphism of projective
systems of bornological algebras.  By Corollary~\ref{cor:anborno},
\(\mathcal{L}\) is isomorphic to the projective system of the
completions~\(\comb{L_{(m)}}\) for \(m\in\N^*\), with the
bornologies described above.

\begin{lemma}
  \label{lem:cal_L}
  The embedding \(\tens L \hookrightarrow \tans L\) extends to an
  isomorphism of projective systems from the projective system of
  completions \(\comb{(\tens L)_{(m)}}\) for \(m\in\N^*\)
  to~\(\tans \mathcal{L}\).
\end{lemma}

\begin{proof}
  For a bounded \(\dvr\)\nb-submodule \(M\subseteq E\), let
  \(M_K \defeq M\cap K\) and let
  \(\comb{\Omega}_\mathcal{L}^0(M) \defeq M_K\) and
  \(\comb{\Omega}_\mathcal{L}^{2 k}(M) \defeq \comb{\Omega}^{2 k-1}
  (M) \hot M_K\) for \(k>1\).  A proof like that for
  Proposition~\ref{pro:completed_analytic_tensor} shows that the
  completion of~\((\tens L)_{(m)}\) is the union of the products
  \begin{multline}
    \label{eq:TL_strange_bornology}
    \prod_{j\ge0,i_0,\dotsc,i_{2 j} \ge 0}^\infty
    \dvgen^{-\floor{(6 j+ 2 i_0 + \dotsb + 2 i_{2 j})/m}}
    \comb{\Omega}_\mathcal{L}^{2 i_0}(M)
    \hot \comb{\Omega}_\mathcal{L}^{2 i_1}(M)
    \hot \dotsb \hot \comb{\Omega}_\mathcal{L}^{2 i_{2 j}}(M)
    \\\times
    \prod_{j\ge0,i_1,\dotsc,i_{2 j} \ge 0}^\infty
    \dvgen^{-\floor{(6 j - 4 + 2 i_0 + \dotsb + 2 i_{2 j})/m}}
    \comb{\Omega}_\mathcal{L}^{2 i_1}(M)
    \hot \comb{\Omega}_\mathcal{L}^{2 i_2}(M)
    \hot \dotsb \hot \comb{\Omega}_\mathcal{L}^{2 i_{2 j}}(M)
  \end{multline}
  taken over all bounded \(\dvr\)\nb-submodules \(M\subseteq E\);
  elementary tensors in a factor of the first product correspond to
  differential forms \(l_0 \,\Diff l_1\dotsc \Diff l_{2 j}\) with
  \(l_0,\dotsc,l_{2 j}\in \mathcal{L}\) and \(\deg(l_j) = 2 i_j\),
  whereas those for the second product correspond to differential
  forms \(\Diff l_1\dotsc \Diff l_{2 j}\).  The exponent
  of~\(\dvgen\) is the total degree defined above.

  Proposition~\ref{pro:completed_analytic_tensor}
  describes~\(\tans E\).  The pro-subalgebra~\(\mathcal{L}\) is
  described similarly, by also asking for the last entry of all
  differential forms to belong to~\(K\).  Then a second application
  of Proposition~\ref{pro:completed_analytic_tensor}
  describes~\(\tans \mathcal{L}\).  The result is very similar to
  the projective system above.  The only difference is that the
  exponent of~\(\dvgen\) in the bornology is replaced by
  \(h \defeq \floor{j/k} + \sum_{l=0}^{2 j} {}\floor{i_l/m}\) for
  each factor in~\eqref{eq:TL_strange_bornology}.  So it remains to
  prove linear estimates between these two notions of ``degree''.
  In one direction, this is the trivial estimate
  \[
    \floor*{\frac{j}{m}} + \sum_{l=0}^{2 j} {}\floor*{\frac{i_l}{m}}
    \le
    \floor*{\frac{j}{m} + \sum_{l=0}^{2 j} {}\frac{i_l}{m}}
    \le
    \floor*{\frac{1}{m} \biggl( 6j + \sum_{l=0}^{2 j} i_l\biggr)}
  \]
  for \(j\ge0\) and a similar estimate with \(6 j-4 = 4(j-1)+2 j\)
  instead of~\(6j\) for \(j\ge1\).  In the other direction, we
  distinguish two cases.  Let \(i\defeq \sum i_l\).  If
  \(i < 4 j\cdot m\), then \(6 j + 2 i < j\cdot (6+ 8 m)\) and we
  simply estimate
  \[
    \floor*{\frac{j}{m}} + \sum_{l=0}^{2 j} {}\floor*{\frac{i_l}{m}}
    \ge \floor*{\frac{j}{m}}
    \ge \floor*{\frac{6j + 2 i}{(6+ 8 m)\cdot m}}.
  \]
  The other case is \(i\ge 4 j\cdot m\).  Each floor operation
  changes a number by at most~\(1\), and
  \(6 j + 2 i \le \frac{3}{2 m} i + 2 i \le 4 i\).  So
  \[
    \floor*{\frac{j}{m}} + \sum_{l=0}^{2 j} {}\floor*{\frac{i_l}{m}}
    \ge \frac{i}{m} - 2 j
    \ge \frac{i}{2 m}
    \ge \floor*{\frac{6j + 2 i}{8 m}}.\qedhere
  \]
\end{proof}

As a result, \(\upsilon\) defines a pro-algebra homomorphism
\(\mathcal{L} \to \tans \mathcal{L}\).  Then~\(\mathcal{L}\) is
analytically quasi-free.  This ends the proof of the excision
theorem.

\section{Stability with respect to algebras of matrices}
\label{sec:stability}
\numberwithin{equation}{section}

A \emph{matricial pair} consists of two torsion-free bornological
modules \(X\) and~\(Y\) and a surjective linear map
\(\braket{{\cdot}}{{\cdot}}\colon Y\otimes X\to \dvr\).  Any such
map is bounded.  A \emph{homomorphism} from \((X,Y)\) to another
matricial pair \((W,Z)\) is a pair \(f=(f_1,f_2)\) of bounded linear
homomorphisms \(f_1\colon X\to W\), \(f_2\colon Y\to Z\) such that
\(\braket{f_2(y)}{f_1(x)}=\braket{y}{x}\) for all \(x\in X\) and
\(y\in Y\).  An \emph{elementary homotopy} is a pair \(H=(H_1,H_2)\)
of bounded linear maps, where \(H_1\colon X\to W[t]\) and
\(H_2\colon Y\to Z\) or \(H_1\colon X\to W\) and
\(H_2\colon Y\to Z[t]\), such that the following diagram commutes:
\[
  \begin{tikzcd}[column sep=1.5cm]
    Y\otimes X \ar[d, "\braket{}{}"] \ar[r, "H_2\otimes H_1"] &
    Z\otimes W[t] \ar[d, "\braket{}{}\otimes\id"] \\
    \dvr \ar[r, rightarrowtail, "\mathrm{inc}"] & \dvr[t]
  \end{tikzcd}
\]

Let \((X,Y)\) be a matricial pair.  Let \(\mathcal{M}=\mathcal{M}(X,Y)\) be
\(X\otimes Y\) with the product
\[
(x_1\otimes y_1)(x_2\otimes y_2)=\braket{y_1}{x_2}x_1\otimes y_2.
\]
This product is associative and bounded, and it even makes~\(\mathcal{M}\) a
semi-dagger algebra.  The bornological algebra~\(\mathcal{M}\) is also
torsion-free by
\cite{Meyer-Mukherjee:Bornological_tf}*{Proposition~4.12}.
Thus the completion~\(\comb{\mathcal{M}}\) is a dagger algebra and
\(\comb{\mathcal{M}} = \mathcal{M}^\updagger\).

Homomorphisms and homotopies of matricial pairs induce homomorphisms
and homotopies of the corresponding algebras.  Any pair
\((\xi,\eta)\in X\times Y\) with \(\braket{\eta}{\xi}=1\) yields a
bounded algebra homomorphism
\[
  \iota=\iota_{\xi,\eta}\colon \dvr\to \mathcal{M},\qquad
  \iota(1)=\xi\otimes\eta.
\]
We shall also write~\(\iota\) for the composite of the map above
with the completion map
\(\mathcal{M}\to\comb{\mathcal{M}} = \mathcal{M}^\updagger\).
If~\(R\) is a torsion-free bornological algebra, then
\(R\hot \mathcal{M}^\updagger\) is torsion-free by
\cite{Meyer-Mukherjee:Bornological_tf}*{Theorem~4.6 and Propositions
  14.11 and~14.12}.  Define
\begin{equation}
  \label{map:iotaR}
  \iota_R\defeq \id_R\otimes\iota\colon R\to R\hot \mathcal{M}^\updagger.
\end{equation}

\begin{proposition}
  \label{prop:stable}
  Let~\(R\) be a complete, torsion-free bornological algebra.  Then
  the map~\(\iota_R\) induces a chain homotopy equivalence
  \(\HAC(R)\simeq \HAC(R\hot \mathcal{M}^\updagger)\) and an isomorphism
  \(\HA_*(R) \cong \HA_*(R \hot \mathcal{M}^\updagger)\).
\end{proposition}

\begin{proof}
  Corollary~\ref{coro:composmallcurvature} yields a natural
  pro-algebra homomorphism
  \(\tans(R\hot\mathcal{M}^\updagger)\to
  \tans(R)\hot\mathcal{M}^{\updagger}\) covering the identity of
  \(R\hot\mathcal{M}^\updagger\).  And any elementary homotopy
  between matricial pairs \((X,Y)\) and \((W,Z)\) yields an
  elementary dagger homotopy
  \(\mathcal{M}(X,Y)^\updagger\to \mathcal{M}(Z,W)^\updagger\hot
  \dvr[t]^\updagger\).  The \(X\)\nb-complex is invariant under
  dagger homotopies by Proposition~\ref{pro:X-homotopy}.  Taking all
  this into account, the argument of the proof of
  \cite{Meyer:HLHA}*{Theorem~5.65} now applies verbatim and proves
  the proposition.
\end{proof}

Let~\(\Lambda\) be a set.  We now describe increasingly complicated
algebras of matrices indexed by the set~\(\Lambda\).

\begin{example}
  \label{exa:finite_matrices}
  Let~\(\Lambda\) be a set and let~\(\dvr^{(\Lambda)}\) be the
  \(\dvr\)\nb-module of finitely supported functions
  \(\Lambda\to \dvr\).  This is the free module with basis
  \(\setgiven{\chi_\lambda}{\lambda\in \Lambda}\) formed by the
  characteristic functions of the singletons.  The algebra
  \(\mathcal{M}(\dvr^{(\Lambda)},\dvr^{(\Lambda)})\) associated to the
  bilinear form
  \(\braket{\chi_\lambda}{\chi_\mu}=\delta_{\lambda,\mu}\) is just
  the algebra~\(M_\Lambda\) of finitely supported matrices indexed
  by \(\Lambda\times\Lambda\), equipped with the fine bornology.
  The latter algebra is already a dagger algebra.
  Proposition~\ref{prop:stable} implies
  \(\HAC(R) \cong \HAC(M_\Lambda \otimes R)\) for all~\(R\).
\end{example}

\begin{example}
  \label{exa:c0_matrices}
  Define~\(\dvr^{(\Lambda)}\) as in
  Example~\ref{exa:finite_matrices}.  Its \(\dvgen\)\nb-adic
  completion is the Banach module \(c_0(\Lambda)\) with the supremum
  norm.  The bilinear form in Example~\ref{exa:finite_matrices}
  extends to \(c_0(\Lambda)\).  The \(\dvgen\)\nb-adic completion of
  \(\mathcal{M}(c_0(\Lambda),c_0(\Lambda))\) is isomorphic to the Banach
  \(\dvr\)\nb-algebra
  \(M^0_\Lambda \cong c_0(\Lambda\times\Lambda)\) of matrices
  indexed by \(\Lambda\times\Lambda\) with entries going to zero at
  infinity.  The Banach \(\dvr\)\nb-modules above become
  bornological by declaring all subsets to be bounded.  Then the
  completions and tensor products as Banach \(\dvr\)\nb-modules and
  as bornological \(\dvr\)\nb-modules are the same.  Therefore,
  Proposition~\ref{prop:stable} implies
  \(\HAC(R) \cong \HAC(M_\Lambda^0 \hot R)\) for all~\(R\).
\end{example}

\begin{example}
  \label{exa:length-controlled_matrices}
  Let \(\ell\colon \Lambda \to \N\) be a proper function, that is,
  for each \(n\in\N\) the set of \(x\in\Lambda\) with
  \(\ell(x) \le n\) is finite.  Define~\(\dvr^{(\Lambda)}\) as in
  Example~\ref{exa:finite_matrices} and give it the bornology that
  is cofinally generated by the \(\dvr\)\nb-submodules
  \[
    S_m \defeq \sum_{\lambda\in\Lambda}
    \dvgen^{\floor{\ell(\lambda)/m}} \chi_\lambda
  \]
  for \(m\in\N^*\).  The bilinear form in
  Example~\ref{exa:finite_matrices} remains bounded for this
  bornology on~\(\dvr^{(\Lambda)}\).  So
  \(\mathcal{M}(\dvr^{(\Lambda)},\dvr^{(\Lambda)})\) with the tensor
  product bornology from the above bornology is a bornological
  algebra as well.  It is torsion-free and semi-dagger.  So its
  dagger completion is the same as its completion.  We denote it
  by~\(M_\Lambda^\ell\).  It is isomorphic to the algebra of
  infinite matrices \((c_{x,y})_{x,y\in\Lambda}\) for which there is
  \(m\in\N^*\) such that
  \(c_{x,y} \in \dvgen^{\floor{(\ell(x)+\ell(y))/m}}\) for all
  \(x,y\in\Lambda\); this is the same as asking for
  \(\lim {}\abs*{c_{x,y} \dvgen^{-\floor{(\ell(x)+\ell(y))/m}}} = 0\)
  because~\(\ell\) is proper.  It makes no difference to replace the
  exponent of~\(\dvgen\) by \(\floor{\ell(x)/m}+\floor{\ell(y)/m}\)
  or \(\floor{\max \{\ell(x),\ell(y)\}/m}\) because we may
  vary~\(m\).  Proposition~\ref{prop:stable} implies
  \(\HAC(R) \cong \HAC(M_\Lambda^\ell \hot R)\) for all~\(R\).
\end{example}

The following completed matrix algebras will be needed in
Section~\ref{sec:Leavitt}.

\begin{example}
  \label{exa:filtered_length-controlled_matrices}
  Let~\(\Lambda\) be a set with a filtration by a directed
  set~\(I\).  That is, there are subsets
  \(\Lambda_S\subseteq \Lambda\) for \(S\in I\) with
  \(\Lambda_S \subseteq \Lambda_T\) for \(S\le T\) and
  \(\Lambda = \bigcup_{S\in I} \Lambda_S\).  Let
  \(\ell\colon \Lambda \to \N\) be a function whose restriction
  to~\(\Lambda_S\) is proper for each \(S\in I\).  For
  \(S\in\Lambda\), form the matrix algebra~\(M_{\Lambda_S}^\ell\) as
  in Example~\ref{exa:length-controlled_matrices}.  These algebras
  for \(S\in I\) form an inductive system.  Let
  \(\varinjlim M_{\Lambda_S}^\ell\) be its bornological inductive
  limit.  This bornological algebra is also associated to a
  matricial pair, namely, the pair based on
  \(\varinjlim \dvr^{(\Lambda_S)}\), where each
  \(\dvr^{(\Lambda_S)}\) carries the bornology described in
  Example~\ref{exa:length-controlled_matrices}.  Thus
  Proposition~\ref{prop:stable} implies
  \(\HAC(R) \cong \HAC(\varinjlim M_{\Lambda_S}^\ell \hot R)\) for
  all~\(R\).
\end{example}

\section{Leavitt path algebras}
\label{sec:Leavitt}
\numberwithin{equation}{section}

Our next goal is to compute the analytic cyclic homology for tensor
products with Leavitt and Cohn path algebras of directed graphs and
their dagger completions.  A \emph{directed graph}~\(E\) consists of
a set~\(E^0\) of vertices and a set~\(E^1\) of edges together with
source and range maps \(s,r\colon E^1\to E^0\).
A vertex \(v\in E^0\) is \emph{regular} if
\(0< \abs{s^{-1}(\{v\})} <\infty\).  Let \(\reg(E)\subseteq E^0\) be
the subset of regular vertices.  Define
\[
  N_E\colon E^0\times\reg(E)\to \Z,\qquad
  (v,w)\mapsto \delta_{v,w}- \abs{s^{-1}(\{w\})\cap r^{-1}(\{v\})}.
\]
Let \(L(E)\) and \(C(E)\) be the Leavitt and Cohn path algebras
over~\(\dvr\), as defined in
\cite{Abrams-Ara-Siles-Molina:Leavitt_path}*{Definitions 1.2.3
  and~1.2.5}.  We consider them as bornological algebras with the
fine bornology.  The following theorem follows easily from the
results in~\cite{Cortinas-Montero:K_Leavitt} and the formal
properties of analytic cyclic homology:

\begin{theorem}
  \label{the:uncompleted_Leavitt}
  Assume \(\car \dvf = 0\).  Let~\(R\) be a complete bornological
  algebra.  Let~\(E\) be a graph with countably many vertices.  Then
  \begin{alignat*}{2}
    \HAC(R \otimes C(E))
    &\simeq \HAC(R \otimes \dvr^{(E^0)}),&\qquad
    \HAC(C(E))
    &\simeq \dvr^{(E^0)},\\
    &&\HAC(L(E))
    &\simeq \coker(N_E)\oplus\ker(N_E)[1],
  \end{alignat*}
  If~\(E^0\) is finite, then
  \begin{align*}
    \HAC(R \otimes C(E))
    &\simeq \bigoplus_{v\in E^0} \HAC(R),\\
    \HAC(R \otimes L(E))
    &\simeq (\coker(N_E)\oplus\ker(N_E)[1]) \otimes \HAC(R).
  \end{align*}
\end{theorem}

\begin{proof}
  We define a functor~\(H\) from the category of
  \(\dvr\)\nb-algebras to the triangulated category of
  pro-supercomplexes by giving~\(A\) the fine bornology and taking
  \(\HAC(R\otimes A)\).  The functor~\(H\) is homotopy invariant for
  polynomial (and even dagger) homotopies by
  Theorem~\ref{the:homotopy_invariance_X}, stable for algebras of
  finite matrices over any set~\(\Lambda\) by
  Proposition~\ref{prop:stable} applied to
  Example~\ref{exa:finite_matrices}, and exact on semi-split
  extensions by Theorem~\ref{the:excision}.
  Theorem~\ref{the:excision} also implies that~\(\HAC\) is finitely
  additive.  It is not countably additive in general, but
  Corollary~\ref{cor:sum_V_quasi-free} shows that it is countably
  additive on the ground ring~\(\dvr\).  Now
  \cite{Cortinas-Montero:K_Leavitt}*{Theorem~4.2} proves a homotopy
  equivalence
  \[
    \HAC(R\otimes C(E)) \simeq \HAC(R \otimes \dvr^{(E^0)}).
  \]
  If~\(E^0\) is finite, then this is homotopy equivalent to
  \(\HAC(R) \otimes \dvr^{(E^0)} = \bigoplus_{v\in E^0} \HAC(R)\) by
  finite additivity.  And if \(R=\dvr\), then
  Corollary~\ref{cor:sum_V_quasi-free} identifies
  \(\HAC(\dvr^{(E^0)}) \simeq \dvr^{(E^0)}\).

  \cite{Cortinas-Montero:K_Leavitt}*{Proposition~5.2} yields a
  distinguished triangle of pro-supercomplexes
  \[
    \HAC(R \otimes \dvr^{(\reg(E))})
    \xrightarrow{f} \HAC(R \otimes \dvr^{(E^0)})
    \to \HAC(R \otimes L(E))
    \to \HAC(R \otimes \dvr^{(\reg(E))})
  \]
  and partly describes the map~\(f\).  If \(R=\dvr\) and~\(E^0\) is
  countable, then Corollary~\ref{cor:sum_V_quasi-free} identifies
  \(\HAC(\dvr^{(E^0)}) \simeq \dvr^{(E^0)}\) and
  \(\HAC(\dvr^{\reg(E)}) \simeq \dvr^{\reg(E)}\), and the
  information about the map~\(f\) in
  \cite{Cortinas-Montero:K_Leavitt}*{Proposition~5.2} shows that it
  multiplies vectors with the matrix~\(N_E\).  If~\(E^0\) is
  finite, then~\(\HAC\) is \(E^0\)\nb-additive and
  \cite{Cortinas-Montero:K_Leavitt}*{Theorem~5.4} gives a
  distinguished triangle
  \[
    \HAC(R)\otimes \dvf^{\reg(E)} \xrightarrow{\id\otimes N_E}
    \HAC(R)\otimes \dvf^{E^0} \to
    \HAC(R\otimes L(E)) \to \dotsb.
  \]
  Since \(\car(\dvf)=0\), there are invertible matrices \(x,y\) with
  entries in~\(\dvf\) such that \(x N_E y\) is a diagonal matrix
  with only zeros and ones in the diagonal.  We may replace the map
  \(N_E\) or \(\id\otimes N_E\) above by \(\id\otimes (x N_E y)\).
  Then the formulas for \(\HAC(L(E))\) in general and for
  \(\HAC(R \otimes L(E))\) for finite~\(E^0\) follow.
\end{proof}

\begin{corollary}
  \label{coro:funda_uncomplete}
  \(\HAC(R\otimes \dvr[t,t^{-1}])\) is chain homotopy equivalent to
  \(\HAC(R)\oplus \HAC(R)[1]\) and
  \(\HA_*(R\otimes \dvr[t,t^{-1}]) \cong \HA_*(R) \oplus
  \HA_*(R)[1]\).
\end{corollary}

\begin{proof}
  Apply Theorem~\ref{the:uncompleted_Leavitt} to the graph
  consisting of one vertex and one loop.
\end{proof}

The following theorem says that
Theorem~\ref{the:uncompleted_Leavitt} remains true for the dagger
completions \(C(E)^\updagger\) and \(L(E)^\updagger\) of \(C(E)\)
and \(L(E)\):

\begin{theorem}
  \label{the:completed_Leavitt}
  Let~\(R\) be a complete bornological algebra and let~\(E\) be a
  graph.  Then
  \[
    \HAC(R \otimes C(E)) \simeq \HAC(R \hot C(E)^\updagger),
    \qquad
    \HAC(R \otimes L(E)) \simeq \HAC(R \hot L(E)^\updagger).
  \]
\end{theorem}

So the formulas in Theorem~\textup{\ref{the:uncompleted_Leavitt}}
also compute \(\HAC(R \hot C(E)^\updagger)\) and
\(\HAC(R \hot L(E)^\updagger)\) -- assuming~\(E^0\) to be countable
or finite or \(R=\dvr\) for the different cases.

\begin{corollary}[Fundamental Theorem]
  \label{coro:funda}
  \(\HAC(R\hot \dvr[t,t^{-1}]^\updagger)\) is chain homotopy
  equivalent to \(\HAC(R)\oplus \HAC(R)[1]\) and
  \(\HA_*(R\hot \dvr[t,t^{-1}]^\updagger) \cong \HA_*(R) \oplus
  \HA_*(R)[1]\).
\end{corollary}

\begin{proof}
  Combine Theorem~\ref{the:completed_Leavitt} and
  Corollary~\ref{coro:funda_uncomplete}.
\end{proof}

We are going to prove Theorem~\ref{the:completed_Leavitt} by showing
that the proofs in~\cite{Cortinas-Montero:K_Leavitt} continue to
work when we suitably complete all algebras that occur there.  We
must be careful, however, because the dagger completion is
\emph{not} an exact functor.  We first recall some basic facts that
are used in~\cite{Cortinas-Montero:K_Leavitt}.  These will be used
to describe the dagger completions \(C(E)^\updagger\) and
\(L(E)^\updagger\).

By definition, \(L(E)\) has the same generators as~\(C(E)\) and more
relations.  This provides a quotient map
\(p\colon C(E) \onto L(E)\).  Let \(K(E) \subseteq C(E)\) be its
kernel.

\begin{lemma}
  \label{lem:Cohn_Leavitt_extension}
  There is a semi-split extension of \(\dvr\)\nb-algebras
  \[
    K(E) \into C(E) \onto L(E).
  \]
\end{lemma}

\begin{proof}
  Let~\(\mathcal{P}\) be the set of finite paths in~\(E\).  For
  \(v\in \reg(E)\), choose \(e_v\in s^{-1}(\{v\})\).  Let
  \begin{align*}
    \mathcal{B}
    &\defeq \setgiven{\alpha\beta^*}
      {\alpha,\beta\in\mathcal{P},\ r(\alpha)=r(\beta)},\\
    \mathcal{B}'
    &\defeq \mathcal{B}\setminus \setgiven{\alpha e_ve_v^*\beta^*}%
      {v\in\reg(E),\ \alpha,\beta\in\mathcal{P},\ r(\alpha)=r(\beta)=v}.
  \end{align*}
  By \cite{Abrams-Ara-Siles-Molina:Leavitt_path}*{Propositions 1.5.6
    and~1.5.11}, \(\mathcal{B}\) is a basis of~\(C(E)\)
  and~\(\mathcal{B}'\) is a basis of~\(L(E)\).  Let
  \(\sigma\colon L(E)\to C(E)\) be the linear map that sends each
  element of~\(\mathcal{B}'\) to itself.  This is a section for the
  quotient map \(p\colon C(E) \to L(E)\).
\end{proof}

Next we describe~\(K(E)\) as in
\cite{Abrams-Ara-Siles-Molina:Leavitt_path}*{Proposition~1.5.11}.
Let \(v\in \reg(E)\).  Define
\[
  q_v \defeq v - \sum_{s(e)=v} e e^*.
\]
Let \(\mathcal{P}_v\subseteq \mathcal{P}\) be the set of all paths
with \(r(\alpha)=v\).  Let~\(\dvr^{(\mathcal{P}_v)}\) be the free
\(\dvr\)\nb-module on the set~\(\mathcal{P}_v\) and
let~\(\mathcal{M}_{\mathcal{P}_v}\) be the algebra of finite matrices
indexed by~\(\mathcal{P}_v\) as in
Example~\ref{exa:finite_matrices}.  The map
\[
  \bigoplus_{v\in \reg(E)} \mathcal{M}_{\mathcal{P}_v} \to K(E),
  \qquad
  \alpha \otimes \beta \mapsto \alpha q_v \beta^*,
\]
is a \(\dvr\)\nb-algebra isomorphism by
\cite{Abrams-Ara-Siles-Molina:Leavitt_path}*{Proposition~1.5.11}.
Each~\(\mathcal{M}_{\mathcal{P}_v}\) with the fine bornology is a dagger
algebra because it is a union of finite-dimensional subalgebras.
Thus~\(K(E)\) is a dagger algebra as well.  In contrast, \(C(E)\)
and~\(L(E)\) with the fine bornology are not semi-dagger algebras.
And the restriction to~\(K(E)\) of the linear growth bornology
of~\(C(E)\) is \emph{not} just the fine bornology: this is visible
in the special case where~\(C(E)\) is the Toeplitz algebra and
\(L(E) = \dvr[t,t^{-1}]\).

We are going to describe the linear growth bornology on~\(C(E)\).
Let~\(\mathcal{F}\) be the set of all finite subsets
\(S\subseteq E^0\cup E^1\) such that
\[
  e\in S\cap E^1 \text{ and } s(e)\in \reg(E)
  \Rightarrow \{s(e)\}\cup s^{-1}(s(e))\subseteq S.
\]
Let~\(S^{(\infty)}\) for \(S\in\mathcal{F}\) be the set of all paths
that consist only of edges in~\(S\).  Let~\(\abs{\alpha}\) be the
length of a path \(\alpha\in\mathcal{P}\).  For \(n\in\N\), let
\[
  S_n\defeq \setgiven{\alpha \beta^*}
  {\alpha,\beta\in S^{(\infty)},\ \abs{\alpha}+\abs{\beta}\le n}
  \subseteq \mathcal{B}.
\]
This is an increasing filtration on the basis~\(\mathcal{B}\)
of~\(C(E)\).

\begin{lemma}
  \label{lem:linear_growth_Cohn_Leavitt}
  A subset of~\(C(E)\) has linear growth if and only if there are
  \(S\in\mathcal{F}\) and \(m\in\N^*\) such that it is contained in
  the \(\dvr\)\nb-linear span of
  \(\bigcup_{n\in\N} \dvgen^{\floor{n/m}} S_n\).
\end{lemma}

\begin{proof}
  It is easy to see that the \(\dvr\)\nb-linear span of
  \(\bigcup_{n\in\N} \dvgen^{\floor{n/m}} S_n\) in~\(C(E)\) has
  linear growth.  Conversely, we claim that any subset of linear
  growth is contained in one of this form.  Every finite subset of
  \(E^0\cup E^1\) is contained in an element of~$\mathcal{F}$.  It
  follows that, for every finitely generated submodule
  \(M\subseteq C(E)\), there are \(S\in \mathcal{F}\) and \(m\ge 1\)
  such that~\(M\) is contained in the \(\dvr\)\nb-submodule
  generated by~\(S_m\).  Then~\(M^j\) is contained in the
  \(\dvr\)\nb-submodule generated by~\(S_{m j}\) for all
  \(j\in\N^*\).  Thus~\(M^\diamond\) is contained in the
  \(\dvr\)\nb-submodule generated by~\(\dvgen^{j-1} S_{m j}\) for
  all \(j\in\N^*\).  This is the \(\dvr\)\nb-linear span of
  \(\bigcup_{n\in\N^*} \dvgen^{\ceil{n/m}-1} S_n\).  Letting~\(m\)
  vary, we may replace~\(\ceil{n/m}-1\) by~\(\floor{n/m}\).
\end{proof}

Constructing linear growth bornologies commutes with taking
quotients.  So a subset of~\(L(E)\) has linear growth if and only if
it is the image of a subset of linear growth in~\(C(E)\).  Next we
show that the section \(\sigma\colon L(E) \to C(E)\) is bounded for
the linear growth bornologies, and we describe the restriction
to~\(K(E)\) of the linear growth bornology on~\(C(E)\):

\begin{lemma}
  \label{lem:linear_growth_Cohn_kernel}
  Give \(\dvr^{(\mathcal{P}_v)}\subseteq \dvr^{(\mathcal{P})}\) the
  bornology where a subset is bounded if and only if it is contained
  in the linear span of
  \(\setgiven{\dvgen^{\floor{\abs{\alpha}/m}} \alpha}{\alpha \in
    S^{(\infty)}}\) for some \(S\in\mathcal{F}\) and some
  \(m\in\N^*\).  Equip the matrix algebra
  \(\mathcal{M}_{\mathcal{P}_v} = \dvr^{(\mathcal{P}_v \times
    \mathcal{P}_v)}\) with the resulting tensor product bornology
  and the multiplication defined by the obvious bilinear pairing as
  in Section~\textup{\ref{sec:stability}}, and give
  \(\bigoplus_{v\in \reg(E)} \mathcal{M}_{\mathcal{P}_v}\) the
  direct sum bornology.  There is a semi-split extension of
  bornological algebras
  \[
    \begin{tikzcd}
      \displaystyle\bigoplus_{v\in \reg(E)} \mathcal{M}_{\mathcal{P}_v}
      \ar[r,rightarrowtail, "i"]&
      \ling{C(E)} \ar[r, twoheadrightarrow, "p"] &
      \ling{L(E)}.  \ar[l, bend left, dotted, "\sigma"]
    \end{tikzcd}
  \]
\end{lemma}

\begin{proof}
  Let \(S\in\mathcal{F}\).  We claim that~\(\sigma\circ p\) maps the
  linear span of~\(S_n\) into itself.  If
  \(\alpha \beta^* \in \mathcal{B}'\), then
  \(\sigma\circ p(\alpha \beta^*) = \alpha\beta^*\).  If
  \(\alpha \beta^* \notin \mathcal{B}'\), then
  \(\alpha = \alpha_0 e_v\), \(\beta = \beta_0 e_v\) for some
  \(v\in \reg(E)\), \(\alpha_0,\beta_0\in \mathcal{P}_v\).  And then
  \[
    p(\alpha \beta^*)
    = p(\alpha_0 \beta_0^*)
    - \sum_{s(e) = v, e\neq e_v} p(\alpha_0 e e^* \beta_0).
  \]
  Since \(\alpha_0 \beta_0^*\) is shorter than \(\alpha \beta^*\)
  and \(\alpha_0 e e^* \beta_0 \in \mathcal{B}'\) for \(e\in E^1\)
  with \(s(e) = v\) and \(e\neq e_v\), an induction over
  \(\abs{\alpha} + \abs{\beta}\) shows that
  \(\sigma\circ p(\alpha \beta^*)\) is always a \(\dvr\)\nb-linear
  combination of shorter words; in addition, all edges in these
  words are again contained in~\(S\) because \(S\in\mathcal{F}\).
  This proves the claim.  Now
  Lemma~\ref{lem:linear_growth_Cohn_Leavitt} implies
  that~\(\sigma\circ p\) preserves linear growth of subsets.
  Equivalently, \(\sigma\) is a bounded map
  \(\ling{L(E)} \to \ling{C(E)}\).  Then a subset of~\(K(E)\) has
  linear growth in~\(C(E)\) if and only if it is of the form
  \((\id - \sigma\circ p)(M)\) for a \(\dvr\)\nb-submodule
  \(M\subseteq C(E)\) that has linear growth.  The projection
  \(\id - \sigma\circ p\) kills \(\alpha \beta^* \in \mathcal{B}'\).
  Thus we may disregard these generators when we describe the
  restriction to~\(K(E)\) of the linear growth bornology
  on~\(C(E)\).  Instead of applying \(\id - \sigma\circ p\) to the
  remaining basis vectors \(\alpha e_v e_v^* \beta^*\) for
  \(r(\alpha) = r(\beta) = v\in \reg(E)\), we may also apply it to
  \(\alpha e_v e_v^* \beta^* - \alpha \beta^*\)
  because~\(\alpha \beta^*\) is a shorter basis vector that involves
  the same edges.  And
  \begin{align*}
    (\id-\sigma\circ p)(\alpha e_v e_v^* \beta^* - \alpha \beta^*)
    &= \alpha e_v e_v^* \beta^* - \alpha \beta^*
    + \sigma\biggl( \sum_{s(e) = v,\ e\neq e_v} p(\alpha e e^* \beta^*)\biggr)
    \\&= - \alpha \beta^* + \sum_{s(e) = v} \alpha e e^* \beta^*
    = - \alpha q_v \beta^*.
  \end{align*}
  Now Lemma~\ref{lem:linear_growth_Cohn_Leavitt} implies that a
  subset of~\(K(E)\) has linear growth in~\(C(E)\) if and only if
  there are \(S\in\mathcal{F}\) and \(m\in\N^*\) so that it belongs
  to the \(\dvr\)\nb-linear span of
  \(\dvgen^{\floor{n/m}} \alpha q_v \beta^*\) with \(v\in \reg(E)\),
  \(\alpha,\beta \in \mathcal{P}_v \cap S^{(\infty)}\), and
  \(\abs{\alpha}+\abs{\beta} + 2 \le n\).  Under the isomorphism
  \(\bigoplus_{v\in \reg(E)} \mathcal{M}_{\mathcal{P}_v} \cong
  K(E)\), this becomes equal to the bornology on
  \(\bigoplus_{v\in \reg(E)} \mathcal{M}_{\mathcal{P}_v}\) specified
  in the statement of the lemma.
\end{proof}

The semi-split extension in
Lemma~\ref{lem:linear_growth_Cohn_kernel} implies a similar
semi-split extension involving the dagger completions
\(C(E)^\updagger\), \(L(E)^\updagger\) and the completion of
\(\bigoplus_{v\in \reg(E)} \mathcal{M}_{\mathcal{P}_v}\) for the
bornology specified in Lemma~\ref{lem:linear_growth_Cohn_kernel}.

Now Theorem~\ref{the:completed_Leavitt} is proven by showing that
all homomorphisms and quasi-homomorphisms that are used
in~\cite{Cortinas-Montero:K_Leavitt} remain bounded and all
homotopies among them remain dagger homotopies when we give all
algebras that occur the suitable ``linear growth'' bornology,
defined using the lengths of paths to define linear growth.  This is
because all maps in~\cite{Cortinas-Montero:K_Leavitt} are described
by explicit formulas in terms of paths, which change the length only
by finite amounts.  We have put linear growth in quotation marks
because the correct bornologies on the ideals \(K(E)\) and
\(\hat{K}(E)\) in~\cite{Cortinas-Montero:K_Leavitt} are restrictions
of linear growth bornologies on larger algebras as in
Lemma~\ref{lem:linear_growth_Cohn_kernel}.  These bornological
algebras are special cases of
Example~\ref{exa:filtered_length-controlled_matrices}, and so
\(\HAC\) is stable for such matrix algebras.  The bornology
on~\(K(E)\) in Lemma~\ref{lem:linear_growth_Cohn_kernel} actually
deserves to be called a ``linear growth bornology''.  But the
relevant length function is specified by hand and not by the length
of products as for the official linear growth bornology in
Definition~\ref{def:semi-dagger}.

\section{Filtered Noetherian rings and analytic quasi-freeness}
\label{sec:filtered_Noether}

\subsection{Finite-degree connections}
\numberwithin{equation}{subsection}

A complete bornological \(\dvr\)\nb-algebra~\(R\) is quasi-free if
the complete bornological \(R\)\nb-bimodule \(\comb{\Omega}^1(R)\) is
projective.  Equivalently, there is a \emph{connection}
on~\(\comb{\Omega}^1(R)\), that is, a linear map
\(\nabla \colon \comb{\Omega}^1(R) \to \comb{\Omega}^2(R)\) satisfying
\[
  \nabla(a \omega) = a \nabla(\omega)
  \quad\text{and}\quad
  \nabla(\omega a)
  = \nabla(\omega)a + \omega \,\diff a,
\]
for all \(a\in R\) and \(\omega \in \comb{\Omega}^1(R)\) (see
\cite{Cuntz-Quillen:Algebra_extensions}*{Proposition~3.4}).

We are going to prove that~\(R\) is analytically quasi-free if the
growth of such a connection may be controlled in a suitable way.
This uses increasing filtrations.
An (increasing) \emph{filtration} on a \(\dvr\)\nb-module~\(M\) is
an increasing sequence of \(\dvr\)\nb-submodules
\((\Fil_n M)_{n\in\N}\) with \(\bigcup \Fil_n M = M\).  For a
\(\dvr\)\nb-algebra~\(R\), we require, in addition, that
\(\Fil_n R \cdot \Fil_m R \subseteq \Fil_{n+m} R\) for all
\(n,m\in\N\).  And for a module~\(M\) over a
\(\dvr\)\nb-algebra~\(R\) with a fixed
filtration~\((\Fil_n R)_{n\in\N}\), we require, in addition, that
\(\Fil_n R \cdot \Fil_m M \subseteq \Fil_{n+m} M\) for all
\(n,m\in\N\).  Then we speak of a \emph{filtered algebra} and a
\emph{filtered module}, respectively.

\begin{definition}
  A map \(f\colon M \to N\) between filtered \(\dvr\)\nb-modules has
  \emph{finite degree} if there is \(a\in \N\) -- the degree -- such
  that \(f(\Fil_nM) \subseteq \Fil_{n+a}(N)\) for all \(n\in \N\).
  Two filtrations \((\Fil_nM)_n\) and \((\Fil_n'M)_n\) on a filtered
  \(\dvr\)\nb-module~\(M\) are called \emph{shift equivalent} if
  there is \(a\in \N\) such that
  \(\Fil_n M \subseteq \Fil_{n+a}' M\) and
  \(\Fil_n' M \subseteq \Fil_{n+a} M\) for all \(n \in \N\).
\end{definition}

\begin{example}
  \label{filtration:differential_forms}
  Let~\(R\) be a torsion-free bornological \(\dvr\)-algebra.
  Define~\(M^{(j)}\) for a complete bounded submodule \(M\subseteq R\) and
  \(j\ge 0\) as in~\eqref{eq:diamond}.  Put
  \begin{equation}
    \label{fil:FM}
    \Fil^M_r\comb{\Omega}^j R \defeq
    \sum_{i_0+\dotsb+i_j\le r} M^{(i_0)}\,\diff M^{(i_1)} \dotsc \diff M^{(i_j)}
    \oplus \sum_{i_1+\dotsb+i_j\le r} \diff M^{(i_1)}\dotsc \diff M^{(i_j)}
  \end{equation}
  for \(r\in \N\).  This is an increasing filtration on the
  differential \(j\)\nb-forms of the subalgebra
  \(M^{(\infty)}\subseteq R\) generated by~\(M\).
\end{example}

The following lemma relates such filtrations to the linear growth
bornology:

\begin{lemma}
  \label{lem:filt}
  Let~\(R\) be a torsion-free bornological algebra, \(M\subseteq R\)
  a bounded \(\dvr\)\nb-submodule and \(n\ge 0\).  Then
  \[
    \sum_{i\ge 0} \dvgen^i \Fil^M_{i+n} \comb{\Omega}^n R
    \subseteq \comb{\Omega}^n (M^\diamond)
    \subseteq \sum_{i\ge 0} \dvgen^i\Fil^M_{i+n+1}\comb{\Omega}^n R.
  \]
\end{lemma}

\begin{proof}
  We compute
  \begin{align*}
    \comb{\Omega}^n (M^\diamond)
    &= M^\diamond \,\diff(M^\diamond)^n\oplus \diff (M^\diamond)^n\\
    &= \sum_{i\ge 0}\dvgen^i \biggl(
      \sum_{i_0+\dots+i_n=i} M^{(i_0+1)}\,\diff M^{(i_1+1)}\dotsc
      \diff M^{(i_n+1)}
    \\&\qquad \oplus \sum_{i_1+\dots+i_n=i}
    \diff M^{(i_1+1)}\dotsc \diff M^{(i_n+1)}\biggr).\qedhere
  \end{align*}
\end{proof}

\begin{lemma}
  \label{lem:filtb}
  Let \(M\subseteq R\) be a bounded submodule, \(r,b\ge 1\) and
  \(s\ge 0\).  Then
  \[
    \Fil^M_r \comb{\Omega}^sR
    \subseteq \Fil^{M^{(b)}}_{\ceil{r/b}+s}\comb{\Omega}^sR
    \subseteq \Fil^M_{r+b(s+1)}\comb{\Omega}^sR.
  \]
\end{lemma}

\begin{proof}
  Straightforward.
\end{proof}

\begin{lemma}
  \label{lem:unibounded}
  Let \(X\) and~\(Y\) be torsion-free bornological modules.
  Let~\((f_n)\) be a sequence of bounded linear maps \(X \to Y\).
  Assume that for each bounded submodule \(M\subseteq X\) there is a
  bounded submodule \(N\subseteq Y\) and a sequence of nonnegative
  integers~\((a_n)\) with \(\lim a_n= \infty\) and
  \(f_n(M)\subseteq \dvgen^{a_n}N\) for all \(n\in\N\).  Then the
  series \(s(x)\defeq \sum_nf_n(x)\) converges in~\(\comb{Y}\) for
  every \(x\in X\), and the assignment \(x\mapsto s(x)\) is bounded
  and linear.  So it extends to a bounded linear map
  \(s\colon \comb{X}\to\comb{Y}\).
\end{lemma}

\begin{proof}
  Straightforward.
\end{proof}

\begin{definition}
  \label{finite_degree_connection}
  Let~\(R\) be a torsion-free bornological \(\dvr\)\nb-algebra.  A
  connection \(\nabla \colon \comb{\Omega}^1(R) \to \comb{\Omega}^2(R)\) has
  \emph{finite degree} on a bounded submodule \(M \subseteq R\) if
  it has finite degree as a \(\dvr\)\nb-module map with respect to
  the filtrations on \(\comb{\Omega}^1(M^{(\infty)})\) and
  \(\comb{\Omega}^2(M^{(\infty)})\) from
  Example~\ref{filtration:differential_forms}.  A
  connection~\(\nabla\) has \emph{finite degree} on~\(R\) if any
  bounded subset is contained in a bounded submodule of~\(R\) on
  which~\(\nabla\) has finite degree.
\end{definition}

\begin{remark}
  \label{rem:bastam}
  Lemma~\ref{lem:filtb} implies that if~$\nabla$ has finite degree
  on~$M$, then it also has finite degree on~$M^{(b)}$ for all~$b$.
  Then~$\nabla$ is a finite degree connection on $M^{(\infty)}$ with
  the bornology that is cofinally generated by~\(M^{(n)}\) for
  \(n\in\N\).
\end{remark}

\begin{theorem}
  \label{thm:qfdagger}
  Let~\(R\) be a complete, torsion-free bornological algebra.
  If~\(\comb{\Omega}^1(R)\) has a connection of finite degree,
  then~\(R^\updagger\) is analytically quasi-free.
\end{theorem}

\begin{proof}
  We introduce some notation on Hochschild cochains.  If~\(X\) is a
  complete, bornological \(R\)\nb-bimodule and
  \(\psi\colon R^{\hot n}\to X\) is an \(n\)\nb-cochain, write
  \(\delta(\psi)\) for its Hochschild coboundary.  If
  \(\xi\colon R^{\hot m}\to Y\) is another cochain, write
  \(\psi\cup\xi\colon R^{\hot n+m}\to X\hot_R Y\) for the cup
  product.  Let
  \(\nabla \colon \comb{\Omega}^1 R \to \comb{\Omega}^2 R\) be a
  connection of finite degree, and let \(M\subseteq R\) be a bounded
  submodule and \(a\ge 0\) an integer such that~\(\nabla\) has
  degree~\(a\) on~\(M\).  The connection~\(\nabla\) is equivalent to
  a \(1\)\nb-cochain \(\varphi_2\colon R\to \comb{\Omega}^2 R\)
  satisfying \(\delta(\varphi_2) = \diff \cup \diff\), via
  \(\nabla(x_0\,\diff x_1) = x_0 \varphi_2(x_1)\) for
  \(x_0\in R^+\), \(x_1\in R\).  Then~\(\varphi_2\) raises the
  \(M\)\nb-filtration degree by at most~\(a\).  If~\(X\) is a
  filtered \(R\)\nb-bimodule and \(\psi\colon R\hot R\to X\) is a
  \(2\)\nb-cocycle of degree at most~\(b\), then
  \[
    \bar{\psi}\colon \comb{\Omega}^2R\to X,\qquad
    \bar{\psi}(x_0\,\diff x_1\,\diff x_2)
    = x_0\psi(x_1,x_2)
  \]
  is a bimodule homomorphism.  And the \(1\)\nb-cochain
  \[
    \psi'=\bar{\psi}\circ\varphi_2
  \]
  raises filtration degree by at most \(a+b\) and satisfies
  \(\delta(\psi')=\psi\).  For \(n\ge 1\), inductively define a
  \(2\)\nb-cocycle and a \(1\)\nb-cochain with values in
  \(\comb{\Omega}^{2(n+1)}R\) as follows:
  \begin{align*}
    \psi_{2(n+1)}
    &\defeq \sum_{j=0}^n \diff \varphi_{2 j}\cup \diff \varphi_{2(n-j)}
      - \sum_{j=1}^n \varphi_{2 j}\cup \varphi_{2(n+1-j)},\\
    \varphi_{2(n+1)}&\defeq \psi_{2(n+1)}'.
  \end{align*}
  Put \(\varphi_0=\id\colon R\to R\).  To see that the maps
  \(\psi_{2n}\) are cocycles, one proves first that
  \[
    \delta(\diff \varphi_{2n})
    = - \sum_{j=0}^{n} \diff (\varphi_{2 j} \cup \varphi_{2(n-j)}).
  \]
  Then a long but straightforward calculation using the Leibniz rule
  for both \(\diff\) and~\(\delta\) shows by induction that
  \(\delta(\psi_{2n})=0\) (see
  \cite{Cortinas:Derived_Cuntz-Quillen}*{Theorem~2.1}).  By
  construction, the bounded linear map
  \(\varphi_{\le 2 n} \defeq \sum_{i=0}^n \varphi_{2 i}\) is a
  section of the canonical projection \(\tens R\to R\), and its
  curvature vanishes modulo \(\jens R^{n+1}\).  So it defines a
  bounded algebra homomorphism \(R \to \tens R/\jens R^{n+1}\).
  Hence the infinite series \(\sum_{i=0}^\infty \varphi_{2 i}\) is
  an algebra homomorphism into the projective limit.  It suffices to
  show that, for each~\(m\), the series
  \(\sum_{i=0}^\infty \varphi_{2 i}\) defines a bounded linear
  homomorphism
  \(\ling{R}\to (\tub{\tens \ling{R}}{\jens \ling{R}^m},
  \tub{\jens \ling{R}}{\jens \ling{R}^m})^\updagger\).

  One checks by induction on~\(n\) that
  \(\varphi_{2n}(M^{(i)})\subseteq
  \Fil^M_{i+(2n-1)a}\comb{\Omega}^{2n}R\).  Hence
  \begin{equation}
    \label{bded:phi}
    \varphi_{2n}(M^\diamond)
    \subseteq \sum_{i=0}^\infty \dvgen^i
    \Fil^M_{i+(2n-1)a+1}\comb{\Omega}^{2 n} R.
  \end{equation}
  Next let \(m\ge 1\) and choose an integer
  \(c> \max\{1, 2am \}\).  Then
  \begin{align}
    \label{ineq:phi}
    i+\floor*{\frac{n}{m}}-\ceil*{\frac{i+(2n-1)a+1}{c}}
    &\ge (1-1/c) i \geq 0
  \end{align}
  for all \(i\ge 0\) and sufficiently large~\(n\).  Then
  \(i \ge \ceil*{\frac{i+(2n-1)a+1}{c}} - \floor*{\frac{n}{m}}\).
  Set \(D(i,n,c) \defeq \ceil*{\frac{i+(2n-1)a+1}{c}}\).  Equations
  \eqref{bded:phi} and~\eqref{ineq:phi} and Lemmas \ref{lem:filt}
  and~\ref{lem:filtb} imply
  \begin{multline*}
    \varphi_{2n}(M^\diamond)
    \subseteq \sum_{i\ge 0} \dvgen^i \Fil_{D(i,n,c) + 2n}^{M^{(c)}}
    \comb{\Omega}^{2n}(R)
    \\ \subseteq \dvgen^{-\floor{\frac{n}{m}}}
    \sum_{i\ge 0} \dvgen^{D(i,n,c)} \Fil_{D(i,n,c) + 2n}^{M^{(c)}} \comb{\Omega}^{2n}(R)
    \subseteq \dvgen^{-\floor{\frac{n}{m}}} \comb{\Omega}^{2n}({(M^{(c)})}^\diamond)
  \end{multline*}
  By Proposition~\ref{pro:completed_analytic_tensor}, the subset of
  infinite series \(\sum_{n=0}^\infty \varphi_{2n}(M^\diamond)\) is
  bounded in
  \((\tub{\tens \ling{R}}{\jens \ling{R}^m},\tub{\jens
    \ling{R}}{\jens \ling{R}^m})^\updagger\).  So
  \(\sum_{n=0}^\infty \varphi_{2n}\) defines a bounded homomorphism
  \[
    R\to (\tub{\tens \ling{R}}{\jens \ling{R}^m},
    \tub{\jens \ling{R}}{\jens \ling{R}^m})^\updagger
  \]
  for each \(m\ge 1\); this completes the proof.  
 \end{proof}

\begin{corollary}
  \label{coro:qfdagger}
  Let~\(R\) be as in Theorem~\textup{\ref{thm:qfdagger}}.  Then the
  natural map \(\HAC(R^\updagger)\to X(R^\updagger \otimes \dvf)\)
  is a chain homotopy equivalence and \(\HA_*(R)\) is isomorphic to
  the homology of \(X(R^\updagger \otimes \dvf)\).
\end{corollary}

\begin{proof}
  Immediate from Theorem~\ref{thm:qfdagger} and
  Corollary~\ref{cor:analytically_quasi_free_computation}.
\end{proof}

\subsection{Filtered Noetherian rings and smooth algebras}

We now show that some quasi-free algebras have a connection of
finite degree.  In particular, this includes smooth, commutative
finitely generated \(\dvr\)-algebras.  For the remainder of this
section, let~\(R\) be a finitely generated
\(\dvr\)\nb-algebra, equipped with the fine bornology.  Let
\(S \subseteq R\) be a finite generating subset and
let~\(S^{\le n}\) be the set of all products of elements of~\(S\) of
length at most~\(n\).  As above, let \(\Fil_n R \subseteq R\) be the
\(\dvr\)\nb-submodule generated by~\(S^{\le n}\).  By convention,
\(S^{\le 0} = \{1\}\) and \(\Fil_0 R = \dvr\cdot 1\).  This is an
increasing filtration on~\(R\).  It induces filtrations on the
bimodules~\(\Omega^l(R)\) as in
Example~\ref{filtration:differential_forms}.  More concretely,
\(\Fil_n(\Omega^l(R))\) is the \(\dvr\)\nb-submodule
of~\(\Omega^l(R)\) generated by \(x_0 \,\diff x_1 \dotsc \diff x_l\)
with \(x_0 \in \Fil_{n_0}(R)\) or \(x_0=1\) and \(n_0=0\), and
\(x_i \in \Fil_{n_i}(R)\) for \(i = 1, \dotsc, l\), and
\(n_0 + \dotsb + n_l \le n\).  By construction, the
\(\dvr\)\nb-module \(\Fil_n R \cdot \Fil_m R\) that is generated by
products \(x\cdot y\) with \(x\in \Fil_n R\), \(y\in \Fil_m R\) is
equal to~\(\Fil_{n+m} R\) for all \(n,m\in\N\).  This is more than
what is required for a filtered algebra, and the extra information
is crucial for the filtration to generate the linear growth
bornology.





Let~\(M\) be an \(R\)\nb-module with a finite generating set
\(S_M\subseteq M\).  Then we define a filtration on~\(M\), called
the \emph{canonical filtration}, by letting~\(\Fil_n M\) be the
\(\dvr\)\nb-submodule generated by \(a\cdot x\) with
\(a\in\Fil_n R\) and \(x\in S_M\).  This satisfies
\(\Fil_m R \cdot \Fil_n M \subseteq \Fil_{n+m}M\) because
\(\Fil_mR \cdot \Fil_nR \subseteq \Fil_{n+m}R\).  The following
proposition characterises canonical filtrations by a universal
property:

\begin{proposition}
  \label{pro:canonical_filtration}
  Let~\(R\) be a filtered \(\dvr\)\nb-algebra and let~\(M\) be a
  finitely generated \(R\)\nb-module.  Equip~\(M\) with the
  filtration described above.  Then any \(R\)\nb-module map
  from~\(M\) to a filtered \(R\)\nb-module~\(Y\) is of finite
  degree.  The canonical filtrations for two different finite
  generating sets of~\(M\) are shift equivalent.
\end{proposition}

\begin{proof}
  Let \(\{m_1,\dotsc,m_n\}\) be a finite generating set for~\(M\) as
  an \(R\)\nb-module.  Let \(h\colon M \to Y\) be an \(R\)\nb-module
  homomorphism into a filtered \(R\)\nb-module~\(Y\).  Since
  \(Y= \bigcup \Fil_l Y\), there is an \(l \in \N\) with
  \(h(m_i) \in \Fil_l Y\) for all \(i=1,\dotsc, m\).  Then
  \(h(a\cdot m_i) \in \Fil_{n+l} R\) for \(a\in \Fil_n R\).  Hence
  \(h(\Fil_n M) \subseteq \Fil_{n+l} Y\) for all \(n\in\N\).  That
  is, \(h\) has finite degree.  In particular, if we equip~\(M\)
  with another filtration~\((\Fil'_n M)_{n\in\N}\), then the
  identity map has finite degree, that is, there is \(l\in\N\) with
  \(\Fil_n M \subseteq \Fil_{n+l}' M\) for all \(n\in\N\).  If the
  other filtration comes from another finite generating set, then we
  may reverse the roles and also get \(l'\in\N\) with inclusions
  \(\Fil_n' M \subseteq \Fil_{n+l'} M\) for all \(n\in\N\).
\end{proof}

\begin{definition}
  \label{def:filtered_Noetherian}
  A filtered \(\dvr\)\nb-algebra~\(R\) is called (left)
  \emph{filtered Noetherian} if every left ideal~\(I\) is finitely
  generated and the filtration \((\Fil_n R \cap I)_{n\in\N}\) is
  shift equivalent to the canonical filtration of
  Proposition~\ref{pro:canonical_filtration} from a finite
  generating set.  In other words, there are finitely many
  \(x_1, \dotsc, x_n \in I\) and \(l \in \N\) such that for all
  \(m\in \N\) and \(y\in \Fil_m R \cap I\), there are
  \(a_i \in \Fil_{m+l} R\) with \(y = \sum_{i=1}^n a_i x_i\).
\end{definition}

\begin{lemma}
  \label{filtered_Noetherian_finite_connection}
  Let~\(R\) be a finitely generated, quasi-free \(\dvr\)\nb-algebra.
  Assume that \(R^+\otimes (R^+)^\op\) is filtered Noetherian.
  Then~\(\Omega^1(R)\) has a connection of finite degree.
\end{lemma}

\begin{proof}
  Since~\(R\) is quasi-free, the left multiplication map
  \(R^+ \otimes \Omega^1(R) \onto \Omega^1(R)\) splits by an
  \(R\)\nb-bimodule homomorphism
  \(s\colon \Omega^1(R) \to R^+ \otimes \Omega^1(R)\).  By
  definition, \(\Omega^1(R)\) is a left ideal in
  \(R^+ \otimes (R^+)^\op\).  By assumption, it is finitely
  generated as such, and the filtration on \(R^+ \otimes (R^+)^\op\)
  restricted to \(\Omega^1(R)\) is the canonical filtration
  on~\(\Omega^1(R)\) as a module over \(R^+ \otimes (R^+)^\op\).
  Now Proposition~\ref{pro:canonical_filtration} shows that the
  section~\(s\) above has finite degree.  The section~\(s\) yields a
  connection \(\nabla \colon \Omega^1(R) \to \Omega^2(R)\), which is
  defined by \(\nabla(\omega) = 1 \otimes \omega - s(\omega)\).  It follows
  that~\(\nabla\) has finite degree.  
\end{proof}

Our next goal is to show that a commutative, finitely generated
\(\dvr\)\nb-algebra with the filtration coming from a finite
generating set is filtered Noetherian.  First consider the
polynomial ring in \(n\)~variables.  The filtration defined by the
obvious generating set is the total degree filtration, where
\(\Fil_m(\dvr[x_1,\dotsc, x_n])\) is the \(\dvr\)\nb-submodule
generated by the monomials of total degree at most~\(m\), that is,
terms of the form
\(x^\alpha = x_1^{\alpha_1} x_2^{\alpha_2} \dotsm x_n^{\alpha_n}\)
with \(\abs{\alpha} \defeq \sum_{i=1}^n \alpha_i \le m\).

\begin{theorem}
  \label{the:polynomial_filtered_Noetherian}
  The polynomial ring \(R = \dvr[x_1,\dotsc, x_n]\) with the total
  degree filtration is filtered Noetherian.
\end{theorem}

\begin{proof}
  Let~\(I\) be any ideal in~\(R\).  Since~\(R\) is Noetherian, \(I\)
  is finitely generated.  Since~\(\dvr\) is a principal ideal
  domain, \(I\) has a finite, strong Gröbner basis with respect to
  any term order on the monomials~\(x^\alpha\) (see
  \cite{Adams-Loustaunau:Introduction_Grobner}*{Theorem~4.5.9}).  We
  use the degree lexicographic order (see
  \cite{Adams-Loustaunau:Introduction_Grobner}*{Definition~1.4.3});
  the only property we need is that \(\abs{\alpha}<\abs{\beta}\)
  implies \(x^\alpha\prec x^\beta\).  The chosen order on monomials
  defines the \emph{leading term} \(\mathrm{lt}(f)\) of a
  polynomial~\(f\).  Let \(G = \{f_1,\dotsc, f_N\}\) be a strong
  Gröbner basis for~\(I\).  By
  \cite{Adams-Loustaunau:Introduction_Grobner}*{Theorem~4.1.12}, any
  \(g\in I\) can be written as \(g = \sum_{j=1}^M c_j t_j f_{i_j}\),
  where \(M\in \N\), \(c_j \in \dvr\), \(t_j\) is a monomial
  in~\(R\), \(i_j \in \{1,\dotsc,N\}\), and
  \(\mathrm{lt}(t_j f_{i_j}) \prec \mathrm{lt}(g)\) for each~\(j\).
  So the total degree of \(t_j f_{i_j}\) is at most the total
  degree of~\(g\) for each \(j=1,\dotsc, M\), and this remains so
  for the total degree of~\(t_j\).  Combining the monomials~\(t_j\)
  with the same~\(i_j\), we write any element \(g\in I\) of total
  degree at most~\(m\) in the form \(\sum_{i=1}^N p_j f_j\) with
  \(p_j \in \Fil_m R\).
\end{proof}

\begin{proposition}
  \label{pro:quotients_filtered_Noetherian}
  A quotient of a filtered Noetherian \(\dvr\)\nb-algebra with the
  induced filtration is again filtered Noetherian.
\end{proposition}

\begin{proof}
  Let~\(R\) be a filtered Noetherian \(\dvr\)\nb-algebra and
  let~\(I\) be an ideal.  Any ideal in the quotient ring~\(R/I\) is
  of the form~\(J/I\) for a unique ideal~\(J\) in~\(R\)
  containing~\(I\).  Let \(x_1, \dotsc, x_n \in J\) and \(l \in \N\)
  be such that for all \(m\in \N\) and \(y\in \Fil_m R \cap I\),
  there are \(a_i \in \Fil_{m+l} R\) with
  \(y = \sum_{i=1}^n a_i x_i\).  Then the images of
  \(x_1,\dotsc,x_n\) in~\(J/I\) and the same~\(l\) will clearly work
  for the ideal~\(J/I\) in the quotient~\(R/I\).
\end{proof}

\begin{corollary}
  \label{finite_type_filtered_Noetherian}
  Any finitely generated, commutative \(\dvr\)\nb-algebra is
  filtered Noetherian.
\end{corollary}

\begin{proof}
  Let~\(A\) be a finitely generated, commutative
  \(\dvr\)\nb-algebra.  Let~\(S\) be any finite generating set.
  Turn it into a surjective homomorphism from the polynomial algebra
  \(R = \dvr[x_1,\dotsc,x_n]\) onto~\(A\).  This identifies
  \(A \cong R/I\) for an ideal~\(I\) in~\(R\).  The filtration
  on~\(A\) defined by~\(S\) is equal to the filtration on the
  quotient~\(R/I\) defined by the degree filtration on~\(R\).  Now
  the claim follows from
  Theorem~\ref{the:polynomial_filtered_Noetherian} and
  Proposition~\ref{pro:quotients_filtered_Noetherian}.
\end{proof}

\begin{proposition}
  \label{pro:smooth_one-dim_quasi-free}
  Let~\(R\) be a smooth, finitely generated commutative
  \(\dvr\)-algebra of relative dimension~\(1\).  Then~\(R\) admits a
  connection of finite degree.
\end{proposition}

\begin{proof}
  The assumptions on~\(R\) imply that~\(\Omega^1(R)\) a projective,
  finitely generated \(R\)\nb-bimodule.  Furthermore, by
  Corollary~\ref{finite_type_filtered_Noetherian}, \(R\) is filtered
  Noetherian.  The result now follows from
  Lemma~\ref{filtered_Noetherian_finite_connection}.
\end{proof}

\begin{remark}
  \label{rem:Monsky-Washnitzer_et_al}
  In their seminal article~\cite{Monsky-Washnitzer:Formal}, Paul
  Monsky and Gerard Washnitzer introduced the so-called
  \textit{Monsky--Washnitzer cohomology} $H_{\MW}^*(A)$ for a smooth
  unital $\resf$\nb-algebra~$A$ that has a ``very smooth'' lift.
  This is a presentation $A=S/\pi S$ where~$S$ is dagger complete
  and very smooth
  (\cite{Monsky-Washnitzer:Formal}*{Definition~2.5}); by definition,
  $H_{\MW}^*(A)=H_\mathrm{dR}(S\otimes \dvf)$ is the de Rham
  cohomology of $S\otimes\dvf$.  As in the current article, Monsky
  and Washnitzer assumed that $\car(\dvf)=0$ but made no assumption
  about the characteristic of~$\resf$.  The very smooth liftability
  assumption in~\cite{Monsky-Washnitzer:Formal} was crucial for
  their proof of the functoriality of $H_{\MW}^*$.  Later on, Marius
  van der Put~\cite{vdp} managed to remove that assumption; for any
  smooth commutative unital $\resf$\nb-algebra~$A$ of finite type,
  he defines $H_{\MW}^*(A)$ as the de Rham cohomology of the dagger
  completion of any smooth $\dvr$\nb-algebra~$R$ with $R/\pi R=A$.
  The existence of such a lift follows from a theorem of Ren\'ee
  Elkik~\cite{Elkik:Solutions}; van der Put proves functoriality of
  $H_{\MW}^*$ using Artin approximation.  However, in his paper he
  assumes that~$\resf$ is finite.  More recently, under very general
  assumptions (in particular, for~$\resf$ of arbitrary
  characteristic) Alberto Arabia~\cite{Arabia:Relevements} proved that every
  smooth $\resf$\nb-algebra admits a very smooth lift, and extended
  the original definition of Monsky and Washnitzer.  In a parallel
  development, Pierre Berthelot introduced rigid cohomology
  $H^*_{\rig}(X)$ of general schemes~$X$ over a field~$\resf$ with
  $\car(\resf)>0$, which for smooth affine $X=\spec(A)$ agrees with
  $H_{\MW}^*(A)$.  With no assumptions on $\car(\resf)$, Elmar
  Grosse-Kl\"onne~\cite{Grosse-Kloenne:De_Rham} introduced the de Rham
  cohomology of dagger spaces over~$\dvr$, and he related it to
  rigid cohomology in the case when $\car(\resf)>0$.
\end{remark}

The following is one of the main applications of our theory:

\begin{theorem}
  \label{the:non-singular_1-dim}
  Let~\(X\) be a smooth affine variety over the residue
  field~\(\resf\) of dimension~\(1\) and let \(A=\mathcal{O}(X)\) be
  its algebra of polynomial functions.  Let~\(R\) be a smooth,
  commutative algebra with \(R/\dvgen R \cong A\).  Equip~\(R\) with
  the fine bornology and let~\(R^\updagger\) be its dagger
  completion.  If \(*=0,1\), then \(\HA_*(R^\updagger)\) is
  naturally isomorphic to the de Rham cohomology of~\(R^\updagger\).
  This is isomorphic to the Monsky--Washnitzer cohomology of~$A$,
  which, if $\car(\resf)>0$, agrees with the rigid cohomology
  \(H^*_\rig(A,\dvf)\) of~\(X\).
\end{theorem}

\begin{proof}
  Since~\(R\) is of finite type over~\(\dvr\), it is Noetherian.  We
  first recall a basic result from commutative algebra:

  \begin{lemma}
    \label{reduction_mod_p:projective}
    Let~\(M\) be a torsion-free, finitely generated \(R\)\nb-module.
    Then~\(M\) is a projective \(R\)\nb-module if and only if
    \(M/\dvgen M\) is a projective \(A\)\nb-module.
  \end{lemma}

  \begin{proof}
    Since the rings \(A\) and~\(R\) are Noetherian, finitely
    generated modules over them are flat if and only if they are
    projective.  Since~\(M\) is torsion-free by hypothesis,
    \cite{Monsky-Washnitzer:Formal}*{Lemma 2.1} shows that~\(M\) is
    projective as an \(R\)\nb-module if and only if it is flat, if
    and only if \(M/\dvgen M\) is flat over~\(A\), if and only
    if~\(M/\dvgen M\) is projective.
  \end{proof}

  In our context, Lemma~\ref{reduction_mod_p:projective} implies
  that~\(R\) is a smooth, finitely generated, commutative
  \(\dvr\)\nb-algebra of relative dimension~\(1\).  By
  Proposition~\ref{pro:smooth_one-dim_quasi-free}, \(R\) is
  quasi-free.  Equipping~\(R\) with the fine bornology, we are in
  the situation of Theorem~\ref{thm:qfdagger}.  Then
  Corollary~\ref{coro:qfdagger} and
  \cite{Cortinas-Cuntz-Meyer-Tamme:Nonarchimedean}*{Theorem~5.5}
  imply the desired isomorphism.
  Remark~\ref{rem:Monsky-Washnitzer_et_al} discusses the generality
  in which different cohomology theories over~\(\resf\) are defined
  and equivalent to the de Rham cohomology of~\(R^\updagger\).
\end{proof}

\end{document}